\documentstyle[amssymb,thmsk,11pt,a4,sw20lart]{article}

\input tcilatex

\makeatletter

\def\diagram{\leftwidth=\z@ \rightwidth=\z@ \topheight=\z@
\botheight=\z@ \setbox\@picbox\hbox\bgroup}

\def\enddiagram{\egroup\wd\@picbox\rightwidth\unitlength
\ht\@picbox\topheight\unitlength \dp\@picbox\botheight\unitlength
\hskip\leftwidth\unitlength\box\@picbox}

\def\bfig{\begin{diagram}}
\def\efig{\end{diagram}}
\newcount\wideness \newcount\leftwidth \newcount\rightwidth
\newcount\highness \newcount\topheight \newcount\botheight

\def\ratchet#1#2{\ifnum#1<#2 \global #1=#2 \fi}

\def\putbox(#1,#2)#3{%
\horsize{\wideness}{#3} \divide\wideness by 2
{\advance\wideness by #1 \ratchet{\rightwidth}{\wideness}}
{\advance\wideness by -#1 \ratchet{\leftwidth}{\wideness}}
\vertsize{\highness}{#3} \divide\highness by 2
{\advance\highness by #2 \ratchet{\topheight}{\highness}}
{\advance\highness by -#2 \ratchet{\botheight}{\highness}}
\put(#1,#2){\makebox(0,0){$#3$}}}

\def\putlbox(#1,#2)#3{%
\horsize{\wideness}{#3}
{\advance\wideness by #1 \ratchet{\rightwidth}{\wideness}}
{\ratchet{\leftwidth}{-#1}}
\vertsize{\highness}{#3} \divide\highness by 2
{\advance\highness by #2 \ratchet{\topheight}{\highness}}
{\advance\highness by -#2 \ratchet{\botheight}{\highness}}
\put(#1,#2){\makebox(0,0)[l]{$#3$}}}

\def\putrbox(#1,#2)#3{%
\horsize{\wideness}{#3}
{\ratchet{\rightwidth}{#1}}
{\advance\wideness by -#1 \ratchet{\leftwidth}{\wideness}}
\vertsize{\highness}{#3} \divide\highness by 2
{\advance\highness by #2 \ratchet{\topheight}{\highness}}
{\advance\highness by -#2 \ratchet{\botheight}{\highness}}
\put(#1,#2){\makebox(0,0)[r]{$#3$}}}

\def\adjust[#1]{} 

\newcount \coefa
\newcount \coefb
\newcount \coefc
\newcount\tempcounta
\newcount\tempcountb
\newcount\tempcountc
\newcount\tempcountd
\newcount\xext
\newcount\yext
\newcount\xoff
\newcount\yoff
\newcount\gap%
\newcount\arrowtypea
\newcount\arrowtypeb
\newcount\arrowtypec
\newcount\arrowtyped
\newcount\arrowtypee
\newcount\height
\newcount\width
\newcount\xpos
\newcount\ypos
\newcount\run
\newcount\rise
\newcount\arrowlength
\newcount\halflength
\newcount\arrowtype
\newdimen\tempdimen
\newdimen\xlen
\newdimen\ylen
\newsavebox{\tempboxa}%
\newsavebox{\tempboxb}%
\newsavebox{\tempboxc}%

\newdimen\w@dth

\def\setw@dth#1#2{\setbox\z@\hbox{$#1$}\w@dth=\wd\z@
\setbox\@ne\hbox{$#2$}\ifnum\w@dth<\wd\@ne \w@dth=\wd\@ne \fi
\advance\w@dth by 1.2em}

\def\t@^#1_#2{\def\n@one{#1}\def\n@two{#2}\mathrel{\setw@dth{#1}{#2}
\mathop{\hbox to \w@dth{\rightarrowfill}}\limits
\ifx\n@one\empty\else ^{\box\z@}\fi
\ifx\n@two\empty\else _{\box\@ne}\fi}}
\def\t@@^#1{\@ifnextchar_ {\t@^{#1}}{\t@^{#1}_{}}}
\def\to{\@ifnextchar^ {\t@@}{\t@@^{}}}

\def\t@left^#1_#2{\def\n@one{#1}\def\n@two{#2}\mathrel{\setw@dth{#1}{#2}
\mathop{\hbox to \w@dth{\leftarrowfill}}\limits
\ifx\n@one\empty\else ^{\box\z@}\fi
\ifx\n@two\empty\else _{\box\@ne}\fi}}
\def\t@@left^#1{\@ifnextchar_ {\t@left^{#1}}{\t@left^{#1}_{}}}
\def\toleft{\@ifnextchar^ {\t@@left}{\t@@left^{}}}

\def\two@^#1_#2{\def\n@one{#1}\def\n@two{#2}\mathrel{\setw@dth{#1}{#2}
\mathop{\vcenter{\hbox to \w@dth{\rightarrowfill}\kern-1.7ex
                 \hbox to \w@dth{\rightarrowfill}}%
       }\limits
\ifx\n@one\empty\else ^{\box\z@}\fi
\ifx\n@two\empty\else _{\box\@ne}\fi}}
\def\tw@@^#1{\@ifnextchar_ {\two@^{#1}}{\two@^{#1}_{}}}
\def\two{\@ifnextchar^ {\tw@@}{\tw@@^{}}}

\def\tofr@^#1_#2{\def\n@one{#1}\def\n@two{#2}\mathrel{\setw@dth{#1}{#2}
\mathop{\vcenter{\hbox to \w@dth{\rightarrowfill}\kern-1.7ex
                 \hbox to \w@dth{\leftarrowfill}}%
       }\limits
\ifx\n@one\empty\else ^{\box\z@}\fi
\ifx\n@two\empty\else _{\box\@ne}\fi}}
\def\t@fr@^#1{\@ifnextchar_ {\tofr@^{#1}}{\tofr@^{#1}_{}}}
\def\tofro{\@ifnextchar^ {\t@fr@}{\t@fr@^{}}}

\def\mon{\mathop{\m@th\hbox to
      14.6\P@{\lasyb\char'51\hskip-2.1\P@$\arrext$\hss
$\mathord\rightarrow$}}\limits} 
\def\leftmono{\mathrel{\m@th\hbox to
14.6\P@{$\mathord\leftarrow$\hss$\arrext$\hskip-2.1\P@\lasyb\char'50%
}}\limits} 
\mathchardef\arrext="0200       

\def\settypes(#1,#2,#3){\arrowtypea#1 \arrowtypeb#2 \arrowtypec#3}
\def\settoheight#1#2{\setbox\@tempboxa\hbox{#2}#1\ht\@tempboxa\relax}%
\def\settodepth#1#2{\setbox\@tempboxa\hbox{#2}#1\dp\@tempboxa\relax}%
\def\settokens[#1`#2`#3`#4]{%
     \def\tokena{#1}\def\tokenb{#2}\def\tokenc{#3}\def\tokend{#4}}
\def\setsqparms[#1`#2`#3`#4;#5`#6]{%
\arrowtypea #1
\arrowtypeb #2
\arrowtypec #3
\arrowtyped #4
\width #5
\height #6
}
\def\setpos(#1,#2){\xpos=#1 \ypos#2}

\def\settriparms[#1`#2`#3;#4]{\settripairparms[#1`#2`#3`1`1;#4]}%

\def\settripairparms[#1`#2`#3`#4`#5;#6]{%
\arrowtypea #1
\arrowtypeb #2
\arrowtypec #3
\arrowtyped #4
\arrowtypee #5
\width #6
\height #6
}

\def\resetparms{\settripairparms[1`1`1`1`1;500]\width 500}

\resetparms

\def\mvector(#1,#2)#3{
\put(0,0){\vector(#1,#2){#3}}%
\put(0,0){\vector(#1,#2){26}}%
}
\def\evector(#1,#2)#3{{
\arrowlength #3
\put(0,0){\vector(#1,#2){\arrowlength}}%
\advance \arrowlength by-30
\put(0,0){\vector(#1,#2){\arrowlength}}%
}}

\def\horsize#1#2{%
\settowidth{\tempdimen}{$#2$}%
#1=\tempdimen
\divide #1 by\unitlength
}

\def\vertsize#1#2{%
\settoheight{\tempdimen}{$#2$}%
#1=\tempdimen
\settodepth{\tempdimen}{$#2$}%
\advance #1 by\tempdimen
\divide #1 by\unitlength
}

\def\putvector(#1,#2)(#3,#4)#5#6{{%
\ifnum3<\arrowtype
\putdashvector(#1,#2)(#3,#4)#5\arrowtype
\else
\ifnum\arrowtype<-3
\putdashvector(#1,#2)(#3,#4)#5\arrowtype
\else
\xpos=#1
\ypos=#2
\run=#3
\rise=#4
\arrowlength=#5
\ifnum \arrowtype<0
    \ifnum \run=0
        \advance \ypos by-\arrowlength
    \else
        \tempcounta \arrowlength
        \multiply \tempcounta by\rise
        \divide \tempcounta by\run
        \ifnum\run>0
            \advance \xpos by\arrowlength
            \advance \ypos by\tempcounta
        \else
            \advance \xpos by-\arrowlength
            \advance \ypos by-\tempcounta
        \fi
    \fi
    \multiply \arrowtype by-1
    \multiply \rise by-1
    \multiply \run by-1
\fi
\ifcase \arrowtype
\or \put(\xpos,\ypos){\vector(\run,\rise){\arrowlength}}%
\or \put(\xpos,\ypos){\mvector(\run,\rise)\arrowlength}%
\or \put(\xpos,\ypos){\evector(\run,\rise){\arrowlength}}%
\fi\fi\fi
}}

\def\putsplitvector(#1,#2)#3#4{
\xpos #1
\ypos #2
\arrowtype #4
\halflength #3
\arrowlength #3
\gap 140
\advance \halflength by-\gap
\divide \halflength by2
\ifnum\arrowtype>0
   \ifcase \arrowtype
   \or \put(\xpos,\ypos){\line(0,-1){\halflength}}%
       \advance\ypos by-\halflength
       \advance\ypos by-\gap
       \put(\xpos,\ypos){\vector(0,-1){\halflength}}%
   \or \put(\xpos,\ypos){\line(0,-1)\halflength}%
       \put(\xpos,\ypos){\vector(0,-1)3}%
       \advance\ypos by-\halflength
       \advance\ypos by-\gap
       \put(\xpos,\ypos){\vector(0,-1){\halflength}}%
   \or \put(\xpos,\ypos){\line(0,-1)\halflength}%
       \advance\ypos by-\halflength
       \advance\ypos by-\gap
       \put(\xpos,\ypos){\evector(0,-1){\halflength}}%
   \fi
\else \arrowtype=-\arrowtype
   \ifcase\arrowtype
   \or \advance \ypos by-\arrowlength
       \put(\xpos,\ypos){\line(0,1){\halflength}}%
       \advance\ypos by\halflength
       \advance\ypos by\gap
       \put(\xpos,\ypos){\vector(0,1){\halflength}}%
   \or \advance \ypos by-\arrowlength
       \put(\xpos,\ypos){\line(0,1)\halflength}%
       \put(\xpos,\ypos){\vector(0,1)3}%
       \advance\ypos by\halflength
       \advance\ypos by\gap
       \put(\xpos,\ypos){\vector(0,1){\halflength}}%
   \or \advance \ypos by-\arrowlength
       \put(\xpos,\ypos){\line(0,1)\halflength}%
       \advance\ypos by\halflength
       \advance\ypos by\gap
       \put(\xpos,\ypos){\evector(0,1){\halflength}}%
   \fi
\fi
}

\def\putmorphism(#1)(#2,#3)[#4`#5`#6]#7#8#9{{%
\run #2
\rise #3
\ifnum\rise=0
  \puthmorphism(#1)[#4`#5`#6]{#7}{#8}#9%
\else\ifnum\run=0
  \putvmorphism(#1)[#4`#5`#6]{#7}{#8}#9%
\else
\setpos(#1)%
\arrowlength #7
\arrowtype #8
\ifnum\run=0
\else\ifnum\rise=0
\else
\ifnum\run>0
    \coefa=1
\else
   \coefa=-1
\fi
\ifnum\arrowtype>0
   \coefb=0
   \coefc=-1
\else
   \coefb=\coefa
   \coefc=1
   \arrowtype=-\arrowtype
\fi
\width=2
\multiply \width by\run
\divide \width by\rise
\ifnum \width<0  \width=-\width\fi
\advance\width by60
\if l#9 \width=-\width\fi
\putbox(\xpos,\ypos){#4}
{\multiply \coefa by\arrowlength
\advance\xpos by\coefa
\multiply \coefa by\rise
\divide \coefa by\run
\advance \ypos by\coefa
\putbox(\xpos,\ypos){#5} }%
{\multiply \coefa by\arrowlength
\divide \coefa by2
\advance \xpos by\coefa
\advance \xpos by\width
\multiply \coefa by\rise
\divide \coefa by\run
\advance \ypos by\coefa
\if l#9%
   \putrbox(\xpos,\ypos){#6}%
\else\if r#9%
   \putlbox(\xpos,\ypos){#6}%
\fi\fi }%
{\multiply \rise by-\coefc
\multiply \run by-\coefc
\multiply \coefb by\arrowlength
\advance \xpos by\coefb
\multiply \coefb by\rise
\divide \coefb by\run
\advance \ypos by\coefb
\multiply \coefc by70
\advance \ypos by\coefc
\multiply \coefc by\run
\divide \coefc by\rise
\advance \xpos by\coefc
\multiply \coefa by140
\multiply \coefa by\run
\divide \coefa by\rise
\advance \arrowlength by\coefa
\ifcase\arrowtype
\or \put(\xpos,\ypos){\vector(\run,\rise){\arrowlength}}%
\or \put(\xpos,\ypos){\mvector(\run,\rise){\arrowlength}}%
\or \put(\xpos,\ypos){\evector(\run,\rise){\arrowlength}}%
\fi}\fi\fi\fi\fi}}

\newcount\numbdashes \newcount\lengthdash \newcount\increment

\def\howmanydashes{
\numbdashes=\arrowlength \lengthdash=40
\divide\numbdashes by \lengthdash
\lengthdash=\arrowlength
\divide\lengthdash by \numbdashes
\increment=\lengthdash
\multiply\lengthdash by 3
\divide\lengthdash by 5
}

\def\putdashvector(#1)(#2,#3)#4#5{%
\ifnum#3=0 \putdashhvector(#1){#4}#5
\else
\ifnum#2=0
\putdashvvector(#1){#4}#5\fi\fi}

\def\putdashhvector(#1,#2)#3#4{{%
\arrowlength=#3 \howmanydashes
\multiput(#1,#2)(\increment,0){\numbdashes}%
{\vrule height .4pt width \lengthdash\unitlength}
\arrowtype=#4 \xpos=#1
\ifnum\arrowtype<0 \advance\arrowtype by 7 \fi
\ifcase\arrowtype
\or \advance\xpos by 10
    \put(\xpos,#2){\vector(-1,0){\lengthdash}}
    \advance\xpos by 40
    \put(\xpos,#2){\vector(-1,0){\lengthdash}}
\or \advance \xpos by 10
    \put(\xpos,#2){\vector(-1,0){\lengthdash}}
    \advance\xpos by  \arrowlength
    \advance\xpos by  -50
    \put(\xpos,#2){\vector(-1,0){\lengthdash}}
\or \advance\xpos by 10
    \put(\xpos,#2){\vector(-1,0){\lengthdash}}
\or \advance\xpos by \arrowlength
    \advance\xpos by -\lengthdash
    \put(\xpos,#2){\vector(1,0){\lengthdash}}
\or {\advance\xpos by 10
    \put(\xpos,#2){\vector(1,0){\lengthdash}}}
    \advance\xpos by \arrowlength
    \advance\xpos by -\lengthdash
    \put(\xpos,#2){\vector(1,0){\lengthdash}}
\or \advance\xpos by \arrowlength
    \advance\xpos by -\lengthdash
    \put(\xpos,#2){\vector(1,0){\lengthdash}}
    \advance\xpos by -40
    \put(\xpos,#2){\vector(1,0){\lengthdash}}
   \fi
}}

\def\putdashvvector(#1,#2)#3#4{{%
\arrowlength=#3 \howmanydashes
\ypos=#2 \advance\ypos by -\arrowlength
\multiput(#1,#2)(0,\increment){\numbdashes}%
    {\vrule width .4pt height \lengthdash\unitlength}
\arrowtype=#4 \ypos=#2
\ifnum\arrowtype<0 \advance\arrowtype by 7 \fi
\ifcase\arrowtype
\or \advance\ypos by \arrowlength \advance\ypos by -40
    \put(#1,\ypos){\vector(0,1){\lengthdash}}
    \advance\ypos by -40
    \put(#1,\ypos){\vector(0,1){\lengthdash}}
\or \advance\ypos by 10
    \put(#1,\ypos){\vector(0,1){\lengthdash}}
    \advance\ypos by \arrowlength \advance\ypos by -40
    \put(#1,\ypos){\vector(0,1){\lengthdash}}
\or \advance\ypos by \arrowlength \advance\ypos by -40
    \put(#1,\ypos){\vector(0,1){\lengthdash}}
\or \advance\ypos by 10
    \put(#1,\ypos){\vector(0,-1){\lengthdash}}
\or \advance\ypos by 10
    \put(#1,\ypos){\vector(0,-1){\lengthdash}}
    \advance\ypos by \arrowlength \advance\ypos by -40
    \put(#1,\ypos){\vector(0,-1){\lengthdash}}
\or \advance\ypos by 10
    \put(#1,\ypos){\vector(0,-1){\lengthdash}}
    \advance\ypos by 40
    \put(#1,\ypos){\vector(0,-1){\lengthdash}}
\fi
}}

\def\puthmorphism(#1,#2)[#3`#4`#5]#6#7#8{{%
\xpos #1
\ypos #2
\width #6
\arrowlength #6
\arrowtype=#7
\putbox(\xpos,\ypos){#3\vphantom{#4}}%
{\advance \xpos by\arrowlength
\putbox(\xpos,\ypos){\vphantom{#3}#4}}%
\horsize{\tempcounta}{#3}%
\horsize{\tempcountb}{#4}%
\divide \tempcounta by2
\divide \tempcountb by2
\advance \tempcounta by30
\advance \tempcountb by30
\advance \xpos by\tempcounta
\advance \arrowlength by-\tempcounta
\advance \arrowlength by-\tempcountb
\putvector(\xpos,\ypos)(1,0)\arrowlength\arrowtype
\divide \arrowlength by2
\advance \xpos by\arrowlength
\vertsize{\tempcounta}{#5}%
\divide\tempcounta by2
\advance \tempcounta by20
\if a#8 %
   \advance \ypos by\tempcounta
   \putbox(\xpos,\ypos){#5}%
\else
   \advance \ypos by-\tempcounta
   \putbox(\xpos,\ypos){#5}%
\fi}}

\def\putvmorphism(#1,#2)[#3`#4`#5]#6#7#8{{%
\xpos #1
\ypos #2
\arrowlength #6
\arrowtype #7
\settowidth{\xlen}{$#5$}%
\putbox(\xpos,\ypos){#3}%
{\advance \ypos by-\arrowlength
\putbox(\xpos,\ypos){#4}}%
{\advance\arrowlength by-140
\advance \ypos by-70
\ifdim\xlen>0pt
   \if m#8%
      \putsplitvector(\xpos,\ypos)\arrowlength\arrowtype
   \else
   \putvector(\xpos,\ypos)(0,-1)\arrowlength\arrowtype
   \fi
\else
   \putvector(\xpos,\ypos)(0,-1)\arrowlength\arrowtype
\fi}%
\ifdim\xlen>0pt
   \divide \arrowlength by2
   \advance\ypos by-\arrowlength
   \if l#8%
      \advance \xpos by-40
      \putrbox(\xpos,\ypos){#5}%
   \else\if r#8%
      \advance \xpos by40
      \putlbox(\xpos,\ypos){#5}%
   \else
      \putbox(\xpos,\ypos){#5}%
   \fi\fi
\fi
}}

\def\putsquarep<#1>(#2)[#3;#4`#5`#6`#7]{{%
\setsqparms[#1]%
\setpos(#2)%
\settokens[#3]%
\puthmorphism(\xpos,\ypos)[\tokenc`\tokend`{#7}]{\width}{\arrowtyped}b%
\advance\ypos by \height
\puthmorphism(\xpos,\ypos)[\tokena`\tokenb`{#4}]{\width}{\arrowtypea}a%
\putvmorphism(\xpos,\ypos)[``{#5}]{\height}{\arrowtypeb}l%
\advance\xpos by \width
\putvmorphism(\xpos,\ypos)[``{#6}]{\height}{\arrowtypec}r%
}}

\def\putsquare{\@ifnextchar <{\putsquarep}{\putsquarep%
   <\arrowtypea`\arrowtypeb`\arrowtypec`\arrowtyped;\width`\height>}}
\def\square{\@ifnextchar< {\squarep}{\squarep
   <\arrowtypea`\arrowtypeb`\arrowtypec`\arrowtyped;\width`\height>}}
\def\squarep<#1>[#2`#3`#4`#5;#6`#7`#8`#9]{{
\setsqparms[#1]
\diagram
\putsquarep<\arrowtypea`\arrowtypeb`\arrowtypec`
\arrowtyped;\width`\height>
(0,0)[#2`#3`#4`{#5};#6`#7`#8`{#9}]
\enddiagram
}}                                                 
\def\putptrianglep<#1>(#2,#3)[#4`#5`#6;#7`#8`#9]{{%
\settriparms[#1]%
\xpos=#2 \ypos=#3
\advance\ypos by \height
\puthmorphism(\xpos,\ypos)[#4`#5`{#7}]{\height}{\arrowtypea}a%
\putvmorphism(\xpos,\ypos)[`#6`{#8}]{\height}{\arrowtypeb}l%
\advance\xpos by\height
\putmorphism(\xpos,\ypos)(-1,-1)[``{#9}]{\height}{\arrowtypec}r%
}}

\def\putptriangle{\@ifnextchar <{\putptrianglep}{\putptrianglep
   <\arrowtypea`\arrowtypeb`\arrowtypec;\height>}}
\def\ptriangle{\@ifnextchar <{\ptrianglep}{\ptrianglep
   <\arrowtypea`\arrowtypeb`\arrowtypec;\height>}}
\def\ptrianglep<#1>[#2`#3`#4;#5`#6`#7]{{
\settriparms[#1]
\diagram
\putptrianglep<\arrowtypea`\arrowtypeb`
\arrowtypec;\height>
(0,0)[#2`#3`#4;#5`#6`{#7}]
\enddiagram
}}                                            

\def\putqtrianglep<#1>(#2,#3)[#4`#5`#6;#7`#8`#9]{{%
\settriparms[#1]%
\xpos=#2 \ypos=#3
\advance\ypos by\height
\puthmorphism(\xpos,\ypos)[#4`#5`{#7}]{\height}{\arrowtypea}a%
\putmorphism(\xpos,\ypos)(1,-1)[``{#8}]{\height}{\arrowtypeb}l%
\advance\xpos by\height
\putvmorphism(\xpos,\ypos)[`#6`{#9}]{\height}{\arrowtypec}r%
}}

\def\putqtriangle{\@ifnextchar <{\putqtrianglep}{\putqtrianglep
   <\arrowtypea`\arrowtypeb`\arrowtypec;\height>}}
\def\qtriangle{\@ifnextchar <{\qtrianglep}{\qtrianglep
   <\arrowtypea`\arrowtypeb`\arrowtypec;\height>}}
\def\qtrianglep<#1>[#2`#3`#4;#5`#6`#7]{{
\settriparms[#1]
\width=\height                                
\diagram
\putqtrianglep<\arrowtypea`\arrowtypeb`
\arrowtypec;\height>
(0,0)[#2`#3`#4;#5`#6`{#7}]
\enddiagram
}}

\def\putdtrianglep<#1>(#2,#3)[#4`#5`#6;#7`#8`#9]{{%
\settriparms[#1]%
\xpos=#2 \ypos=#3
\puthmorphism(\xpos,\ypos)[#5`#6`{#9}]{\height}{\arrowtypec}b%
\advance\xpos by \height \advance\ypos by\height
\putmorphism(\xpos,\ypos)(-1,-1)[``{#7}]{\height}{\arrowtypea}l%
\putvmorphism(\xpos,\ypos)[#4``{#8}]{\height}{\arrowtypeb}r%
}}

\def\putdtriangle{\@ifnextchar <{\putdtrianglep}{\putdtrianglep
   <\arrowtypea`\arrowtypeb`\arrowtypec;\height>}}
\def\dtriangle{\@ifnextchar <{\dtrianglep}{\dtrianglep
   <\arrowtypea`\arrowtypeb`\arrowtypec;\height>}}
\def\dtrianglep<#1>[#2`#3`#4;#5`#6`#7]{{
\settriparms[#1]
\width=\height                                
\diagram
\putdtrianglep<\arrowtypea`\arrowtypeb`
\arrowtypec;\height>
(0,0)[#2`#3`#4;#5`#6`{#7}]
\enddiagram
}}

\def\putbtrianglep<#1>(#2,#3)[#4`#5`#6;#7`#8`#9]{{%
\settriparms[#1]%
\xpos=#2 \ypos=#3
\puthmorphism(\xpos,\ypos)[#5`#6`{#9}]{\height}{\arrowtypec}b%
\advance\ypos by\height
\putmorphism(\xpos,\ypos)(1,-1)[``{#8}]{\height}{\arrowtypeb}r%
\putvmorphism(\xpos,\ypos)[#4``{#7}]{\height}{\arrowtypea}l%
}}

\def\putbtriangle{\@ifnextchar <{\putbtrianglep}{\putbtrianglep
   <\arrowtypea`\arrowtypeb`\arrowtypec;\height>}}
\def\btriangle{\@ifnextchar <{\btrianglep}{\btrianglep
   <\arrowtypea`\arrowtypeb`\arrowtypec;\height>}}
\def\btrianglep<#1>[#2`#3`#4;#5`#6`#7]{{
\settriparms[#1]
\width=\height                               
\diagram
\putbtrianglep<\arrowtypea`\arrowtypeb`
\arrowtypec;\height>
(0,0)[#2`#3`#4;#5`#6`{#7}]
\enddiagram
}}

\def\putAtrianglep<#1>(#2,#3)[#4`#5`#6;#7`#8`#9]{{%
\settriparms[#1]%
\xpos=#2 \ypos=#3
{\multiply \height by2
\puthmorphism(\xpos,\ypos)[#5`#6`{#9}]{\height}{\arrowtypec}b}%
\advance\xpos by\height \advance\ypos by\height
\putmorphism(\xpos,\ypos)(-1,-1)[#4``{#7}]{\height}{\arrowtypea}l%
\putmorphism(\xpos,\ypos)(1,-1)[``{#8}]{\height}{\arrowtypeb}r%
}}

\def\putAtriangle{\@ifnextchar <{\putAtrianglep}{\putAtrianglep
   <\arrowtypea`\arrowtypeb`\arrowtypec;\height>}}
\def\Atriangle{\@ifnextchar <{\Atrianglep}{\Atrianglep
   <\arrowtypea`\arrowtypeb`\arrowtypec;\height>}}
\def\Atrianglep<#1>[#2`#3`#4;#5`#6`#7]{{
\settriparms[#1]
\width=\height                                     
\diagram
\putAtrianglep<\arrowtypea`\arrowtypeb`
\arrowtypec;\height>
(0,0)[#2`#3`#4;#5`#6`{#7}]
\enddiagram
}}

\def\putAtrianglepairp<#1>(#2)[#3;#4`#5`#6`#7`#8]{{%
\settripairparms[#1]%
\setpos(#2)%
\settokens[#3]%
\puthmorphism(\xpos,\ypos)[\tokenb`\tokenc`{#7}]{\height}{\arrowtyped}b%
\advance\xpos by\height
\puthmorphism(\xpos,\ypos)[\phantom{\tokenc}`\tokend`{#8}]%
{\height}{\arrowtypee}b%
\advance\ypos by\height
\putmorphism(\xpos,\ypos)(-1,-1)[\tokena``{#4}]{\height}{\arrowtypea}l%
\putvmorphism(\xpos,\ypos)[``{#5}]{\height}{\arrowtypeb}m%
\putmorphism(\xpos,\ypos)(1,-1)[``{#6}]{\height}{\arrowtypec}r%
}}

\def\putAtrianglepair{\@ifnextchar <{\putAtrianglepairp}{\putAtrianglepairp%
   <\arrowtypea`\arrowtypeb`\arrowtypec`\arrowtyped`\arrowtypee;\height>}}
\def\Atrianglepair{\@ifnextchar <{\Atrianglepairp}{\Atrianglepairp%
   <\arrowtypea`\arrowtypeb`\arrowtypec`\arrowtyped`\arrowtypee;\height>}}

\def\Atrianglepairp<#1>[#2;#3`#4`#5`#6`#7]{{
\settripairparms[#1]
\settokens[#2]
\width=\height                                
\diagram
\putAtrianglepairp                            
<\arrowtypea`\arrowtypeb`\arrowtypec`
\arrowtyped`\arrowtypee;\height>
(0,0)[{#2};#3`#4`#5`#6`{#7}]
\enddiagram
}}

\def\putVtrianglep<#1>(#2,#3)[#4`#5`#6;#7`#8`#9]{{%
\settriparms[#1]%
\xpos=#2 \ypos=#3
\advance\ypos by\height
{\multiply\height by2
\puthmorphism(\xpos,\ypos)[#4`#5`{#7}]{\height}{\arrowtypea}a}%
\putmorphism(\xpos,\ypos)(1,-1)[`#6`{#8}]{\height}{\arrowtypeb}l%
\advance\xpos by\height
\advance\xpos by\height
\putmorphism(\xpos,\ypos)(-1,-1)[``{#9}]{\height}{\arrowtypec}r%
}}

\def\putVtriangle{\@ifnextchar <{\putVtrianglep}{\putVtrianglep
   <\arrowtypea`\arrowtypeb`\arrowtypec;\height>}}
\def\Vtriangle{\@ifnextchar <{\Vtrianglep}{\Vtrianglep
   <\arrowtypea`\arrowtypeb`\arrowtypec;\height>}}
\def\Vtrianglep<#1>[#2`#3`#4;#5`#6`#7]{{
\settriparms[#1]
\width=\height                                 
\diagram
\putVtrianglep<\arrowtypea`\arrowtypeb`
\arrowtypec;\height>
(0,0)[#2`#3`#4;#5`#6`{#7}]
\enddiagram
}}

\def\putVtrianglepairp<#1>(#2)[#3;#4`#5`#6`#7`#8]{{
\settripairparms[#1]%
\setpos(#2)%
\settokens[#3]%
\advance\ypos by\height
\putmorphism(\xpos,\ypos)(1,-1)[`\tokend`{#6}]{\height}{\arrowtypec}l%
\puthmorphism(\xpos,\ypos)[\tokena`\tokenb`{#4}]{\height}{\arrowtypea}a%
\advance\xpos by\height
\puthmorphism(\xpos,\ypos)[\phantom{\tokenb}`\tokenc`{#5}]%
{\height}{\arrowtypeb}a%
\putvmorphism(\xpos,\ypos)[``{#7}]{\height}{\arrowtyped}m%
\advance\xpos by\height
\putmorphism(\xpos,\ypos)(-1,-1)[``{#8}]{\height}{\arrowtypee}r%
}}

\def\putVtrianglepair{\@ifnextchar <{\putVtrianglepairp}{\putVtrianglepairp%
    <\arrowtypea`\arrowtypeb`\arrowtypec`\arrowtyped`\arrowtypee;\height>}}
\def\Vtrianglepair{\@ifnextchar <{\Vtrianglepairp}{\Vtrianglepairp%
    <\arrowtypea`\arrowtypeb`\arrowtypec`\arrowtyped`\arrowtypee;\height>}}
\def\Vtrianglepairp<#1>[#2;#3`#4`#5`#6`#7]{{
\settripairparms[#1]
\settokens[#2]
\diagram
\putVtrianglepairp                             
<\arrowtypea`\arrowtypeb`\arrowtypec`
\arrowtyped`\arrowtypee;\height>
(0,0)[{#2};#3`#4`#5`#6`{#7}]
\enddiagram
}}

\def\putCtrianglep<#1>(#2,#3)[#4`#5`#6;#7`#8`#9]{{%
\settriparms[#1]%
\xpos=#2 \ypos=#3
\advance\ypos by\height
\putmorphism(\xpos,\ypos)(1,-1)[``{#9}]{\height}{\arrowtypec}l%
\advance\xpos by\height
\advance\ypos by\height
\putmorphism(\xpos,\ypos)(-1,-1)[#4`#5`{#7}]{\height}{\arrowtypea}l%
{\multiply\height by 2
\putvmorphism(\xpos,\ypos)[`#6`{#8}]{\height}{\arrowtypeb}r}%
}}

\def\putCtriangle{\@ifnextchar <{\putCtrianglep}{\putCtrianglep
    <\arrowtypea`\arrowtypeb`\arrowtypec;\height>}}
\def\Ctriangle{\@ifnextchar <{\Ctrianglep}{\Ctrianglep
    <\arrowtypea`\arrowtypeb`\arrowtypec;\height>}}
\def\Ctrianglep<#1>[#2`#3`#4;#5`#6`#7]{{
\settriparms[#1]
\width=\height                               
\diagram
\putCtrianglep<\arrowtypea`\arrowtypeb`
\arrowtypec;\height>
(0,0)[#2`#3`#4;#5`#6`{#7}]
\enddiagram
}}                                           
\def\putDtrianglep<#1>(#2,#3)[#4`#5`#6;#7`#8`#9]{{%
\settriparms[#1]%
\xpos=#2 \ypos=#3
\advance\xpos by\height \advance\ypos by\height
\putmorphism(\xpos,\ypos)(-1,-1)[``{#9}]{\height}{\arrowtypec}r%
\advance\xpos by-\height \advance\ypos by\height
\putmorphism(\xpos,\ypos)(1,-1)[`#5`{#8}]{\height}{\arrowtypeb}r%
{\multiply\height by 2
\putvmorphism(\xpos,\ypos)[#4`#6`{#7}]{\height}{\arrowtypea}l}%
}}

\def\putDtriangle{\@ifnextchar <{\putDtrianglep}{\putDtrianglep
    <\arrowtypea`\arrowtypeb`\arrowtypec;\height>}}
\def\Dtriangle{\@ifnextchar <{\Dtrianglep}{\Dtrianglep
   <\arrowtypea`\arrowtypeb`\arrowtypec;\height>}}
\def\Dtrianglep<#1>[#2`#3`#4;#5`#6`#7]{{
\settriparms[#1]
\width=\height                              
\diagram
\putDtrianglep<\arrowtypea`\arrowtypeb`
\arrowtypec;\height>
(0,0)[#2`#3`#4;#5`#6`{#7}]
\enddiagram
}}                                          
\def\setrecparms[#1`#2]{\width=#1 \height=#2}%

\def\recursep<#1`#2>[#3;#4`#5`#6`#7`#8]{{%
\width=#1 \height=#2
\settokens[#3]
\settowidth{\tempdimen}{$\tokena$}
\ifdim\tempdimen=0pt
  \savebox{\tempboxa}{\hbox{$\tokenb$}}%
  \savebox{\tempboxb}{\hbox{$\tokend$}}%
  \savebox{\tempboxc}{\hbox{$#6$}}%
\else
  \savebox{\tempboxa}{\hbox{$\hbox{$\tokena$}\times\hbox{$\tokenb$}$}}%
  \savebox{\tempboxb}{\hbox{$\hbox{$\tokena$}\times\hbox{$\tokend$}$}}%
  \savebox{\tempboxc}{\hbox{$\hbox{$\tokena$}\times\hbox{$#6$}$}}%
\fi
\ypos=\height
\divide\ypos by 2
\xpos=\ypos
\advance\xpos by \width
\bfig
\putCtrianglep<-1`1`1;\ypos>(0,0)[`\tokenc`;#5`#6`{#7}]%
\puthmorphism(\ypos,0)[\tokend`\usebox{\tempboxb}`{#8}]{\width}{-1}b%
\puthmorphism(\ypos,\height)[\tokenb`\usebox{\tempboxa}`{#4}]{\width}{-1}a%
\advance\ypos by \width
\putvmorphism(\ypos,\height)[``\usebox{\tempboxc}]{\height}1r%
\efig
}}

\def\recurse{\@ifnextchar <{\recursep}{\recursep<\width`\height>}}

\def\puttwohmorphisms(#1,#2)[#3`#4;#5`#6]#7#8#9{{%
\puthmorphism(#1,#2)[#3`#4`]{#7}0a
\ypos=#2
\advance\ypos by 20
\puthmorphism(#1,\ypos)[\phantom{#3}`\phantom{#4}`#5]{#7}{#8}a
\advance\ypos by -40
\puthmorphism(#1,\ypos)[\phantom{#3}`\phantom{#4}`#6]{#7}{#9}b
}}

\def\puttwovmorphisms(#1,#2)[#3`#4;#5`#6]#7#8#9{{%
\putvmorphism(#1,#2)[#3`#4`]{#7}0a
\xpos=#1
\advance\xpos by -20
\putvmorphism(\xpos,#2)[\phantom{#3}`\phantom{#4}`#5]{#7}{#8}l
\advance\xpos by 40
\putvmorphism(\xpos,#2)[\phantom{#3}`\phantom{#4}`#6]{#7}{#9}r
}}

\def\puthcoequalizer(#1)[#2`#3`#4;#5`#6`#7]#8#9{{%
\setpos(#1)%
\puttwohmorphisms(\xpos,\ypos)[#2`#3;#5`#6]{#8}11%
\advance\xpos by #8
\puthmorphism(\xpos,\ypos)[\phantom{#3}`#4`#7]{#8}1{#9}
}}

\def\putvcoequalizer(#1)[#2`#3`#4;#5`#6`#7]#8#9{{%
\setpos(#1)%
\puttwovmorphisms(\xpos,\ypos)[#2`#3;#5`#6]{#8}11%
\advance\ypos by -#8
\putvmorphism(\xpos,\ypos)[\phantom{#3}`#4`#7]{#8}1{#9}
}}

\def\putthreehmorphisms(#1)[#2`#3;#4`#5`#6]#7(#8)#9{{%
\setpos(#1) \settypes(#8)
\if a#9 %
     \vertsize{\tempcounta}{#5}%
     \vertsize{\tempcountb}{#6}%
     \ifnum \tempcounta<\tempcountb \tempcounta=\tempcountb \fi
\else
     \vertsize{\tempcounta}{#4}%
     \vertsize{\tempcountb}{#5}%
     \ifnum \tempcounta<\tempcountb \tempcounta=\tempcountb \fi
\fi
\advance \tempcounta by 60
\puthmorphism(\xpos,\ypos)[#2`#3`#5]{#7}{\arrowtypeb}{#9}
\advance\ypos by \tempcounta
\puthmorphism(\xpos,\ypos)[\phantom{#2}`\phantom{#3}`#4]{#7}{\arrowtypea}{#9}
\advance\ypos by -\tempcounta \advance\ypos by -\tempcounta
\puthmorphism(\xpos,\ypos)[\phantom{#2}`\phantom{#3}`#6]{#7}{\arrowtypec}{#9}
}}

\def\setarrowtoks[#1`#2`#3`#4`#5`#6]{%
\def\toka{#1}
\def\tokb{#2}
\def\tokc{#3}
\def\tokd{#4}
\def\toke{#5}
\def\tokf{#6}
}
\def\hex{\@ifnextchar <{\hexp}{\hexp<1000`400>}}
\def\hexp<#1`#2>[#3`#4`#5`#6`#7`#8;#9]{%
\setarrowtoks[#9]
\yext=#2 \advance \yext by #2
\xext=#1 \advance\xext by \yext
\bfig
\putCtriangle<-1`0`1;#2>(0,0)[`#5`;\tokb``\tokd]
\xext=#1 \yext=#2 \advance \yext by #2
\putsquare<1`0`0`1;\xext`\yext>(#2,0)[#3`#4`#7`#8;\toka```\tokf]
\advance \xext by #2
\putDtriangle<0`1`-1;#2>(\xext,0)[`#6`;`\tokc`\toke]
\efig
}

\QQQ{Language}{
American English
}

\begin{document}

\title{{\bf \ Higher algebraic }$K${\bf -theory of group actions}\\
{\bf with finite stabilizers}}
\author{Gabriele Vezzosi$^{*}$, Angelo Vistoli\thanks{%
Partially supported by the University of Bologna, funds for selected
research topics.} \\
Dipartimento di Matematica, Universit\`a di Bologna,\\
Piazza di Porta S. Donato 5, Bologna 40127 - Italy}
\maketitle

\begin{abstract}
We prove a decomposition theorem for the equivariant $K$-theory of actions
of affine group schemes $G$ of finite type over a field on regular separated
noetherian algebraic spaces, under the hypothesis that the actions have
finite geometric stabilizers and satisfy a rationality condition together
with a technical condition which holds e.g. for $G$ abelian or smooth.

We reduce the problem to the case of a ${\rm GL}_n$-action and finally to a
split torus action.

\smallskip\ 

\end{abstract}

\tableofcontents

\section{Introduction}

The purpose of this paper is to prove a decomposition theorem for the
equivariant $K$-theory of actions of affine group schemes of finite type
over a field on a regular separated noetherian algebraic spaces.{\sl \ }Let $%
X$ be a regular connected separated noetherian scheme with an ample line
bundle, $K_0(X)$ its Grothendieck ring of vector bundles. Then the kernel of
the rank morphism $K_0\left( X\right) \rightarrow {\Bbb Z}$ is nilpotent (%
{\bf \cite{SGA6}}, VI, Th\'eor\`eme 6.9), so the ring $K_0\left( X\right) $
is indecomposable, and remains such after tensoring with any indecomposable $%
{\Bbb Z}$-algebra.{\sl {\bf \ }}

The situation is quite different when we consider the equivariant case. Let $%
G$ be an algebraic group acting on a noetherian separated regular scheme, or
algebraic space, $X$ over a field $k$, and consider the Grothendieck ring $%
K_0(X,G)$ of $G$-equivariant perfect complexes. This is the same as the
Grothendieck group of $G$-equivariant coherent sheaves on $X$, and coincides
with the Grothendieck ring of $G$-equivariant vector bundles if all $G$%
-coherent sheaves are quotients of locally free coherent sheaves (which is
the case, e.g. when $G$ is finite or smooth and $X$ is a scheme). Assume
that the action of $G$ on $X$ is connected, that is, there are no nontrivial
invariant open and closed subschemes of $X$. Still, $K_0(X,G)$ will usually
decompose, after inverting some primes; for example, if $G$ is a finite
group and $X={\rm Spec}{\Bbb C}${\sl {\bf , }}then $K_0(X,G)$ is the ring of
complex representations of $G$, which becomes a product of fields after
tensoring with ${\Bbb Q}$.

In {\bf \cite{Vi1}} the second author analyzes the case that the action of $%
G $ on $X$ has finite reduced geometric stabilizers. Consider the ring of
representations ${\rm R}\left( G\right) $, and the kernel{\bf \ }${\frak m}$
of the rank morphism ${\rm rk:}K_0(X,G)\to {\Bbb Z}$. Then $K_0(X,G)$ is an $%
{\rm R}\left( G\right) $-algebra; he shows that the localization morphism

\[
K_0(X,G)\otimes {\Bbb Q}\longrightarrow K_0(X,G){\rm _{{\frak m}}} 
\]
is surjective, and that the kernel of the rank morphism $K_0(X,G)_{{\frak m}%
}\otimes {\Bbb Q}\to {\Bbb Q}$ is nilpotent. Furthermore he conjectures that 
$K_0(X,G)\otimes {\Bbb Q}$ splits as a product of the localization $K_0(X,G)%
{\rm _{{\frak m}}}$ and some other ring, and formulates a conjecture about
what the other factor ring should be when $G$ is abelian and the field is
algebraically closed of characteristic 0. The proofs of the results in {\bf 
\cite{Vi1}} depend on an equivariant Riemann-Roch theorem, whose proof was
never published by the author; however, all of the results have been proved
and generalized in {\bf \cite{EG-GRR}}.

The case that $G$ is a finite group is studied in {\bf \cite{Vi2}}. Assume
that $k$ contains all $n^{{\rm th}}$ roots of 1, where $n$ is the order of
the group $G$. Then the author shows that, after inverting the order of $G$,
the $K$-theory ring $K_{*}(X,G)$ of $G$-equivariant vector bundles on $X$
(which is assumed to be a scheme in that paper) is canonically the product
of a finite number of rings, expressible in terms of ordinary $K$-theory of
appropriate subschemes of fixed points of $X$. Here $K_{*}(X,G)=%
\bigoplus_iK_i(X,G)$ is the graded higher $K$-theory ring. The precise
formula is as follows.

Let $\sigma $ be a cyclic subgroup of $G$ whose order is prime to the
characteristic of $k$; then the subscheme $X^\sigma $ of fixed points of $X$
under the actions of $\sigma $ is also regular. The representation ring $%
{\rm R}\left( \sigma \right) $ is isomorphic to the ring ${\Bbb Z}%
[t]/(t^n-1) $, where $t$ is a generator of the group of characters $\hom
(\sigma ,k^{*})$. We call ${\rm \widetilde{R}}\left( \sigma \right) $ the
quotient of the ring ${\rm R}\left( \sigma \right) $ by the ideal generated
by the element $\Phi _n(t)$, where $\Phi _n$ is the $n^{{\rm th}}$
cyclotomic polynomial; this is independent of $t$. The ring ${\rm \widetilde{%
R}}\left( \sigma \right) $ is isomorphic to the ring of integers in the $n^{%
{\rm th}}$ cyclotomic field. Call ${\rm N}_G\left( \sigma \right) $ the
normalizer of $\sigma $ in $G$; the group ${\rm N}_G\left( \sigma \right) $
acts on the scheme $X^\sigma $, and, by conjugation, on the group $\sigma $.
Consider the induced actions of ${\rm N}_G\left( \sigma \right) $ on the $K$%
-theory ring $K_{*}(X^\sigma )$, and on the ring ${\rm \widetilde{R}}\left(
\sigma \right) $.

Choose a set ${\cal C}\left( G\right) $ of representatives for the conjugacy
classes of cyclic subgroups of $G$ whose order is prime to the
characteristic of the field. The statement of the main result of {\bf \cite
{Vi2}} is as follows.

{\bf Theorem} {\em There is a canonical ring isomorphism } 
\[
K_{*}\left( X,G\right) \otimes {\Bbb Z}\left[ 1/\left| G\right| \right]
\simeq \prod_{\sigma \in {\cal C}\left( G\right) }\left( K_{*}\left(
X^\sigma \right) \otimes \widetilde{{\rm R}}\left( \sigma \right) \right) ^{%
{\rm N}_G\left( \sigma \right) }\otimes {\Bbb Z}\left[ 1/\left| G\right|
\right] . 
\]

In the present paper we generalize this decomposition to the case in which $%
G $ is an algebraic group scheme of finite type over a field $k$, acting on
a noetherian regular separated algebraic space $X$ over $k$. Of course, we
cannot expect a statement exactly like the one for finite groups, expressing
equivariant $K$-theory simply in terms of ordinary $K$-theory of the fixed
point sets. For example, when $X$ is the Stiefel variety of $r$-frames in $n$%
-space, then the quotient of $X$ by the natural free action of ${\rm GL}_r$
is the Grassmannian of $r$-planes in $n$-space, and $K_0(X,{\rm GL}_r)=K_0(X/%
{\rm GL}_r)$ is nontrivial, while $K_0(X)={\Bbb Z}$.

Let $X$ be a noetherian regular separated algebraic space of over $k$ with
an action of an affine group scheme $G$ of finite type over $k$. We consider
the Waldhausen $K$-theory group $K_{*}(X,G)$ of complexes of quasicoherent $%
G $-equivariant sheaves on $X$ with coherent bounded cohomology. This
coincides on one hand with the Waldhausen $K$-theory group $K_{*}(X,G)$ of
the subcategory of complexes of quasicoherent $G$-equivariant flat sheaves
on $X$ with coherent bounded cohomology and hence has a natural ring
structure given by the total tensor product and on the other with the
Quillen group $K_{*}^{\prime }(X,G)$ of coherent equivariant sheaves on $X$;
furthermore, if every coherent equivariant sheaf on $X$ is the quotient of a
locally free equivariant coherent sheaf, it also coincides with the Quillen
group $K_{*}^{naive}(X,G)$ of coherent locally free equivariant sheaves on $%
X $. These $K$-theories and their relationships are discussed in the
Appendix.

Our result is as follows. First we have to see what plays the role of the
cyclic subgroups of a finite group. This is easy; the group schemes whose
rings of representations are of the form ${\Bbb Z}[t]/(t^n-1)$ are not the
cyclic groups, in general, but their Cartier duals, that is, the group
schemes which are isomorphic to the group scheme $\mu _n$ of $n^{{\rm th}}$
roots of 1 for some $n$. We call these group schemes ``dual cyclic''. If $%
\sigma $ is a dual cyclic group, we can define ${\rm \widetilde{R}}\sigma $
as before. A dual cyclic subgroup $\sigma $ of $G$ is called essential if $%
X^\sigma \neq \emptyset $.

The correct substitute for the ordinary $K$-theory of the subspaces of
invariants is the {\em geometric equivariant }$K${\em -theory\/} $%
K_{*}(X,G)_{{\rm geom}}$, which is defined as follows. Call $N$ the least
common multiple of the orders of all the essential dual cyclic subgroups of $%
G$. Call $S_1$ be the multiplicative subset of the ring ${\rm R}\left(
G\right) $ consisting of elements whose virtual rank is a power of $N$; then 
$K_{*}(X,G)_{{\rm geom}}$ is the localization $S_1^{-1}K_{*}(X,G)$. Notice
that $K_{*}(X,G)_{{\rm geom}}\otimes {\Bbb Q}=K_{*}(X,G)_{{\frak m}}$, with
the notation above. Moreover, if every coherent equivariant sheaf on $X$ is
the quotient of a locally free equivariant coherent sheaf, by {\bf \cite
{EG-GRR}}, we have an isomorphism of rings 
\[
K_0\left( X,G\right) _{{\rm geom}}\otimes {\Bbb Q}=A_G^{*}\left( X\right)
\otimes {\Bbb Q} 
\]
where $A_G^{*}\left( X\right) $ denotes the direct sum of $G$-equivariant
Chow groups of $X.$

We prove the following. Assume that the action of $G$ on $X$ is connected.
Then the kernel of the rank morphism $K_0(X,G)_{{\rm geom}}\to {\Bbb Z}[1/N]$
is nilpotent (see Corollary \ref{nilpotent0}). This is remarkable; we have
made what might look like a small step toward making the equivariant $K$%
-theory ring indecomposable, and immediately get an indecomposable ring.
Indeed, $K_{*}(X,G)_{{\rm geom}}$ ``feels like'' the $K$-theory ring of a
scheme; we want to think of $K_{*}(X,G)_{{\rm geom}}$ as what the $K$-theory
of the quotient $X/G$ should be, if $X/G$ were smooth, after inverting $N$
(see Conjecture~\ref{conj}).

Furthermore, consider the centralizer ${\rm C}_G\left( \sigma \right) $ and
the normalizer ${\rm N}_G\left( \sigma \right) $ of $\sigma $ inside $G$.
The quotient ${\rm w}_G\left( \sigma \right) ={\rm N}_G\left( \sigma \right)
/{\rm C}_G\left( \sigma \right) $ is contained inside the group scheme of
automorphisms of $\sigma $, which is a discrete group, so it is also a
discrete group. It acts on $\widetilde{{\rm R}}\left( \sigma \right) $, by
conjugation, and it also acts on the equivariant $K$-theory ring $%
K_{*}(X^\sigma ,{\rm C}_G\left( \sigma \right) )$, and on the geometric
equivariant $K$-theory ring $K_{*}(X^\sigma ,{\rm C}_G\left( \sigma \right)
)_{{\rm geom}}$ (see Corollary \ref{azione}).

We say that the action of $G$ on $X$ is sufficiently rational when the
following two conditions are satisfied. Let $\overline{k}$ be the algebraic
closure of $k$.

\begin{enumerate}
\item  each essential dual cyclic subgroup $\sigma \subseteq G_{\overline{k}}
$ is conjugate by an element of $G(\,\overline{k}\,)$ to a dual cyclic
subgroup of $G$;

\item  if two essential dual cyclic subgroups of $G$ are conjugate by an
element of $G(\,\overline{k}\,)$, they are also conjugate by an element of $%
G(k)$.
\end{enumerate}

\smallbreak

Obviously every action over an algebraically closed field is sufficiently
rational. Furthermore if $G$ is ${\rm GL}_m$, ${\rm SL}_m$, ${\rm Sp}_m$ or
a totally split torus, then any action of $G$ is sufficiently rational over
an arbitrary base field (Proposition \ref{angelo}). If $G$ is a finite
group, then the action is sufficiently rational when $k$ contains all $n$-th
roots of 1, where $n$ is the least common multiple of the orders of the
cyclic subgroups of $k$ of order prime to the characteristic, whose fixed
point subscheme is nonempty. Denote by ${\cal C}\left( G\right) $ a set of
representatives for essential dual cyclic subgroup schemes, under
conjugation by elements of the group $G(k)$. Here is the statement of our
result.

{\bf Main Theorem.} {\em Let }$G${\em \ be an affine group scheme of finite
type over a field }$k${\em , acting on a noetherian separated regular
algebraic space }$X${\em . Assume the following three conditions.}

\begin{enumerate}
\item[(a)]  {\em The action has finite geometric stabilizers. }

\item[(b)]  {\em The action is sufficiently rational. }

\item[(c)]  {\em For any essential cyclic subgroup }$\sigma ${\em \ of }$G$%
{\em , the quotient }$G/{\rm C}_G\left( \sigma \right) ${\em \ is smooth. }%
\smallbreak 
\end{enumerate}

{\em Then }${\cal C}\left( G\right) ${\em \ is finite, and there is a
canonical isomorphism of }${\rm R}\left( G\right) ${\em -algebras} 
\[
K_{*}\left( X,G\right) \otimes {\Bbb Z}\left[ 1/N\right] \simeq
\prod_{\sigma \in {\cal C}\left( G\right) }(K_{*}(X^\sigma ,{\rm C}_G(\sigma
))_{{\rm geom}}\otimes \widetilde{{\rm R}}\left( \sigma \right) )^{{\rm w}%
_G\left( \sigma \right) }. 
\]

Conditions (a) and (b) are clearly necessary for the theorem to hold. We are
not sure about (c). It is rather mild, as it is satisfied, for example, when 
$G$ is smooth (this is automatically true in characteristic $0$) or when $G$
is abelian. A weaker version of condition (c) is given in Subsection \ref
{hypothesesonG}.

In case that $G$ is abelian over an algebraically closed field of
characteristic $0$, the Main Theorem implies Conjecture~3.6 in {\bf \cite
{Vi1}}. When $G$ is a finite group, and the base field contains enough roots
of 1, as in the statement of Theorem~1, then the conditions of the Main
Theorem are satisfied; since the natural maps $K_{*}\left( X^\sigma ,{\rm C}%
_G\left( \sigma \right) \right) _{{\rm geom}}\to K_{*}(X^\sigma )^{{\rm C}%
_G\left( \sigma \right) }$ become isomorphisms after inverting the order of $%
G$ (Proposition \ref{finitegroup}) the Main Theorem implies {\bf \cite{Vi2} }%
Theorem~1. However, the proof of the main theorem here is completely
different from the proof of Theorem~1 in {\bf \cite{Vi2}}.

As B. Toen pointed out to us, a weaker version of our main theorem (with $%
{\Bbb Q}$-coefficients and assuming $G$ smooth, acting with finite reduced
stabilizers) follows from his Th\'eor\`eme 3.15 in {\bf \cite{To1}}; the
\'etale techniques he uses in proving this result make it impossible to
avoid tensoring with ${\Bbb Q}$. See also {\bf \cite{To2}.}

\smallskip\ 

Here is an outline of the paper.

First we define the homomorphism (Section 2.2). Next, in Section 3, we prove
the result when $G$ is a totally split torus. Here the basic tool is the
result of Thomason, which gives a generic description of the action of a
torus on a noetherian separated algebraic space, and we prove the result by
noetherian induction, using the localization sequence for the $K$-theory of
equivariant coherent sheaves. Like in {\bf \cite{Vi2}}, the difficulty here
is that the homomorphism is defined via pullbacks, and thus it does not
commute with the pushforwards intervening in the localization sequence. This
is solved by producing a different isomorphism between the two groups in
questions, using pushforwards instead of pullbacks, and then relating this
to our map, via the self-intersection formula.

The next step is to prove the result in case $G={\rm GL}_n$; here the key
point is a result of Merkurjev which links the equivariant $K$-theory of a
scheme with a ${\rm GL}_n$-action to the equivariant $K$-theory of the
action of a maximal torus. This is carried out in Section 4.

Finally (Section 5) we reduce the general result to the case of ${\rm GL}_n$%
, by considering an embedding $G\subseteq {\rm GL}_n$, and the induced
action of ${\rm GL}_n$ on $Y={\rm GL}_n\times ^GX$. It is at this point that
condition (c) enters, allowing a clear description of $Y^\sigma $ where $%
\sigma $ is an essential dual cyclic subgroup of $G$ (Proposition \ref
{disjointunion}).

\smallskip\ 

{\bf Acknowledgments}. We wish to thank the referee for her or his useful
and precise remarks. We also thank Bertrand Toen who pointed out the content
of Remark \ref{bert} to us.

\section{General constructions}

{\it Notation}. If $S$ is a separated noetherian scheme, $X$ is a noetherian
separated $S$-algebraic space (which will be most of the time assumed to be
regular) and $G$ is a flat affine group scheme over $S$ acting on $X$, we
denote by $K_{*}(X,G)$ (respectively, $K^{\prime }(X,G)$) the Waldhausen $K$%
-theory of the complicial biWaldhausen (\cite{Th-Tr}) category ${\cal W}%
_{1,X,G}$ of complexes of quasicoherent $G$-equivariant ${\cal O}_X$-Modules
with bounded coherent cohomology (respectively, the Quillen $K$-theory of $G$%
-equivariant coherent ${\cal O}_X$-Modules).

As shown in the Appendix, if $X$ is regular $K_{*}(X,G)$ is isomorphic to $%
K_{*}^{\prime }(X,G)$ and has a canonical graded ring structure. When $X$ is
regular, the isomorphism $K_{*}(X,G)\simeq K_{*}^{\prime }(X,G)$ will then
allow us to switch between the two theories when needed.

\smallskip\ 

\subsection{Morphisms of actions and induced maps on $K$-theory}

Let $S$ be a scheme. By an {\em action} over $S$ we mean a triple $\left(
X,G,\rho \right) $ where $X$ is an $S$-algebraic space, $G$ is a group
scheme over $S$ and $\rho :G\times _SX\rightarrow X$ is an action of $G$ on $%
X$ over $S$. A morphism of actions 
\[
\left( f,\phi \right) :\left( X,G,\rho \right) \longrightarrow \left(
X^{\prime },G^{\prime },\rho ^{\prime }\right) 
\]
is a pair of $S$-morphisms $f:X\rightarrow X^{\prime }$ and $\phi
:G\rightarrow G^{\prime }$, where $\phi $ is a morphism of $S$-group
schemes, such that the following diagram commutes 
\[
\begin{tabular}{lll}
$G\times _SX$ & $\stackrel{\rho }{\longrightarrow }$ & $X$ \\ 
$^{\phi \times f}\downarrow $ &  & $\downarrow ^f$ \\ 
$G^{\prime }\times _SX^{\prime }$ & $\stackunder{\rho ^{\prime }}{%
\longrightarrow }$ & $X^{\prime }$%
\end{tabular}
\]
Equivalently, $f$ is required to be $G$-equivariant with respect to the
given $G$-action on $X$ and the $G$-action on $X^{\prime }$ induced by
composition with $\phi $.

A morphism of actions $\left( f,\phi \right) :\left( X,G,\rho \right)
\rightarrow \left( X^{\prime },G^{\prime },\rho ^{\prime }\right) $ induces
an exact functor $\left( f,\phi \right) ^{*}:{\cal W}_{3,X^{\prime
},G^{\prime }}\rightarrow {\cal W}_{3,X,G}$, where ${\cal W}_{3,Y,H}$
denotes the complicial biWaldhausen category of complexes of $H$-equivariant
flat quasi-coherent Modules with bounded coherent cohomology on the $H$%
-algebraic space $Y$ (see Appendix). Let $({\cal E}^{*},\varepsilon ^{*})$
be an object of ${\cal W}_{3,X^{\prime },G^{\prime }}$, i.e. ${\cal E}^{*}$
is a complex of $G^{\prime }$-equivariant flat quasi-coherent ${\cal O}%
_{X^{\prime }}$-Modules with bounded coherent cohomology and for any $i$ 
\[
\varepsilon ^i:{\rm pr}_2^{\prime *}{\cal E}^i\widetilde{\longrightarrow }%
\rho ^{\prime *}{\cal E}^i 
\]
is an isomorphism satisfying the usual cocycle condition. Here, ${\rm pr}%
_2^{\prime }:G^{\prime }\times _SX^{\prime }\rightarrow X^{\prime }$ denotes
the obvious projection and similarly for ${\rm pr}_2:G\times _SX\rightarrow
X $. Since 
\[
\rho ^{*}f^{*}{\cal E}^{*}=\left( f\rho \right) ^{*}{\cal E}^{*}=\left( \phi
\times f\right) ^{*}\rho ^{\prime *}{\cal E}^{*} 
\]
and 
\[
{\rm pr}_2^{*}f^{*}{\cal E}^{*}{\bf =}\left( \phi \times f\right) ^{*}{\rm pr%
}_2^{\prime *}{\cal E}^{*} 
\]
we define $\left( f,\phi \right) ^{*}({\cal E}^{*},\varepsilon ^{*})\doteq
(f^{*}{\cal E}^{*},\left( \phi \times f\right) ^{*}\left( \varepsilon
^{*}\right) )$ $\in {\cal W}_{3,X,G}$ (the cocycle condition for each $%
\left( \phi \times f\right) ^{*}\left( \varepsilon ^i\right) $ following
from the same condition for $\varepsilon ^i$); $\left( f,\phi \right) ^{*}$
is then defined on morphisms in the only natural way. $\left( f,\phi \right)
^{*}$ is an exact functor and, if $X$ and $X^{\prime }$ are regular so that
the Waldhausen $K$-theory of ${\cal W}_{3,X,G}$ (respectively, of ${\cal W}%
_{3,X^{\prime },G^{\prime }}$) coincides with $K_{*}\left( X,G\right) $
(resp., $K_{*}\left( X^{\prime },G^{\prime }\right) $) (see Appendix), it
defines a ring morphism 
\[
\left( f,\phi \right) ^{*}:K_{*}\left( X^{\prime },G^{\prime }\right)
\longrightarrow K_{*}\left( X,G\right) . 
\]
A similar argument shows that if $f$ is flat, $\left( f,\phi \right) $
induces a morphism 
\[
\left( f,\phi \right) ^{*}:K_{*}^{\prime }\left( X^{\prime },G^{\prime
}\right) \longrightarrow K_{*}^{\prime }\left( X,G\right) . 
\]

\begin{example}
\label{action}Let $G$ and $H$ be group schemes over $S$ and $X$ be an $S$%
-algebraic space. Suppose moreover:

\begin{enumerate}
\item  $G$ and $H$ act on $X$:

\item  $G$ acts on $H$ by $S$-group schemes automorphisms (i.e., it is given
a morphism $G\rightarrow {\rm Aut}_{\left( {\tt GrSch}\right) _{/S}}\left(
H\right) $ of group functors over $S$);

\item  The two preceding actions are compatible, i.e. for any $S$-scheme $T$%
, any $g\in G\left( T\right) $, $h\in H\left( T\right) $ and $x\in X\left(
T\right) $ we have 
\[
g\cdot \left( h\cdot x\right) =h^g\cdot \left( g\cdot x\right) 
\]
where $\left( g,h\right) \mapsto $ $h^g$ denotes the action of $G\left(
T\right) $ on $H\left( T\right) $.
\end{enumerate}

If $g\in G\left( S\right) $ and $g_T$ denotes its image via $G\left(
S\right) \rightarrow G\left( T\right) $, let us define a morphism of actions 
$\left( f_g,\phi _g\right) :\left( X,H\right) \rightarrow \left( X,H\right) $
as: 
\[
f_g\left( T\right) :X\left( T\right) \longrightarrow X\left( T\right)
:x\longmapsto g_T\cdot x
\]
\[
\phi _g\left( T\right) :H\left( T\right) \longrightarrow H\left( T\right)
:h\longmapsto h^{g_T}
\]
This is an isomorphism of actions and induces an action of the group $%
G\left( S\right) $ on $K_{*}^{\prime }\left( X,H\right) $ and on $K_{*}(X,H)$%
. This applies, in particular, to the case where $X$ is an algebraic space
with a $G$ action and $G\vartriangleright H$, $G$ acting on $H$ by
conjugation.
\end{example}

\subsection{The basic definitions and results}

Let $G$ be a linear algebraic $k$-group scheme $G$ acting with finite
geometric stabilizers on a regular noetherian separated algebraic space $X$
over $k$.

We denote by ${\rm \QTR{mathrm}{R}}\left( G\right) $ the representation ring
of $G$.

A (Cartier) {\em dual cyclic subgroup} of $G$ over $k$ is a $k$-subgroup
scheme $\sigma \subseteq G$ such that there exists an $n>0$ and an
isomorphism of $k$-groups $\sigma \simeq {\bf \mu }_{n,k}$. If $\sigma ,\rho
\ $are dual cyclic subgroups of $G$ and $L$ is an extension of $k$, we say
that $\sigma $ and $\rho $ are conjugate over $L$ if there exists $g\in
G\left( L\right) $ such that $g\sigma _{\left( L\right) }g^{-1}=\rho
_{\left( L\right) }$ (where $H_{\left( L\right) }\doteq H\times _{{\rm Spec}%
k}{\rm Spec}L$, for any $k$-group scheme $H$) as $L$-subgroup schemes of $%
G_{\left( L\right) }$.

A dual cyclic subgroup $\sigma \subseteq G$ is said to be {\em essential} if 
$X^\sigma \neq \emptyset $.

We say that the action of $G$ on $X$ is {\em sufficiently rational} if:

\begin{enumerate}
\item  any two essential dual cyclic subgroups of $G$ are conjugated over $k$
iff they are conjugated over an algebraic closure $\overline{k}$ of $k$;

\item  any essential dual cyclic subgroup $\overline{\rho }$ of $G_{\left( 
\overline{k}\right) }$ is conjugated over $\overline{k}$ to a dual cyclic
subgroup of the form $\sigma _{\left( \overline{k}\right) }$ where $\sigma
\subseteq G$ is (essential) dual cyclic.
\end{enumerate}

We denote by ${\cal \QTR{mathcal}{C}}\left( G\right) $ a set of
representatives for essential dual cyclic subgroups of $G$ with respect to
the relation of conjugacy over $k$.

\begin{remark}
\label{rmkRC}Note that if the action is sufficiently rational and $\rho $, $%
\sigma $ are essential dual cyclic subgroups of $G$ which are conjugate over
an algebraically closed extension $\Omega $ of $k$, then they are also
conjugate over $k$.
\end{remark}

\begin{proposition}
\label{angelo}Any action of ${\rm GL}_n$, ${\rm SL}_n,{\rm Sp}_{2n}$ or of a
split torus is sufficiently rational.
\end{proposition}

\TeXButton{Proof}{\proof}If $G$ is a split torus, Condition~1 is clear,
because $G$ is abelian, while it follows from the rigidity of diagonalizable
groups that any subgroup scheme of $G_{\overline{k}}$ is in fact defined
over $k$. Let $\sigma \subseteq {\rm GL}_m$ be a dual cyclic subgroup. Since 
$\sigma $ is diagonalizable, we have an eigenspace decomposition

\[
V=k^m=\bigoplus_{\chi \in \widehat{\sigma }}V_\chi ^\sigma 
\]
such that the $\chi $ with $V_\chi \neq 0$ generate $\widehat{\sigma }$.
Conversely, given a cyclic group $C$ and a decomposition 
\[
V=\bigoplus_{\chi \in \widehat{C}}V_\chi 
\]
such that the $\chi $ with $V_\chi \neq 0$ generate $\widehat{C}$, there is
a corresponding embedding of the Cartier dual $\sigma $ of $C$ into ${\rm GL}%
_n$ with $V_\chi =V_\chi ^\sigma $ for each $\chi \in C=\widehat{\sigma }$.
Now, if $\sigma \subseteq {\rm GL}_{m,\overline{k}}$ is a dual cyclic
subgroup defined over $\overline{k}$, we can apply an element of ${\rm GL}_m(%
\overline{k})$ to make the $V_\chi $ defined over $k$, and then $g\sigma
g^{-1}$ will be defined over $k$. If $\sigma \subseteq {\rm GL}_m$ and $\tau
\subseteq {\rm GL}_m$ are dual cyclic subgroups which are conjugate over $%
\overline{k}$, pick an element of ${\rm GL}_m(\overline{k})$ sending $\sigma 
$ to $\tau $. This induces an isomorphism $\phi :\sigma _{\overline{k}%
}\simeq \tau _{\overline{k}}$, which by rigidity will be defined over $k$.
Then if $\chi $ and $\chi ^{\prime }$ are characters which correspond under
the isomorphism of $\widehat{\sigma }$ and $\widehat{\tau }$ induced by $%
\phi $, then the dimension of $V_\chi ^\sigma $ is equal to the dimension of 
$V_{\chi ^{\prime }}^\tau $, so we can find an element $g$ of ${\rm GL}_m$
which carries each $V_\chi ^\sigma $ onto the corresponding $V_{\chi \prime
}^\tau $; conjugation by this element carries $\sigma $ onto $\tau $. For $%
{\rm SL}_m$ the proof is very similar, if we remark that to give a dual
cyclic subgroup $\sigma \subseteq {\rm SL}_m\subseteq {\rm GL}_m$ correspond
to giving a decomposition 
\[
V=k^m=\bigoplus_{\chi \in \widehat{\sigma }}V_\chi ^\sigma 
\]
such that the $\chi $ with $V_\chi ^\sigma \neq 0$ generate $\widehat{\sigma 
}$, with the condition $\prod_{\chi \in \widehat{\sigma }}\chi ^{\dim V_\chi
^\sigma }=1\in \widehat{\sigma }$. For ${\rm Sp}_m\subseteq {\rm GL}_{2m}$,
a dual cyclic subgroup $\sigma \subseteq {\rm Sp}_m$ gives a decomposition 
\[
V=k^{2m}=\bigoplus_{\chi \in \widehat{\sigma }}V_\chi ^\sigma 
\]
such that the $\chi $ with $V_\chi ^\sigma \neq 0$ generate $\widehat{\sigma 
}$, with the condition that for $v\in V_\chi ^\sigma $ and $v^{\prime }\in
V_{\chi \prime }^\sigma $ the symplectic product of $v$ and $v^{\prime }$ is
always $0$, unless $\chi \chi ^{\prime }=1\in \widehat{\sigma }$. Both
conditions then follow rather easily from the fact that any two symplectic
forms over a vector space are isomorphic. \TeXButton{End Proof}{\endproof}

\smallskip\ 

Let $N_{\left( G,X\right) }$ denote the least common multiple of the orders
of essential dual cyclic subgroups of $G$. Notice that $N_{(G,X)}$ is
finite: since the action has finite stabilizers, the group scheme of
stabilizers is quasifinite over $X$, therefore the orders of the stabilizers
of the geometric points of $X$ are globally bounded.

We define $\Lambda _{\left( G,X\right) }\doteq {\Bbb Z}\left[ 1/N_{\left(
G,X\right) }\right] $. If $H\subseteq G$ is finite, we also write $\Lambda
_H $ for ${\Bbb Z}\left[ 1/\left| H\right| \right] $. Note that, if $\sigma
\subseteq G$ is dual cyclic, then $\Lambda _\sigma =\Lambda _{\left( \sigma ,%
{\rm Spec}k\right) }$ and if moreover $\sigma $ is essential $\Lambda
_\sigma \subseteq \Lambda _{\left( G,X\right) }$.

If $H\subseteq G$ is a subgroup scheme and $A$ is a ring, we write ${\rm 
\QTR{mathrm}{R}}\left( H\right) _A$ for ${\rm \QTR{mathrm}{R}}\left(
H\right) \otimes _{{\Bbb Z}}A$. We denote by {\rm rk}$_H:{\rm \QTR{mathrm}{R}%
}\left( H\right) \longrightarrow {\Bbb Z}$ and by {\rm rk}$_{H,\Lambda
_{\left( G,X\right) }}:{\rm \QTR{mathrm}{R}}\left( H\right) _{\Lambda
_{\left( G,X\right) }}\longrightarrow \Lambda _{\left( G,X\right) } $ the
rank ring homomorphisms.

We let $K_{*}^{\prime }\left( X,G\right) _{\Lambda _{\left( G,X\right)
}}=K_{*}^{\prime }\left( X,G\right) \otimes \Lambda _{G,X}$ and $K_{*}\left(
X,G\right) _{\Lambda _{\left( G,X\right) }}=K_{*}\left( X,G\right) \otimes
\Lambda _{G,X}$ (for the notation, see the beginning of this Section).
Recall that $K_{*}\left( X,G\right) _{\Lambda _{\left( G,X\right) }}$ is an $%
{\rm R}\left( G\right) $-algebra via the pullback ${\rm \QTR{mathrm}{R}}%
\left( G\right) \simeq K_0\left( {\rm Spec}k,G\right) \rightarrow K_0\left(
X,G\right) $ and that $K_{*}\left( X,G\right) \simeq K_{*}^{\prime }\left(
X,G\right) $ since $X$ is regular (Appendix).

If $\sigma $ is a dual cyclic subgroup of $G$ of order $n$, the choice of a
generator $t$ for the dual group $\widehat{\sigma }\doteq Hom_{{\tt GrSch}%
_{/k}}\left( \sigma ,{\Bbb G}_{m,k}\right) $ determines an isomorphism 
\[
{\rm \QTR{mathrm}{R}}\left( \sigma \right) \simeq \frac{{\Bbb Z}\left[
t\right] }{\left( t^n-1\right) }. 
\]
Let $p_\sigma $ be the canonical ring surjection 
\[
\frac{{\Bbb Z}\left[ t\right] }{\left( t^n-1\right) }\longrightarrow
\prod_{d\mid n}\frac{{\Bbb Z}\left[ t\right] }{\left( \Phi _d\right) } 
\]
and $\widetilde{p_\sigma }$ the induced surjection 
\[
\frac{{\Bbb Z}\left[ t\right] }{\left( t^n-1\right) }\longrightarrow \frac{%
{\Bbb Z}\left[ t\right] }{\left( \Phi _n\right) }, 
\]
where $\Phi _d$ is the $d$-th cyclotomic polynomial. If ${\frak m}_\sigma $
is the kernel of the composition 
\[
{\rm \QTR{mathrm}{R}}\left( \sigma \right) \simeq \frac{{\Bbb Z}\left[
t\right] }{\left( t^n-1\right) }\longrightarrow \frac{{\Bbb Z}\left[
t\right] }{\left( \Phi _n\right) }, 
\]
the quotient ring ${\rm \QTR{mathrm}{R}}\left( \sigma \right) /{\frak m}%
_\sigma $ does not depend on the choice of the generator $t$ for $\widehat{%
\sigma }$.

{\bf Notation.} We denote by $\widetilde{{\rm R}}\left( \sigma \right) $ the
quotient ${\rm \QTR{mathrm}{R}}\left( \sigma \right) /{\frak m}_\sigma $.

We remark that if $\sigma $ is dual cyclic of order $n$ and $t$ is a
generator of $\widehat{\sigma }$, there are isomorphisms 
\begin{equation}
{\rm R}\left( \sigma \right) _{\Lambda _\sigma }\simeq \frac{\Lambda _\sigma
\left[ t\right] }{\left( t^n-1\right) }\simeq \prod_{d\mid n}\frac{\Lambda
_\sigma \left[ t\right] }{\left( \Phi _d\right) }.  \label{decomposition}
\end{equation}

Let $\widetilde{\pi _\sigma }:{\rm \QTR{mathrm}{R}}\left( G\right) _{\Lambda
_{\left( G,X\right) }}\rightarrow \widetilde{{\rm R}}\left( \sigma \right)
_{\Lambda _{\left( G,X\right) }}$ be the canonical ring homomorphism. The $%
\sigma ${\em -localization} $K_{*}^{\prime }\left( X,G\right) _\sigma $ of $%
K_{*}^{\prime }\left( X,G\right) _{\Lambda _{\left( G,X\right) }}$ is the
localization of the ${\rm \QTR{mathrm}{R}}\left( G\right) _{\Lambda _{\left(
G,X\right) }}$-module $K_{*}^{\prime }\left( X,G\right) _{\Lambda _{\left(
G,X\right) }}$ at the multiplicative subset $S_\sigma \doteq \widetilde{\pi
_\sigma }^{-1}\left( 1\right) $. The $\sigma $-localization $K_{*}\left(
X,G\right) _\sigma $ is defined in the same way. If $H\subseteq G$ is a
subgroup scheme, we also write ${\rm R}\left( H\right) _\sigma $ for $%
S_\sigma ^{-1}\left( {\rm \QTR{mathrm}{R}}\left( H\right) _{\Lambda _{\left(
G,X\right) }}\right) $.

If $\sigma $ {\it is the trivial group}, we denote by $K_{*}^{\prime }\left(
X,G\right) _{{\rm geom}}$ the $\sigma $-localization of $K_{*}^{\prime
}\left( X,G\right) _{\Lambda _{\left( G,X\right) }}$ and call it the {\em %
geometric part} or {\em geometric localization} of $K_{*}^{\prime }\left(
X,G\right) _{\Lambda _{\left( G,X\right) }}$. Note that $\widetilde{\pi _1}$
coincides with the rank morphism ${\rm rk}_{G,\Lambda _{\left( G,X\right) }}:%
{\rm \QTR{mathrm}{R}}\left( G\right) _{\Lambda _{\left( G,X\right)
}}\longrightarrow \Lambda _{\left( G,X\right) }$. Same definition for $%
K_{*}\left( X,G\right) _{{\rm geom}}$.

Let ${\rm N}_G\left( \sigma \right) $ (respectively, ${\rm C}_G\left( \sigma
\right) \subseteq {\rm N}_G\left( \sigma \right) $) be the normalizer
(resp., the centralizer) of $\sigma $ in $G$; since ${\rm Aut}\left( \sigma
\right) $ is a finite constant group scheme, 
\[
{\rm W}_G\left( \sigma \right) \doteq \frac{{\rm N}_G\left( \sigma \right) }{%
{\rm C}_G\left( \sigma \right) } 
\]
is also a constant group scheme over $k$ associated to a finite group $%
w_G\left( \sigma \right) $.

\begin{lemma}
\label{lemma}Let $H$ be a $k$-linear algebraic group, $\sigma \simeq \mu
_{n,k}$ a normal subgroup and $Y$ an algebraic space with an action of $%
H/\sigma $. Then there is a canonical action of $w_H\left( \sigma \right) $
on $K_{*}^{\prime }\left( Y,{\rm C}_H\left( \sigma \right) \right) $.
\end{lemma}

\TeXButton{Proof}{\proof}Let us first assume that the natural map 
\begin{equation}
H\left( k\right) \longrightarrow w_H\left( \sigma \right)  \label{surj}
\end{equation}
is surjective (which is true, for example, if $k$ is algebraically closed).
Since ${\rm C}_H\left( \sigma \right) \left( k\right) $ acts trivially on $%
K_{*}^{\prime }\left( Y,{\rm C}_H\left( \sigma \right) \right) $ and, by
Example \ref{action}, $H\left( k\right) $ acts naturally on $K_{*}^{\prime
}\left( Y,{\rm C}_H\left( \sigma \right) \right) $, we may use (\ref{surj})
to define the desired action.

In general (\ref{surj}) is not surjective and we proceed as follows. Suppose
we can find a closed immersion of $k$-linear algebraic groups $%
H\hookrightarrow H^{\prime }$ such that

\begin{enumerate}
\item[(i)]  $\sigma $ is normal in $H^{\prime }$;

\item[(ii)]  $H^{\prime }/{\rm C}_{H^{\prime }}\left( \sigma \right) \simeq
W_H\left( \sigma \right) $;

\item[(iii)]  $H^{\prime }\left( k\right) \rightarrow w_H\left( \sigma
\right) $ is surjective.
\end{enumerate}

Consider the open and closed immersion $Y\times {\rm C}_H\left( \sigma
\right) \hookrightarrow Y\times H$; this induces an open and closed
immersion $Y\times ^{{\rm C}_H\left( \sigma \right) }{\rm C}_{H^{\prime
}}\left( \sigma \right) \hookrightarrow Y\times ^{{\rm C}_H\left( \sigma
\right) }H^{\prime }$ whose composition with the \'etale cover $Y\times ^{%
{\rm C}_H\left( \sigma \right) }H^{\prime }\rightarrow Y\times ^HH^{\prime }$
is easily checked (e.g. on geometric points) to be an isomorphism. Therefore 
\[
K_{*}^{\prime }\left( Y\times ^HH^{\prime },{\rm C}_{H^{\prime }}\left(
\sigma \right) \right) \simeq K_{*}^{\prime }\left( Y\times ^{{\rm C}%
_H\left( \sigma \right) }{\rm C}_{H^{\prime }}\left( \sigma \right) ,{\rm C}%
_{H^{\prime }}\left( \sigma \right) \right) \simeq K_{*}^{\prime }\left( Y,%
{\rm C}_H\left( \sigma \right) \right) 
\]
where the last isomorphism is given by Morita equivalence theorem ({\bf \cite
{Th3}}, Proposition 6.2). By (i) and (iii) we can apply the argument at the
beginning of the proof and get an action of $w_H\left( \sigma \right) $ on $%
K_{*}^{\prime }\left( Y\times ^HH^{\prime },{\rm C}_{H^{\prime }}\left(
\sigma \right) \right) $ and therefore on $K_{*}^{\prime }\left( Y,{\rm C}%
_H\left( \sigma \right) \right) $ as desired. It is not difficult to check
that this action does not depend on the chosen immersion $H\hookrightarrow
H^{\prime }$.

Finally, let us prove that there exists a closed immersion $H\hookrightarrow
H^{\prime }$ satisfying conditions (i)-(iii) above. First choose a closed
immersion $j:H\hookrightarrow {\rm GL}_{n,k}$ for some $n$. Clearly 
\[
H/{\rm C}_H\left( \sigma \right) \hookrightarrow {\rm GL}_{n,k}/{\rm C}_{%
{\rm GL}_{n,k}}\left( \sigma \right) 
\]
and, embedding $\sigma $ in a maximal torus of ${\rm GL}_{n,k}$, it is easy
to check that 
\[
{\rm GL}_{n,k}\left( k\right) \rightarrow \left( {\rm GL}_{n,k}/{\rm C}_{%
{\rm GL}_{n,k}}\left( \sigma \right) \right) \left( k\right) 
\]
is surjective. Now define $H^{\prime }$ as the inverse image of $H/{\rm C}%
_H\left( \sigma \right) $ in the normalizer ${\rm N}_{{\rm GL}_{n,k}}\left(
\sigma \right) $. \TeXButton{End Proof}{\endproof}

\begin{corollary}
\label{azione}There is a canonical action of $w_G\left( \sigma \right) $ on $%
K_{*}^{\prime }\left( X^\sigma ,{\rm C}_G\left( \sigma \right) \right) $
which induces an action on $K_{*}^{\prime }\left( X^\sigma ,{\rm C}_G\left(
\sigma \right) \right) _{{\rm geom}}.$
\end{corollary}

\TeXButton{Proof}{\proof}Since ${\rm C}_G\left( \sigma \right) ={\rm C}%
_{N_G\left( \sigma \right) }\left( \sigma \right) $, lemma \ref{lemma}
applied to $Y=X^\sigma $ (respectively, $Y={\rm Spec}k$) and $H=N_G\left(
\sigma \right) $, yields an action of $w_G\left( \sigma \right) $ on $%
K_{*}^{\prime }\left( X^\sigma ,{\rm C}_G\left( \sigma \right) \right) $
(respectively, on $K_0\left( {\rm Spec}k,{\rm C}_G\left( \sigma \right)
\right) ={\rm R}\left( {\rm C}_G\left( \sigma \right) \right) $). The
multiplicative system $S_1=${\rm rk}$^{-1}\left( 1\right) $ is preserved by
this action so that there is an induced action on the ring $S_1^{-1}{\rm R}%
\left( {\rm C}_G\left( \sigma \right) \right) $. The pullback 
\[
K_0\left( {\rm Spec}k,{\rm C}_G\left( \sigma \right) \right) \rightarrow
K_0\left( X^\sigma ,{\rm C}_G\left( \sigma \right) \right) 
\]
is $w_G\left( \sigma \right) $-equivariant and then $w_G\left( \sigma
\right) $ acts on $K_{*}^{\prime }\left( X^\sigma ,{\rm C}_G\left( \sigma
\right) \right) _{{\rm geom}}$. \TeXButton{End Proof}{\endproof}

\begin{remark}
\label{python}If $Y$ is regular, Lemma \ref{lemma} gives also an action of $%
w_H\left( \sigma \right) $ on $K_{*}\left( Y,{\rm C}_H\left( \sigma \right)
\right) $, since $K_{*}\left( Y,{\rm C}_H\left( \sigma \right) \right)
\simeq K_{*}^{\prime }\left( Y,{\rm C}_H\left( \sigma \right) \right) $
(Appendix). In particular, since by {\bf \cite{Th1}}, Prop. 3.1, $X^\sigma $
is regular, Corollary \ref{azione} still holds for $K_{*}\left( X^\sigma ,%
{\rm C}_G\left( \sigma \right) \right) _{{\rm geom}}$.
\end{remark}

Note also that the embedding of $k$-group schemes $W_G\left( \sigma \right)
\hookrightarrow {\rm Aut}_k\left( \sigma \right) $ induces, by Example \ref
{action}, an action of $w_G\left( \sigma \right) $ on $K_0\left( {\rm Spec}%
k,\sigma \right) ={\rm \QTR{mathrm}{R}}\left( \sigma \right) $.

The product in $\sigma $ induces a morphism of $k$-groups 
\[
\sigma \times {\rm C}_G\left( \sigma \right) \longrightarrow {\rm C}_G\left(
\sigma \right) 
\]
which in its turn induces a morphism 
\[
m_\sigma :K_{*}\left( X^\sigma ,{\rm C}_G\left( \sigma \right) \right)
\longrightarrow K_{*}\left( X^\sigma ,\sigma \times {\rm C}_G\left( \sigma
\right) \right) . 
\]

\begin{lemma}
If $H\subseteq G$ is a subgroup scheme and $\sigma $ is contained in the
center of $H$, \label{sga}there is a canonical ring isomorphism 
\[
K_{*}\left( X^\sigma ,\sigma \times H\right) \simeq K_{*}\left( X^\sigma
,H\right) \otimes {\rm \QTR{mathrm}{R}}\left( \sigma \right) .
\]
\end{lemma}

\TeXButton{Proof}{\proof}Since $\sigma $ acts trivially on $X^\sigma $, we
have an equivalence ({\bf \cite{SGA3}}, I 4.7.3) 
\begin{equation}
\left( \sigma \times H\right) {\tt -Coh}_{X^\sigma }\simeq \bigoplus_{%
\widehat{\sigma }}\left( H{\tt -Coh}_{X^\sigma }\right)  \label{splitting}
\end{equation}
(where $\widehat{\sigma }$ is the character group of $\sigma $) which
induces an isomorphism 
\[
K_{*}^{\prime }\left( X^\sigma ,\sigma \times H\right) \simeq K_{*}^{\prime
}\left( X^\sigma ,H\right) \otimes {\rm R}\left( \sigma \right) . 
\]
and we conclude since $K_{*}\left( Y,H\right) \simeq K_{*}^{\prime }\left(
Y,H\right) $ and $K_{*}\left( X^\sigma ,\sigma \times H\right) \simeq
K_{*}^{\prime }\left( X^\sigma ,\sigma \times H\right) $ (Appendix). 
\TeXButton{End Proof}{\endproof}

For any essential dual cyclic subgroup $\sigma \subseteq G$, let $\Lambda
\doteq \Lambda _{\left( G,X\right) }$ and consider the composition 
\[
K_{*}\left( X,G\right) _\Lambda \rightarrow K_{*}\left( X,{\rm C}_G\left(
\sigma \right) \right) _\Lambda \rightarrow K_{*}\left( X^\sigma ,{\rm C}%
_G\left( \sigma \right) \right) _\Lambda \stackrel{m_\sigma }{%
\longrightarrow } 
\]
$\ $%
\begin{equation}
\longrightarrow K_{*}\left( X^\sigma ,{\rm C}_G\left( \sigma \right) \right)
_\Lambda \otimes _\Lambda {\rm R}\left( \sigma \right) _\Lambda
\longrightarrow K_{*}\left( X^\sigma ,{\rm C}_G\left( \sigma \right) \right)
_{{\rm geom}}\otimes _\Lambda \widetilde{{\rm R}}\left( \sigma \right)
_\Lambda  \label{uno}
\end{equation}
where the first map is induced by group restriction the last one is the
geometric localization map tensored with the projection ${\rm R}\left(
\sigma \right) _\Lambda \rightarrow \widetilde{{\rm R}}\left( \sigma \right)
_\Lambda $ and we have used Lemma \ref{sga} with $H={\rm C}_G\left( \sigma
\right) $; the second map is induced by restriction along $X^\sigma
\hookrightarrow X$ which is a regular closed immersion ({\bf \cite{Th1}},
Prop. 3.1) therefore has finite $Tor$-dimension so that the pullback on $K$%
-groups is well defined (see the Appendix). It is not difficult to show that
the image of (\ref{uno}) is actually contained in the invariant submodule 
\[
\left( K_{*}\left( X^\sigma ,{\rm C}_G\left( \sigma \right) \right) _{{\rm %
geom}}\otimes _\Lambda \widetilde{{\rm R}}\left( \sigma \right) _\Lambda
\right) ^{w_G\left( \sigma \right) } 
\]
so that we get a map 
\[
\psi _{\sigma ,X}:K_{*}\left( X,G\right) _\Lambda \longrightarrow \left(
K_{*}\left( X^\sigma ,{\rm C}_G\left( \sigma \right) \right) _{{\rm geom}%
}\otimes _\Lambda \widetilde{{\rm R}}\left( \sigma \right) _\Lambda \right)
^{w_G\left( \sigma \right) }. 
\]
Our basic map will be: 
\begin{equation}
\Psi _{X,G}\doteq \prod_{\sigma \in {\cal C}\left( G\right) }\psi _{\sigma
,X}:K_{*}\left( X,G\right) _\Lambda \longrightarrow \prod_{\sigma \in {\cal C%
}\left( G\right) }\left( K_{*}\left( X^\sigma ,{\rm C}_G\left( \sigma
\right) \right) _{{\rm geom}}\otimes _\Lambda \widetilde{{\rm R}}\left(
\sigma \right) _\Lambda \right) ^{w_G\left( \sigma \right) }.  \label{due}
\end{equation}
Note that $\Psi _{X,G}$ is a morphism of ${\rm R}(G)$-algebras as a
composition of morphisms of $R(G)$-algebras.

The following technical lemma will be used in Propositions \ref{end'} and 
\ref{pincopallino}.

\begin{lemma}
\label{construction}Let $G$ be a linear algebraic $k$-group acting with
finite stabilizers on a noetherian separated $k$-algebraic space $X$ and $%
\Lambda \doteq \Lambda _{\left( G,X\right) }$. Let $H\subseteq G$ be a
subgroup and $\sigma $ an essential dual cyclic subgroup contained in the
center of $H$. Consider the composition 
\begin{equation}
K_{*}^{\prime }\left( Y^\sigma ,H\right) _\Lambda \longrightarrow
K_{*}^{^{\prime }}\left( Y^\sigma ,H\right) _\Lambda \otimes _\Lambda {\rm R}%
\left( \sigma \right) _\Lambda \longrightarrow K_{*}^{^{\prime }}\left(
Y^\sigma ,H\right) _{{\rm geom}}\otimes _\Lambda \widetilde{{\rm R}}\left(
\sigma \right) _\Lambda   \label{preloc}
\end{equation}
where the first morphism is induced by the product morphism $\sigma \times
H\rightarrow H$ (recall Lemma \ref{sga}) and the second is the tensor
product of the geometric localization morphism with the projection ${\rm R}%
\left( \sigma \right) _\Lambda \rightarrow \widetilde{{\rm R}}\left( \sigma
\right) _\Lambda $. Then, (\ref{preloc}) factors through $K_{*}^{\prime
}\left( Y^\sigma ,H\right) _\Lambda \rightarrow K_{*}^{\prime }\left(
Y^\sigma ,H\right) _\sigma $, yielding a morphism 
\begin{equation}
\theta _{H,\sigma }:K_{*}^{\prime }\left( Y^\sigma ,H\right) _\sigma
\longrightarrow K_{*}^{^{\prime }}\left( Y^\sigma ,H\right) _{{\rm geom}%
}\otimes _\Lambda \widetilde{{\rm R}}\left( \sigma \right) _\Lambda \text{.}
\label{LOC}
\end{equation}
\end{lemma}

\TeXButton{Proof}{\proof} Let $S_1$ (resp. $S_\sigma $) be the
multiplicative subset in ${\rm R}\left( H\right) _\Lambda $ consisting of
elements going to $1$ via the homomorphism $rk_{H,\Lambda }:{\rm R}\left(
H\right) _\Lambda \rightarrow \Lambda $ (resp., ${\rm R}\left( H\right)
_\Lambda \rightarrow \widetilde{{\rm R}}\left( \sigma \right) _\Lambda $).
Observe that $K_{*}^{^{\prime }}\left( X^\sigma ,H\right) _\Lambda \otimes
_\Lambda {\rm R}\left( \sigma \right) _\Lambda $ (resp. $K_{*}^{^{\prime
}}\left( X^\sigma ,H\right) _{{\rm geom}}\otimes _\Lambda \widetilde{{\rm R}}%
\left( \sigma \right) _\Lambda $) is canonically an ${\rm R}\left( H\right)
_\Lambda \otimes {\rm R}\left( \sigma \right) _\Lambda $-module (resp. an $%
S_1^{-1}{\rm R}\left( H\right) _\Lambda \otimes \widetilde{{\rm R}}\left(
\sigma \right) _\Lambda $-module) and therefore an ${\rm R}\left( H\right) $
module via the coproduct ring morphism 
\[
\Delta _\sigma :{\rm R}\left( H\right) _\Lambda \longrightarrow {\rm R}%
\left( H\right) _\Lambda \otimes {\rm R}\left( \sigma \right) _\Lambda 
\]
(resp. via the ring morphism 
\[
f_\sigma :{\rm R}\left( H\right) _\Lambda \stackrel{\Delta _\sigma }{%
\longrightarrow }{\rm R}\left( H\right) _\Lambda \otimes {\rm R}\left(
\sigma \right) _\Lambda \longrightarrow S_1^{-1}{\rm R}\left( H\right)
_\Lambda \otimes \widetilde{{\rm R}}\left( \sigma \right) _\Lambda \text{).} 
\]
If we denote by $A^{\prime }$ the ${\rm R}\left( H\right) _\Lambda $-algebra 
$f_\sigma :{\rm R}\left( H\right) _\Lambda \longrightarrow S_1^{-1}{\rm R}%
\left( H\right) _\Lambda \otimes \widetilde{{\rm R}}\left( \sigma \right)
_\Lambda $, it is enough to show that the localization homomorphism 
\[
l_\sigma ^{\prime }:A^{\prime }\longrightarrow S_\sigma ^{-1}\left(
A^{\prime }\right) 
\]
is an isomorphism, because in this case the morphism (\ref{LOC}) will be
induced by the $S_\sigma $-localization of (\ref{preloc}).

Let $A$ denote the ${\rm R}\left( H\right) _\Lambda $-algebra 
\[
\lambda _1\otimes 1:{\rm R}\left( H\right) _\Lambda \longrightarrow S_1^{-1}%
{\rm R}\left( H\right) _\Lambda \otimes \widetilde{{\rm R}}\left( \sigma
\right) _\Lambda 
\]
where $\lambda _1:{\rm R}\left( H\right) _\Lambda \rightarrow S_1^{-1}{\rm R}%
\left( H\right) _\Lambda $ denotes the localization homomorphism. It is a
well known fact that there is an isomorphism of ${\rm R}\left( H\right)
_\Lambda $-algebras $\varphi :A^{\prime }\rightarrow A$: this is exactly the
dual assertion to ``the action $H\times \sigma \rightarrow \sigma $ is
isomorphic to the projection on the second factor $H\times \sigma
\rightarrow \sigma $ ''. Therefore, we have a commutative diagram 
\[
\begin{tabular}{ccc}
$A^{\prime }$ & $\stackrel{\varphi }{\longrightarrow }$ & $A$ \\ 
$^{l_\sigma ^{\prime }}\downarrow $ &  & $\downarrow ^{l_\sigma }$ \\ 
$S_\sigma ^{-1}A^{\prime }$ & $\stackunder{S_\sigma ^{-1}\varphi }{%
\longrightarrow }$ & $S_\sigma ^{-1}A$%
\end{tabular}
\]
where $l_\sigma $ denotes the localization homomorphism and it is enough to
prove that $l_\sigma $ is an isomorphism. To see this, note that the ring $%
\widetilde{{\rm R}}\left( \sigma \right) _\Lambda $ is a free $\Lambda $%
-module of finite rank (equal to $\phi \left( \left| \sigma \right| \right) $%
, $\phi $ being the Euler function) and there is a norm homomorphism 
\[
{\rm N}:\widetilde{{\rm R}}\left( \sigma \right) _\Lambda \longrightarrow
\Lambda 
\]
sending an element to the determinant of the $\Lambda $-endomorphism of $%
\widetilde{{\rm R}}\left( \sigma \right) _\Lambda $ induced by
multiplication by this element; obviously, we have 
\[
{\rm N}^{-1}\left( \Lambda ^{*}\right) =\left( \widetilde{{\rm R}}\left(
\sigma \right) _\Lambda \right) ^{*}. 
\]
Analogously, there is a norm homomorphism 
\[
{\rm N}^{\prime }:A^{\prime }=S_1^{-1}{\rm R}\left( H\right) _\Lambda
\otimes \widetilde{{\rm R}}\left( \sigma \right) _\Lambda \longrightarrow
S_1^{-1}{\rm R}\left( H\right) _\Lambda 
\]
and 
\[
{\rm N}^{-1}((S_1^{-1}{\rm R}\left( H\right) _\Lambda )^{*})=\left( S_1^{-1}%
{\rm R}\left( H\right) _\Lambda \otimes \widetilde{{\rm R}}\left( \sigma
\right) _\Lambda \right) ^{*}. 
\]
There is a commutative diagram 
\[
\begin{tabular}[t]{ccc}
$S_1^{-1}{\rm R}\left( H\right) _\Lambda \otimes \widetilde{{\rm R}}\left(
\sigma \right) _\Lambda $ & $\stackrel{{\rm N}^{\prime }}{\longrightarrow }$
& $S_1^{-1}{\rm R}\left( H\right) _\Lambda $ \\ 
$^{{\rm rk}_{H,\Lambda }\otimes {\rm id}}\downarrow $ &  & $\quad \downarrow
^{{\rm rk}_{H,\Lambda }}$ \\ 
$\widetilde{{\rm R}}\left( \sigma \right) _\Lambda $ & $\stackunder{N}{%
\longrightarrow }$ & $\Lambda $%
\end{tabular}
\]
and, by definition of $S_1$, we get ${\rm rk}_{H,\Lambda }^{-1}\left(
\Lambda ^{*}\right) =\left( S_1^{-1}{\rm R}\left( H\right) _\Lambda \right)
^{*}$. Therefore, by definition of $S_\sigma $, $S_\sigma /1$ consist of
units in $A$ and we conclude.

\TeXButton{End Proof}{\endproof}

The following Lemma, which is an easy consequence of a result of Merkurjev,
will be the main tool in reducing the proof of the main theorem from $G={\rm %
GL}_{n,k}$ to its maximal torus $T$.

\begin{lemma}
\label{kurjev}Let $X$ be a noetherian separated algebraic space over $k$
with an action of a split reductive group $G$ over $k$ such that $\pi
_1\left( G\right) $ ({\bf \cite{Me}}, 1.1) is torsion free. Then, if $T$
denotes a maximal torus in $G$, the canonical morphism 
\[
K_{*}^{\prime }\left( X,G\right) \otimes _{{\rm R}\left( G\right) }{\rm R}%
\left( T\right) \longrightarrow K_{*}^{\prime }\left( X,T\right) 
\]
is an isomorphism.
\end{lemma}

\TeXButton{Proof}{\proof} Let $B\supseteq T$ be a Borel subgroup of $G$.
Since ${\rm R}\left( B\right) \simeq {\rm R}(T)$ and $K_{*}^{\prime }\left(
X,B\right) \simeq K_{*}^{\prime }\left( X,T\right) $ ({\bf \cite{Th5}},
proof of Th. 1.13, p. 594), by {\bf \cite{Me},} Prop. 4.1, the canonical
ring morphism 
\[
K_{*}^{\prime }\left( X,G\right) \otimes _{{\rm R}\left( G\right) }{\rm R}%
\left( T\right) \longrightarrow K_{*}^{\prime }\left( X,T\right) 
\]
is an isomorphism.

Since Merkurjev states his theorem for a scheme, we briefly indicate how it
extends to a noetherian separated algebraic space $X$ over $k$. By {\bf \cite
{Th2}} Lemma 4.3, there exists an open dense $G$-invariant separated
subscheme $U\subset X$. Since Merkurjev's map commutes with localization, by
the localization sequence and noetherian induction it is enough to know the
result for $U$. And this is given in {\bf \cite{Me},} Prop. 4.1. Note that
by {\bf \cite{Me}, }1.22, ${\rm R}\left( T\right) $ is flat (actually free)
over ${\rm R}\left( G\right) $ and therefore the localization sequence
remain exact after tensoring with ${\rm R}\left( T\right) $. 
\TeXButton{End Proof}{\endproof}

The following is Lemma 3.2 in {\bf \cite{Vi2}}. It will be used frequently
in the rest of the paper and it is stated here for the convenience of the
reader:

\begin{lemma}
\label{smalllemma}Let $W$ be a finite group acting on the left on a set $%
{\cal A}$ and let ${\cal B}\subseteq {\cal A}$ be a set of representatives
for the orbits. Assume that $W$ acts on the left on a product of abelian
groups of the type $\prod_{\alpha \in {\cal A}}M_\alpha $ in such a way that 
\[
sM_\alpha =M_{s\alpha }
\]
for any $s\in W$.

For each $\alpha \in {\cal B}$, let us denote by $W_\alpha $ the stabilizer
of $\alpha $ in $W$. Then the canonical projection 
\[
\prod_{\alpha \in {\cal A}}M_\alpha \longrightarrow \prod_{\alpha \in {\cal B%
}}M_\alpha 
\]
induces an isomorphism 
\[
\left( \prod_{\alpha \in {\cal A}}M_\alpha \right) ^W\longrightarrow
\prod_{\alpha \in {\cal B}}\left( M_\alpha \right) ^{W_\alpha }.
\]
\end{lemma}

\section{The main theorem: the split torus case}

In this section $T$ will be a split torus over $k$.

\begin{proposition}
\label{above}Let $T^{\prime }\subset T$ a closed subgroup scheme
(diagonalizable, by {\bf \cite{SGA3},} IX 8.1), finite over $k$. Then the
canonical morphism 
\[
\delta :{\rm R}\left( T^{\prime }\right) _{\Lambda _{T^{\prime
}}}\longrightarrow \prod_{_{_{\Sb \sigma \text{{\rm \ }dual cyclic} \\ %
\sigma \subseteq T^{\prime } \endSb }}}\widetilde{{\rm R}}\left( \sigma
\right) _{\Lambda _{T^{\prime }}}
\]
is a ring isomorphism.
\end{proposition}

\TeXButton{Proof}{\proof} Since both ${\rm R}\left( T^{\prime }\right)
_{\Lambda _{T^{\prime }}}$ and $\prod \widetilde{{\rm R}}\left( \sigma
\right) _{\Lambda _{T^{\prime }}}$ are free $\Lambda _{T^{\prime }}$-modules
of finite rank, it is enough to prove that for any nonzero prime $p\nmid
\left| T^{\prime }\right| $ the induced morphism of ${\Bbb F}_p$-vector
spaces 
\begin{equation}
{\rm R}\left( T^{\prime }\right) _{\Lambda _{T^{\prime }}}\otimes _{{\Bbb Z}}%
{\Bbb F}_p\longrightarrow \prod_{_{\Sb \sigma \text{{\rm \ }dual cyclic} \\ %
\sigma \subseteq T^{\prime } \endSb }}\widetilde{{\rm R}}\left( \sigma
\right) _{\Lambda _{T^{\prime }}}\otimes _{{\Bbb Z}}{\Bbb F}_p
\label{effepi}
\end{equation}
is an isomorphism. Now, for any finite abelian group $A$ we have an equality 
$\left| A\right| =\sum_{A\twoheadrightarrow C}\varphi (\left| C\right| )$
where $\varphi $ denotes the Euler function, $\left| H\right| $ denotes the
order of the group $H$ and the sum is extended to all cyclic quotients of $A$%
; applying this to the group of characters $\widehat{T^{\prime }}$ (so that
the corresponding cyclic quotients $C$ are exactly the group of characters $%
\widehat{\sigma }$ for $\sigma $ dual cyclic subgroups of $T^{\prime }$) we
see that the ranks of both sides in (\ref{effepi}) coincide with $\left|
T^{\prime }\right| $ and it is then enough to prove that (\ref{effepi}) is
injective. Define a morphism 
\[
f:\prod_{_{\Sb \tau \text{{\rm \ }dual cyclic} \\ \tau \subseteq T^{\prime }
\endSb }}\widetilde{{\rm R}}\left( \tau \right) _{\Lambda _{T^{\prime
}}}\longrightarrow \prod_{\Sb \sigma \text{{\rm \ }dual cyclic} \\ \sigma
\subseteq T^{\prime } \endSb }{\rm R}\left( \sigma \right) _{\Lambda
_{T^{\prime }}}
\]
of ${\rm R}\left( T^{\prime }\right) _{\Lambda _{T^{\prime }}}$-modules by
requiring for any dual cyclic subgroup $\sigma \subseteq T^{\prime }$, the
commutativity of the following diagram 
\[
\begin{tabular}{lll}
$\prod_{\Sb \tau \text{{\rm \ }dual cyclic} \\ \tau \subseteq T^{\prime }
\endSb }\widetilde{{\rm R}}\left( \tau \right) _{\Lambda _{T^{\prime }}}$ & $%
\stackrel{f}{\longrightarrow }$ & $\prod_{_{_{\Sb \sigma \text{{\rm \ }dual
cyclic} \\ \sigma \subseteq T^{\prime } \endSb }}}{\rm R}\left( \sigma
\right) _{\Lambda _{T^{\prime }}}$ \\ 
\qquad $^{_{{\rm Pr}_\sigma }}\downarrow $ &  & $\qquad \qquad \downarrow ^{%
{\rm pr}_\sigma }$ \\ 
$\prod_{_{_{\tau \subseteq \sigma }}}\widetilde{{\rm R}}\left( \tau \right)
_{\Lambda _{T^{\prime }}}$ & $\stackunder{\varphi }{\widetilde{%
\longleftarrow }}$ & $\quad \qquad {\rm R}\left( \sigma \right) _{\Lambda
_{T^{\prime }}}$%
\end{tabular}
\]
where ${\rm Pr}_\sigma $ and ${\rm pr}_\sigma $ are the obvious projections
and $\varphi $ is the isomorphism 
\[
{\rm R}\left( \sigma \right) _{\Lambda _{T^{\prime }}}\stackrel{\tprod_{\tau
\subseteq \sigma }{\rm res}_\tau ^\sigma }{\longrightarrow }\prod_{\tau
\subseteq \sigma }{\rm R}\left( \tau \right) _{\Lambda _{T^{\prime }}}%
\stackrel{\left( \widetilde{{\rm pr}}_\tau \right) _\tau }{\longrightarrow }%
\prod_{\tau \subseteq \sigma }\widetilde{{\rm R}}\left( \tau \right)
_{\Lambda _{T^{\prime }}}
\]
induced by (\ref{decomposition}). Obviously, $f\circ \delta $ coincides with
the map 
\[
\prod_{_{\Sb \sigma \text{{\rm \ }dual cyclic} \\ \sigma \subseteq T^{\prime
} \endSb }}{\rm res}_\sigma ^{T^{\prime }}:{\rm R}\left( T^{\prime }\right)
_{\Lambda _{T^{\prime }}}\longrightarrow \prod_{\Sb \sigma \text{{\rm \ }%
dual cyclic} \\ \sigma \subseteq T^{\prime } \endSb }{\rm R}\left( \sigma
\right) _{\Lambda _{T^{\prime }}}
\]
so we are reduced to proving that 
\[
{\rm R}\left( T^{\prime }\right) _{\Lambda _{T^{\prime }}}\otimes _{{\Bbb Z}}%
{\Bbb F}_p\longrightarrow \prod_{_{_{\Sb \sigma \text{{\rm \ }dual cyclic}
\\ \sigma \subseteq T^{\prime } \endSb }}}{\rm R}\left( \sigma \right)
_{\Lambda _{T^{\prime }}}\otimes _{{\Bbb Z}}{\Bbb F}_p
\]
is injective i.e. that if $A$ is a finite abelian group and $p\nmid \left|
A\right| $ then 
\begin{equation}
\varphi :{\Bbb F}_p\left[ A\right] \longrightarrow \prod_{C\in \left\{ \text{%
cyclic quotients of }A\right\} }{\Bbb F}_p\left[ C\right]   \label{artin}
\end{equation}
is injective. If $\widehat{A}=Hom_{{\tt AbGrps}}\left( A,{\Bbb C}^{*}\right) 
$ denotes the complex characters group of $A$, ${\rm R}\left( \widehat{A}%
\right) \simeq {\Bbb Z}\left[ A\right] $ and 
\[
\varphi =\prod_{\widehat{C}}{\rm res}_{\widehat{C}}^{\widehat{A}}:{\rm R}(%
\widehat{A})\longrightarrow \prod_{\widehat{C}\in \left\{ \text{cyclic
subgroups of }\widehat{A}\right\} }{\rm R}(\widehat{C})
\]
Since $p\nmid \left| A\right| $ it is enough to prove that if $\xi \in {\rm R%
}(\widehat{A})\dot \otimes _{{\Bbb Z}}{\Bbb Z}\left[ 1/\left| A\right|
\right] $ has image via 
\[
{\rm res}_{\widehat{C}}^{\widehat{A}}\otimes {\Bbb Z}\left[ 1/\left|
A\right| \right] :{\rm R}(\widehat{A})\otimes _{{\Bbb Z}}{\Bbb Z}\left[
1/\left| A\right| \right] \longrightarrow {\rm R}(\widehat{C})\otimes _{%
{\Bbb Z}}{\Bbb Z}\left[ 1/\left| A\right| \right] 
\]
contained in $p\left( {\rm R}(\widehat{C})\otimes _{{\Bbb Z}}{\Bbb Z}\left[
1/\left| A\right| \right] \right) $ for each cyclic $\widehat{C}\subseteq 
\widehat{A}$, then $\xi \in p\left( {\rm R}(\widehat{A})\dot \otimes _{{\Bbb %
Z}}{\Bbb Z}\left[ 1/\left| A\right| \right] \right) $.

By {\bf \cite{Se}} p. 73, there exists $(\theta _{\widehat{C}}^{\prime })_{%
\widehat{C}}\in \prod_{\widehat{C}\in \left\{ \text{cyclic subgroups of }%
\widehat{A}\right\} }{\rm R}(\widehat{C})\otimes _{{\Bbb Z}}{\Bbb Z}\left[
1/\left| A\right| \right] $ such that 
\[
1=\sum_{\widehat{C}}\left( {\rm ind}_{\widehat{A}}^{\widehat{C}}\otimes 
{\Bbb Z}\left[ 1/\left| A\right| \right] \right) (\theta _{\widehat{C}%
}^{\prime }); 
\]
therefore 
\[
\xi =\sum_{\widehat{C}}\xi \left( {\rm ind}_{\widehat{A}}^{\widehat{C}%
}\otimes {\Bbb Z}\left[ 1/\left| A\right| \right] \right) (\theta _{\widehat{%
C}}^{\prime })=\sum_{\widehat{C}}\left( {\rm ind}_{\widehat{A}}^{\widehat{C}%
}\otimes {\Bbb Z}\left[ 1/\left| A\right| \right] \right) \left( \theta _{%
\widehat{C}}^{\prime }\left( {\rm res}_{\widehat{C}}^{\widehat{A}}\otimes 
{\Bbb Z}\left[ 1/\left| A\right| \right] \right) \left( \xi \right) \right) 
\]
(by the projection formula) and we conclude. \TeXButton{End Proof}{\endproof}

\begin{remark}
The proof of Prop. \ref{above} is similar to the proof of Prop. (1.5) of 
{\bf \cite{Vi2}} which is however incomplete; that is why we have decided to
give all details here.
\end{remark}

\begin{corollary}
\label{bunch}(i) If $\sigma \neq \sigma ^{\prime }$ are dual cyclic
subgroups of $T$, we have $\widetilde{{\rm R}}\left( \sigma \right) _{\sigma
^{\prime }}=0$ and $\widetilde{{\rm R}}\left( \sigma \right) _\sigma =%
\widetilde{{\rm R}}\left( \sigma \right) $;

(ii) If $T^{\prime }\subset T$ is a closed subgroup scheme, finite over $k$
and $\sigma $ is a dual cyclic subgroup of $T$, we have ${\rm R}\left(
T^{\prime }\right) _\sigma =0$ if $\sigma \nsubseteq T^{\prime }$;

(iii) If $T^{\prime }\subset T$ is a closed subgroup scheme, finite over $k$%
, the canonical morphism of ${\rm R}(T)$-algebras 
\[
{\rm R}\left( T^{\prime }\right) _{\Lambda _{T^{\prime }}}\longrightarrow
\prod_{_{\Sb \sigma \text{{\rm dual cyclic}} \\ \sigma \subseteq T^{\prime }
\endSb }}{\rm R}\left( T^{\prime }\right) _\sigma 
\]
is an isomorphism.
\end{corollary}

\TeXButton{Proof}{\proof}(i) Suppose $\sigma \neq \sigma ^{\prime }$ and let 
$T^{\prime }\subset T$ be the closed subgroup scheme of $T$ generated by $%
\sigma $ and $\sigma ^{\prime }$. The obvious morphism $\pi :{\rm R}%
(T)_{\Lambda _{T^{\prime }}}\rightarrow \widetilde{{\rm R}}\left( \sigma
\right) _{\Lambda _{T^{\prime }}}\times \widetilde{{\rm R}}\left( \sigma
^{\prime }\right) _{\Lambda _{T^{\prime }}}$ factors through ${\rm R}%
(T^{\prime })_{\Lambda _{T^{\prime }}}\rightarrow \widetilde{{\rm R}}\left(
\sigma \right) _{\Lambda _{T^{\prime }}}\times \widetilde{{\rm R}}\left(
\sigma ^{\prime }\right) _{\Lambda _{T^{\prime }}}$ which is an epimorphism
by Proposition \ref{above}. If $\xi \in {\rm R}(T)_{\Lambda _{T^{\prime }}}$
with $\pi \left( \xi \right) =\left( 0,1\right) \otimes 1$, we have 
\[
\xi \in S_{\sigma ^{\prime }}\cap \ker ({\rm R}(T)_{\Lambda _{T^{\prime
}}}\rightarrow \widetilde{{\rm R}}\left( \sigma \right) _{\Lambda
_{T^{\prime }}}). 
\]
Then $\widetilde{{\rm R}}\left( \sigma \right) _{\sigma ^{\prime }}=0$. The
second assertion is obvious.

(ii) and (iii) follow immediately from (i) and Proposition \ref{above}. 
\TeXButton{End Proof}{\endproof}

\smallskip\ 

Now let $X$ be a regular noetherian separated $k$-algebraic space on which $%
T $ acts with finite stabilizers and let $\Lambda \doteq \Lambda _{\left(
T,X\right) }$. Obviously, ${\cal C}\left( T\right) $ is just the set of
essential dual cyclic subgroups of $T$, since $T$ is abelian.

\begin{proposition}
\label{concentration}(i) If $j_\sigma :X^\sigma \hookrightarrow X$ denotes
the inclusion, the pushforward $\left( j_\sigma \right) _{*}$ induces an
isomorphism 
\[
K_{*}^{\prime }\left( X^\sigma ,T\right) _\sigma \longrightarrow
K_{*}^{\prime }\left( X,T\right) _\sigma 
\]

(ii) The canonical localization morphism 
\[
K_{*}^{\prime }\left( X,T\right) _\Lambda \longrightarrow \prod_{\sigma \in 
{\cal \QTR{mathcal}{C}}\left( T\right) }K_{*}^{\prime }\left( X,T\right)
_\sigma 
\]
is an isomorphism and the product on the left is finite.
\end{proposition}

\TeXButton{Proof}{\proof} (i) The proof is the same as that of {\bf \cite
{Th1}} Th. 2.1 where the use of {\bf \cite{Th1}} Cor. 1.3 has to be
substituted by that of our Cor. \ref{bunch} (ii) above since we use a
localization different from Thomason's.

(ii) By the generic slice theorem for torus actions ({\bf \cite{Th2}}, Prop.
4.10), there exists a $T$-invariant nonempty open subspace $U\subset X$, a
closed (necessarily diagonalizable) subgroup $T^{\prime }$ of $T$ and a $T$%
-equivariant isomorphism 
\[
U\simeq T/T^{\prime }\times \left( U/T\right) \simeq (U/T)\times ^{T^{\prime
}}T 
\]
Since $U$ is nonempty and $T$ acts on $X$ with finite stabilizers, $%
T^{\prime }$ is finite over $k$ and $K_{*}^{\prime }\left( U,T\right) \simeq
K_{*}^{\prime }\left( U/T\right) \otimes _{{\Bbb Z}}{\rm R}\left( T^{\prime
}\right) $, by Morita equivalence theorem ({\bf \cite{Th3}}, Proposition
6.2) and \cite{Th2} Lemma 5.6. By Cor. \ref{bunch} (ii), the Proposition for 
$X=U$ follows from Cor. \ref{bunch} (iii).

By noetherian induction and the localization sequence for $K^{\prime }$%
-groups ({\bf \cite{Th3}}, Theorem 2.7), the statement for $U$ implies that
for $X$.

Again using noetherian induction, Thomason's generic slice theorem for torus
actions and (i), one similarly shows that the product $\prod_{\sigma \in 
{\cal C}\left( T\right) }K_{*}^{\prime }\left( X,T\right) _\sigma $ is
finite. \TeXButton{End Proof}{\endproof}

By Prop. \ref{concentration}, there is an induced isomorphism (of ${\rm R}%
\left( T\right) $-modules, not a ring isomorphism due to the composition
with pushforwards) 
\begin{equation}
\prod_{\sigma \in {\cal C}\left( T\right) }K_{*}^{\prime }\left( X^\sigma
,T\right) _\sigma \longrightarrow K_{*}^{\prime }\left( X,T\right) _\Lambda
\label{step}
\end{equation}
As shown in Lemma \ref{construction}, the product morphism $\sigma \times
T\rightarrow T$ induces a morphism 
\begin{equation}
\theta _{T,\sigma }:K_{*}^{\prime }\left( X^\sigma ,T\right) _\sigma
\longrightarrow K_{*}^{^{\prime }}\left( X^\sigma ,T\right) _{{\rm geom}%
}\otimes \widetilde{{\rm R}}\left( \sigma \right) _\Lambda .  \label{brick}
\end{equation}

\begin{proposition}
\label{end'}For any $\sigma \in {\cal C}\left( T\right) $, $\theta
_{T,\sigma }$ is an isomorphism.
\end{proposition}

\TeXButton{Proof}{\proof} We will write $\theta _{X,\sigma }$ for $\theta
_{T,\sigma }$ in order to emphasize the dependence of the map on the space.
We proceed by noetherian induction on $X^\sigma $. Let $X^{\prime }\subseteq
X^\sigma $ be a $T$-invariant closed subspace and let us suppose that (\ref
{brick}) is an isomorphism with $X$ replaced by any $T$-invariant proper
closed subspace $Z$ of $X^{\prime }$. By Thomason's generic slice theorem
for torus actions ({\bf \cite{Th2}}, Prop. 4.10), there exists a $T$%
-invariant nonempty open subscheme $U\subset X^{\prime }$, a (necessarily
diagonalizable) subgroup $T^{\prime }$ of $T$ and a $T$-equivariant
isomorphism 
\[
U^\sigma \equiv U\simeq T/T^{\prime }\times \left( U/T\right) \simeq
(U/T)\times ^{T^{\prime }}T. 
\]
Since $U$ is nonempty and $T$ acts on $X$ with finite stabilizers, $%
T^{\prime }$ is finite over $k$ and, obviously, $\Lambda _{T^{\prime
}}\subseteq \Lambda $. Let $Z^\sigma \equiv Z\doteq X^{\prime }\diagdown U$.
Since 
\[
\begin{tabular}{ccccccc}
$\rightarrow $ & $K_{*}^{^{\prime }}\left( Z^\sigma ,T\right) _\sigma $ & $%
\longrightarrow $ & $K_{*}^{\prime }\left( X^{\prime \sigma },T\right)
_\sigma $ & $\longrightarrow $ & $K_{*}^{\prime }\left( U^\sigma ,T\right)
_\sigma $ & $\rightarrow $ \\ 
& $^{\theta _{Z,\sigma }}\downarrow $ &  & $^{\theta _{Y^{\prime },\sigma
}}\downarrow $ &  & $^{\theta _{U,\sigma }}\downarrow $ &  \\ 
$\rightarrow $ & $K_{*}^{^{\prime }}\left( Z^\sigma ,T\right) _{{\rm geom}%
}\otimes \widetilde{{\rm R}}\left( \sigma \right) _\Lambda $ & $%
\longrightarrow $ & $K_{*}^{^{\prime }}\left( X^{\prime \sigma },T\right) _{%
{\rm geom}}\otimes \widetilde{{\rm R}}\left( \sigma \right) _\Lambda $ & $%
\longrightarrow $ & $K_{*}^{^{\prime }}\left( U^\sigma ,T\right) _{{\rm geom}%
}\otimes \widetilde{{\rm R}}\left( \sigma \right) _\Lambda $ & $\rightarrow $%
\end{tabular}
\]
is commutative, by induction hypothesis and the five-lemma it will be enough
to show that $\theta _{U,\sigma }$ is an isomorphism. By Morita equivalence
theorem ({\bf \cite{Th3}}, Proposition 6.2) and \cite{Th2} Lemma 5.6, $%
K_{*}^{\prime }\left( U,T\right) \simeq K_{*}^{\prime }\left( U/T\right)
\otimes _{{\Bbb Z}}{\rm R}\left( T^{\prime }\right) $, so it is enough to
prove that 
\[
\theta _{{\rm Spec}k,\sigma }:K_0^{\prime }\left( {\rm Spec}k,T^{\prime
}\right) _\sigma ={\rm R}\left( T^{\prime }\right) _\sigma \rightarrow
K_0^{^{\prime }}\left( {\rm Spec}k,T^{\prime }\right) _{{\rm geom}}\otimes 
\widetilde{{\rm R}}\left( \sigma \right) _\Lambda ={\rm R}\left( T^{\prime
}\right) _{{\rm geom}}\otimes \widetilde{{\rm R}}\left( \sigma \right)
_\Lambda 
\]
is an isomorphism. But this follows immediately from Cor. \ref{bunch} (i)
and (iii). \TeXButton{End Proof}{\endproof}

Combining Proposition \ref{end'} with (\ref{step}) we get an isomorphism 
\begin{equation}
\Phi _{X,T}:\prod_{\sigma \in {\cal C}\left( T\right) }K_{*}^{\prime }\left(
X^\sigma ,T\right) _{{\rm geom}}\otimes \widetilde{{\rm R}}\left( \sigma
\right) _\Lambda \longrightarrow K_{*}^{^{\prime }}\left( X,T\right)
_\Lambda .  \label{fidix}
\end{equation}
The following lemma is a variant of {\bf \cite{Th1}}, Lemme 3.2, that
already proves the statement below after tensoring with ${\Bbb Q}$.

\begin{lemma}
\label{unit}Let $X$ be a noetherian regular separated algebraic space over $k
$ on which a split $k$-torus acts with finite stabilizers and let $\sigma
\in {\cal C}\left( T\right) $. Let $X^\sigma $ denote the regular $\sigma $%
-fixed subscheme, $j_\sigma :X^\sigma \hookrightarrow X$ the regular closed
immersion ({\bf \cite{Th1}}, Proposition 3.1) and ${\cal N}\left( j_\sigma
\right) $ the corresponding locally free conormal sheaf. Then, for any $T$%
-invariant algebraic subspace $Y$ of $X^\sigma $, the cap-product 
\[
\lambda _{-1}({\cal N}\left( j_\sigma \right) )\cap \left( -\right)
:K_{*}^{\prime }\left( Y,T\right) _\sigma \longrightarrow K_{*}^{\prime
}\left( Y,T\right) _\sigma 
\]
is an isomorphism.
\end{lemma}

\TeXButton{Proof}{\proof} We proceed by noetherian induction on closed $T$%
-invariant subspaces $Y$ of $X^\sigma $. The statement is trivial for $%
Y=\emptyset $, so let us suppose $Y$ nonempty and 
\[
\lambda _{-1}({\cal N}\left( j_\sigma \right) )\cap \left( -\right)
:K_{*}^{\prime }\left( Z,T\right) _\sigma \longrightarrow K_{*}^{\prime
}\left( Z,T\right) _\sigma 
\]
an isomorphism for any proper $T$-invariant closed subspace $Z$ of $Y$. By
Thomason's generic slice theorem for torus actions ({\bf \cite{Th2}}, Prop.
4.10), there exists a $T$-invariant nonempty open subscheme $U\subset Y$, a
closed (necessarily diagonalizable) subgroup $T^{\prime }$ of $T$ and a $T$%
-equivariant isomorphism 
\[
U^\sigma \equiv U\simeq T/T^{\prime }\times \left( U/T\right) .
\]
Since $U$ is nonempty and $T$ acts on $X$ with finite stabilizers, $%
T^{\prime }$ is finite over $k$. Using the localization sequence and the
five-lemma, we reduce ourselves to showing that 
\[
\lambda _{-1}({\cal N}\left( j_\sigma \right) )\cap \left( -\right)
:K_{*}^{\prime }\left( U,T\right) _\sigma \longrightarrow K_{*}^{\prime
}\left( U,T\right) _\sigma 
\]
is an isomorphism. For this, it is enough to show that (the restriction of) $%
\lambda _{-1}\left( {\cal N}\left( j_\sigma \right) \right) $ is a unit in $%
K_0\left( U,T\right) _\sigma \simeq K_0\left( U/T\right) _\Lambda \otimes 
{\rm R}\left( T^{\prime }\right) _\sigma $ ({\bf \cite{Th3}}, Proposition
6.2). Decomposing ${\cal N}\left( j_\sigma \right) $ according to the
characters of $T^{\prime }$ we may write, shrinking $U$ if necessary, 
\[
{\cal N}\left( j_\sigma \right) =\bigoplus_{\rho \in \widehat{T^{\prime }}}%
{\cal O}_{U/T}^{r_\rho }\otimes {\cal L}_\rho 
\]
where ${\cal L}_\rho $ is the line bundle attached to the $T^{\prime }$%
-character $\rho $ and $r_\rho \geq 0$ and therefore $\lambda _{-1}({\cal N}%
\left( j_\sigma \right) )=\prod_{\rho \in \widehat{T^{\prime }}}\left(
1-\rho \right) ^{r_\rho }$ in $K_0\left( U/T\right) \otimes {\rm R}\left(
T^{\prime }\right) $. The localization map ${\rm R}\left( T^{\prime }\right)
_\Lambda \rightarrow {\rm R}\left( T^{\prime }\right) _\sigma \simeq 
\widetilde{{\rm R}}\left( \sigma \right) _\Lambda $ coincides with the
composition 
\[
{\rm R}\left( T^{\prime }\right) _\Lambda \stackrel{\pi _\sigma }{%
\longrightarrow }{\rm R}\left( \sigma \right) _\Lambda \stackrel{p_\sigma }{%
\longrightarrow }\widetilde{{\rm R}}\left( \sigma \right) _\Lambda 
\]
of the restriction to $\sigma $ followed by the projection (Cor. \ref{bunch}%
) and then 
\[
\left( {\rm id}_{K_0\left( U/T\right) _\Lambda }\otimes \pi _\sigma \right)
\left( {\cal N}\left( j_\sigma \right) \right) =\bigoplus_{\chi \in \widehat{%
\sigma }\backslash \left\{ 0\right\} }{\cal O}_{U/T}^{r_\chi }\otimes {\cal L%
}_\chi 
\]
in $K_0\left( U/T\right) _\Lambda \otimes {\rm R}\left( \sigma \right)
_\Lambda $, where the summand omits the trivial character since the
decomposition of ${\cal N}\left( j_\sigma \right) $ according to the
characters of $\sigma $ has vanishing fixed subsheaf ${\cal N}\left(
j_\sigma \right) _0$ (see, e.g., {\bf \cite{Th1}}, Prop. 3.1). Therefore 
\[
\lambda _{-1}\left( \left( {\rm id}_{K_0\left( U/T\right) _\Lambda }\otimes
\pi _\sigma \right) \left( {\cal N}\left( j_\sigma \right) \right) \right)
=\prod_{\chi \in \widehat{\sigma }\backslash \left\{ 0\right\} }\left(
1-\chi \right) ^{r_\chi }
\]
and it is enough to show that the image of $1-\chi $ in $\widetilde{{\rm R}}%
\left( \sigma \right) _\Lambda $ via $p_\sigma $ is a unit for any
nontrivial character $\chi $ of $\sigma $. Now, the image of such a $\chi $
in 
\[
\widetilde{{\rm R}}\left( \sigma \right) _\Lambda \simeq \frac{\Lambda
\left[ t\right] }{(\Phi _{\left| \sigma \right| })}
\]
($\Phi _{\left| \sigma \right| }$ being the $\left| \sigma \right| $-th
cyclotomic polynomial) is of the form $\left[ t^l\right] $ for some $1\leq
l<\left| \sigma \right| $, where $\left[ -\right] $ denotes the class $%
\limfunc{mod}\Phi _{\left| \sigma \right| }$; therefore the cokernel of the
multiplication by $1-\left[ t^l\right] $ in $\Lambda \left[ t\right] \diagup
\left( \Phi _{\left| \sigma \right| }\right) $ is 
\[
\frac{\Lambda \left[ t\right] }{(\Phi _{\left| \sigma \right| },1-t^l)}=0
\]
since $\Phi _{\left| \sigma \right| }$ and $\left( 1-t^l\right) $ are
relatively prime in $\Lambda \left[ t\right] $, for $1\leq l<\left| \sigma
\right| $. Thus $1-\left[ t^l\right] $ is a unit in $\Lambda \left[ t\right]
\diagup (\Phi _{\left| \sigma \right| })$ and we conclude. 
\TeXButton{End Proof}{\endproof}

\smallskip\ 

We are now able to prove our main theorem for $G=T$:

\begin{theorem}
\label{torus}If $X$ is a regular noetherian separated $k$-algebraic space ,
then 
\[
\Psi _{X,T}:K_{*}\left( X,T\right) _\Lambda \longrightarrow \prod_{\sigma
\in {\cal C}\left( T\right) }K_{*}\left( X^\sigma ,T\right) _{{\rm geom}%
}\otimes \widetilde{{\rm R}}\left( \sigma \right) _\Lambda 
\]
is an isomorphism of ${\rm R}\left( T\right) $-algebras.
\end{theorem}

\TeXButton{Proof}{\proof}Recall (Appendix) that $K_{*}\left( X,T\right)
\simeq K_{*}^{\prime }\left( X,T\right) $ and $K_{*}\left( X^\sigma
,T\right) \simeq K_{*}^{\prime }\left( X^\sigma ,T\right) $, since both $X$
and $X^\sigma $ are regular ({\bf \cite{Th1}}, Proposition 3.1). Since $\Phi
_{X,T}$ is an isomorphism of ${\rm R}\left( T\right) $-modules, it is enough
to show that the composition $\Psi _{X,T}\circ \Phi _{X,T}$ is an
isomorphism. A careful inspection of the definitions of $\Psi _{X,T}$ $\ $%
and $\Phi _{X,T}$, easily reduce the problem to proving that, for any $%
\sigma $ $\in {\cal \QTR{mathcal}{C}}\left( T\right) $, the composition 
\[
K_{*}^{\prime }\left( X^\sigma ,T\right) _\sigma \stackrel{j_{\sigma *}}{%
\longrightarrow }K_{*}^{\prime }\left( X,T\right) _\sigma \stackrel{j_\sigma
^{*}}{\longrightarrow }K_{*}^{\prime }\left( X^\sigma ,T\right) _\sigma 
\]
is an isomorphism, $j_\sigma :X^\sigma \hookrightarrow X$ being the natural
inclusion. Since $j_\sigma $ is regular, there is a self-intersection
formula 
\begin{equation}
j_\sigma ^{*}\circ j_{\sigma *}\left( -\right) =\lambda _{-1}({\cal N}\left(
j_\sigma \right) )\cap \left( -\right)   \label{self}
\end{equation}
${\cal N}\left( j_\sigma \right) $ being the conormal sheaf associated to $%
j_\sigma $ and we conclude by lemma \ref{unit}. To prove the
self-intersection formula (\ref{self}), we adapt {\bf \cite{Th1}}, proof of
Lemme 3.3. First of all, by Proposition \ref{concentration} (i), $j_{\sigma
*}$ is an isomorphism so it is enough to prove that $j_{\sigma *}j_\sigma
^{*}j_{\sigma *}\left( -\right) =j_{\sigma *}(\lambda _{-1}\left( {\cal N}%
\left( j_\sigma \right) \right) \cap \left( -\right) )$. By projection
formula (Appendix, Proposition \ref{a4}), we have 
\[
j_{\sigma *}j_\sigma ^{*}j_{\sigma *}\left( -\right) =j_{\sigma *}j_\sigma
^{*}\left( 1\right) \cap j_{\sigma *}\left( -\right) =j_{\sigma *}j_\sigma
^{*}({\cal O}_X)\cap j_{\sigma *}\left( -\right) =
\]
\[
=j_{\sigma *}({\cal O}_{X^\sigma })\cap j_{\sigma *}\left( -\right)
=j_{\sigma *}(j_\sigma ^{*}({\cal O}_{X^\sigma })\cap \left( -\right) )
\]
Now, as explained in the Appendix, to compute $j_\sigma ^{*}({\cal O}%
_{X^\sigma })$ we choose a complex $F^{*}$ of flat quasi-coherent $G$%
-equivariant Modules on $X$ which is quasi-isomorphic to ${\cal O}_{X^\sigma
}$ and then 
\[
j_\sigma ^{*}({\cal O}_{X^\sigma })=[j_\sigma ^{*}(F^{*})]=[F^{*}\otimes 
{\cal O}_{X^\sigma }]=\sum_i\left( -1\right) ^i[H^i(F^{*}\otimes {\cal O}%
_{X^\sigma })].
\]
But $F^{*}$ is a flat resolution of ${\cal O}_{X^\sigma }$, so $%
H^i(F^{*}\otimes {\cal O}_{X^\sigma })={\rm Tor}_i^{{\cal O}_X}\left( {\cal O%
}_{X^\sigma },{\cal O}_{X^\sigma }\right) \simeq {\rm \bigwedge^i}{\cal N}%
\left( j_\sigma \right) ,$ where the last isomorphism (\cite{SGA6}, VII,
2.5) is natural hence $T$-equivariant. Therefore $j_\sigma ^{*}({\cal O}%
_{X^\sigma })=\lambda _{-1}({\cal N}\left( j_\sigma \right) )$ and we
conclude. \TeXButton{End Proof}{\endproof}

\section{The main theorem: the case of $G={\rm GL}_{n,k}$}

In this section we will use the result for $\Psi _{X,T}$ to deduce the same
result for $\Psi _{X,{\rm GL}_{nk}}$.

\smallskip\ 

\begin{theorem}
\label{GLn}Let $X$ be a noetherian regular separated algebraic space over a
field $k$ on which $G={\rm GL}_{n,k}$ acts with finite stabilizers. Then the
map defined in (\ref{due}) 
\begin{equation}
\Psi _{X,G}:K_{*}\left( X,G\right) _{\Lambda _{\left( G,X\right)
}}\longrightarrow \prod_{\sigma \in {\cal C}\left( G\right) }\left(
K_{*}\left( X^\sigma ,C\left( \sigma \right) \right) _{{\rm geom}}\otimes
_{\Lambda _{\left( G,X\right) }}\widetilde{{\rm R}}\left( \sigma \right)
_{\Lambda _{\left( G,X\right) }}\right) ^{w_G\left( \sigma \right) }
\label{mainGLn}
\end{equation}
is an isomorphism of ${\rm R}\left( G\right) $-algebras and the product on
the right is finite.
\end{theorem}

\smallskip\ 

Throughout this section, entirely devoted to the proof of Theorem \ref{GLn},
we will simply write $G$ for ${\rm GL}_{n,k}$, $\Lambda $ for $\Lambda
_{\left( G,X\right) }$ and $T$ for the maximal torus of diagonal matrices in 
${\rm GL}_{n,k}$.

First of all, let us observe that we can choose each $\sigma \in {\cal C}%
\left( G\right) $ contained in $T$. Moreover, $\Lambda _{\left( T,X\right)
}=\Lambda _{(G,X)}$.

\smallskip\ 

We need the following three preliminary lemmas (\ref{lemma?}, \ref{decompose}
and \ref{galois}).

If $\sigma ,\tau \subset T$ are dual cyclic subgroups, they are conjugate
under the $G\left( k\right) $-action iff they are conjugated via an element
in the Weyl group $S_n$. For any group scheme $H$ with a dual cyclic
subgroup $\sigma \subseteq H$, we denote by ${\frak m}_\sigma ^H$ the kernel
of ${\rm R}\left( H\right) _\Lambda \rightarrow \widetilde{{\rm R}}\left(
\sigma \right) _\Lambda $ and by $\widehat{{\rm R}\left( H\right) _{\Lambda
,\sigma }}$ the completion of ${\rm R}\left( H\right) _\Lambda $ with
respect to the ideal ${\frak m}_\sigma ^H$.

The following Lemma is essentially a variant of Lemma \ref{kurjev} for $%
\sigma $-localizations.

\begin{lemma}
\label{lemma?}Let $G={\rm GL}_{n,k}$, $T$ the maximal torus of $G$
consisting of diagonal matrices and $X$ an algebraic space on which $G$ acts
with finite stabilizers.

(i) for any essential dual cyclic subgroup $\sigma \subseteq T$, the
morphisms 
\[
\omega _{\sigma ,{\rm geom}}:K_{*}^{\prime }\left( X^\sigma ,{\rm C}_G\left(
\sigma \right) \right) _{{\rm geom}}\otimes _{{\rm R}\left( {\rm C}_G\left(
\sigma \right) \right) _\Lambda }{\rm R}\left( T\right) _\Lambda
\longrightarrow K_{*}^{\prime }\left( X^\sigma ,T\right) _{{\rm geom}}
\]
\[
\omega _\sigma :K_{*}^{\prime }\left( X^\sigma ,{\rm C}_G\left( \sigma
\right) \right) _\sigma \otimes _{{\rm R}\left( {\rm C}_G\left( \sigma
\right) \right) _\Lambda }{\rm R}\left( T\right) _\Lambda \longrightarrow
K_{*}^{\prime }\left( X^\sigma ,T\right) _\sigma 
\]
induced by $T\hookrightarrow {\rm C}_G\left( \sigma \right) $ are
isomorphisms;

(ii) for any essential dual cyclic subgroup $\sigma \subseteq T$, 
\[
({\frak m}_\sigma ^{{\rm C}_G\left( \sigma \right) })^N\cdot K_{*}^{^{\prime
}}\left( X^\sigma ,{\rm C}_G\left( \sigma \right) \right) _\sigma =0\text{,
\quad }N\gg 0
\]
and the morphism induced by $T\hookrightarrow {\rm C}_G\left( \sigma \right) 
$%
\[
\widehat{\omega _\sigma }:K_{*}^{\prime }\left( X^\sigma ,{\rm C}_G\left(
\sigma \right) \right) _\sigma \otimes _{\widehat{{\rm R}\left( {\rm C}%
_G\left( \sigma \right) \right) _{\Lambda ,\sigma }}}\widehat{{\rm R}\left(
T\right) _{\Lambda ,\sigma }}\longrightarrow K_{*}^{\prime }\left( X^\sigma
,T\right) _\sigma 
\]
is an isomorphism.
\end{lemma}

\TeXButton{Proof}{\proof}(i) Since ${\rm C}_G\left( \sigma \right) $ is
isomorphic to a product of general linear groups over $k$ and $T$ is a
maximal torus in ${\rm C}_G\left( \sigma \right) $, by Lemma \ref{kurjev},
the canonical ring morphism 
\begin{equation}
K_{*}^{\prime }\left( X,{\rm C}_G\left( \sigma \right) \right) \otimes _{%
{\rm R}\left( {\rm C}_G\left( \sigma \right) \right) }{\rm R}\left( T\right)
\longrightarrow K_{*}^{\prime }\left( X,T\right)  \label{merk}
\end{equation}
is an isomorphism. If $H\subseteq G$ is a subgroup scheme, we denote by $%
S_\sigma ^H$ the multiplicative subset of ${\rm R}\left( H\right) _\Lambda $
consisting of the elements sent to $1$ by the canonical ring homomorphism $%
{\rm R}\left( H\right) _\Lambda \rightarrow \widetilde{{\rm R}}\left( \sigma
\right) _\Lambda $. By (\ref{merk}), $\omega _\sigma $ coincides with the
composition

\[
K_{*}^{\prime }\left( X^\sigma ,{\rm C}_G\left( \sigma \right) \right)
_\sigma \otimes _{{\rm R}\left( {\rm C}_G\left( \sigma \right) \right)
_\Lambda }{\rm R}\left( T\right) _\Lambda \simeq 
\]
\[
K_{*}^{\prime }\left( X^\sigma ,T\right) \otimes _{{\rm R}\left( {\rm C}%
_G\left( \sigma \right) \right) _\Lambda }\left( \left( S_\sigma ^{{\rm C}%
_G\left( \sigma \right) }\right) ^{-1}{\rm R}\left( {\rm C}_G\left( \sigma
\right) \right) _\Lambda \right) \otimes _{{\rm R}\left( {\rm C}_G\left(
\sigma \right) \right) _\Lambda }{\rm R}\left( T\right) _\Lambda \rightarrow 
\]

\[
\stackrel{{\rm id}\otimes \nu _\sigma }{\longrightarrow }K_{*}^{\prime
}\left( X^\sigma ,{\rm C}_G\left( \sigma \right) \right) \otimes _{{\rm R}%
\left( {\rm C}_G\left( \sigma \right) \right) _\Lambda }\left( S_\sigma
^T\right) ^{-1}{\rm R}\left( T\right) _\Lambda \simeq K_{*}^{\prime }\left(
X^\sigma ,T\right) _\sigma 
\]
where 
\begin{equation}
\nu _\sigma :\left( S_\sigma ^{{\rm C}_G\left( \sigma \right) }\right) ^{-1}%
{\rm R}\left( {\rm C}_G\left( \sigma \right) \right) _\Lambda \otimes _{{\rm %
R}\left( {\rm C}_G\left( \sigma \right) \right) _\Lambda }{\rm R}\left(
T\right) _\Lambda \rightarrow \left( S_\sigma ^T\right) ^{-1}{\rm R}\left(
T\right) _\Lambda  \label{nusigma}
\end{equation}
is induced by $T\hookrightarrow {\rm C}_G\left( \sigma \right) $ and the
last isomorphism follows from (\ref{merk}); the same is true for $\omega
_{\sigma ,{\rm geom}}$. Therefore it is enough to prove that $\nu _\sigma $
and 
\[
\nu _{\sigma ,{\rm geom}}:\left( S_1^{{\rm C}_G\left( \sigma \right)
}\right) ^{-1}{\rm R}\left( {\rm C}_G\left( \sigma \right) \right) _\Lambda
\otimes _{{\rm R}\left( {\rm C}_G\left( \sigma \right) \right) _\Lambda }%
{\rm R}\left( T\right) _\Lambda \rightarrow \left( S_1^T\right) ^{-1}{\rm R}%
\left( T\right) _\Lambda 
\]
are isomorphisms, i.e. if $S_\tau $ denotes the image of $S_1^{{\rm C}%
_G\left( \sigma \right) }$ via the restriction map 
\[
{\rm R}\left( {\rm C}_G\left( \sigma \right) \right) _\Lambda
\longrightarrow {\rm R}\left( T\right) _\Lambda 
\]
that $S_\tau ^T/1$ consists of units in $\left( S_\tau \right) ^{-1}{\rm R}%
\left( T\right) _\Lambda $, for $\tau =1$ and $\tau =\sigma $.

If $\Delta _\sigma $ denotes the Weyl group of ${\rm C}_G\left( \sigma
\right) $, which is a product of symmetric groups, we have ${\rm R}\left( 
{\rm C}_G\left( \sigma \right) \right) \simeq {\rm R}\left( T\right)
^{\Delta _\sigma }$ and therefore 
\[
\left( S_\tau ^{{\rm C}_G\left( \sigma \right) }\right) ^{-1}{\rm R}\left( 
{\rm C}_G\left( \sigma \right) \right) _\Lambda \simeq \left( \left( S_\tau
\right) ^{-1}{\rm R}\left( T\right) _\Lambda \right) ^{\Delta _\sigma }, 
\]
since ${\rm R}\left( T\right) $ is torsion free. Moreover, there is a
commutative diagram 
\[
\begin{tabular}[t]{ccc}
$\left( S_\tau ^{{\rm C}_G\left( \sigma \right) }\right) ^{-1}{\rm R}\left( 
{\rm C}_G\left( \sigma \right) \right) _\Lambda $ & $\hookrightarrow $ & $%
\left( S_\tau \right) ^{-1}{\rm R}\left( T\right) _\Lambda $ \\ 
$^\psi \downarrow $ & $\swarrow _{\doteq \varphi }$ & $\downarrow $ \\ 
$\widetilde{{\rm R}}\left( \tau \right) _\Lambda \simeq \left( S_\tau
^T\right) ^{-1}{\rm R}\left( \tau \right) _\Lambda $ & $\stackunder{\left(
S_\tau ^T\right) ^{-1}{\rm res}_\tau ^T}{\longleftarrow }$ & $\left( S_\tau
^T\right) ^{-1}{\rm R}\left( T\right) _\Lambda $%
\end{tabular}
\]
where $\psi $ is induced by $\widetilde{\pi }_\tau $ and the isomorphism $%
\widetilde{{\rm R}}\left( \tau \right) _\Lambda \simeq \left( S_\tau
^T\right) ^{-1}{\rm R}\left( \tau \right) _\Lambda $ is obtained from
Proposition \ref{above} and Corollary \ref{bunch}. If we define the map 
\[
M:\left( S_\tau \right) ^{-1}{\rm R}\left( T\right) _\Lambda \longrightarrow
\left( \left( S_\tau \right) ^{-1}{\rm R}\left( T\right) _\Lambda \right)
^{\Delta _\sigma } 
\]
\[
\xi \longmapsto \prod_{g\in \Delta _\sigma }g\cdot \xi \text{ }, 
\]
it is easily checked that if for $\xi \in \left( S_\tau \right) ^{-1}{\rm R}%
\left( T\right) _\Lambda $, $\xi $ is a unit if $M\left( \xi \right) $is a
unit and that $\psi \left( M\left( \xi \right) \right) =1$ implies $\xi $ is
a unit in $\left( \left( S_\tau \right) ^{-1}{\rm R}\left( T\right) _\Lambda
\right) ^{\Delta _\sigma }$. But $\varphi $ is $\Delta _\sigma $-equivariant
and therefore $S_\tau ^T/1$ consists of units in $\left( S_\tau \right) ^{-1}%
{\rm R}\left( T\right) _\Lambda $, for $\tau =1$ or $\sigma $.

(ii) Since ${\rm R}\left( {\rm C}_G\left( \sigma \right) \right) \rightarrow 
{\rm R}\left( T\right) $ is faithfully flat, by (i), it is enough to prove
that 
\begin{equation}
\left( {\frak m}_\sigma ^T\right) ^NK_{*}^{\prime }\left( X^\sigma ,T\right)
_\sigma =0\text{ \quad for }N\gg 0.  \label{compl}
\end{equation}
But (\ref{compl}) can be proved using the same technique as in the proof of,
e.g., Proposition \ref{end'}, i.e. noetherian induction together with
Thomason's generic slice theorem for torus actions.

The second part of (ii) follows, arguing as in (i), from the fact that (\ref
{nusigma}) is an isomorphism since 
\[
K_{*}^{\prime }\left( X^\sigma ,{\rm C}_G\left( \sigma \right) \right)
_\sigma \otimes _{{\rm R}\left( {\rm C}_G\left( \sigma \right) \right)
_\Lambda }\widehat{{\rm R}\left( {\rm C}_G\left( \sigma \right) \right)
_{\Lambda ,\sigma }}\simeq K_{*}^{\prime }\left( X^\sigma ,{\rm C}_G\left(
\sigma \right) \right) _\sigma , 
\]
\[
K_{*}^{\prime }\left( X^\sigma ,T\right) _\sigma \otimes _{{\rm R}\left(
T\right) _\Lambda }\widehat{{\rm R}\left( T\right) _{\Lambda ,\sigma }}%
\simeq K_{*}^{\prime }\left( X^\sigma ,T\right) _\sigma . 
\]
\TeXButton{End Proof}{\endproof}

If $\sigma ,\tau \subset T$ are dual cyclic subgroups conjugated under $%
G\left( k\right) $, they are conjugate through an element of the Weyl group $%
S_n$ and we write $\tau \approx _{S_n}\sigma $; moreover, we have ${\frak m}%
_\sigma ^G={\frak m}_\tau ^G$ because conjugation by an element in $S_n$
(actually, by any element in $G\left( k\right) $) induces the identity
morphism on $K$-theory and in particular on the representation ring. Then
there are canonical maps 
\begin{equation}
\widehat{{\rm R}\left( G\right) _{\Lambda ,\sigma }}\otimes _{{\rm R}\left(
G\right) _\Lambda }{\rm R}\left( T\right) _\Lambda \longrightarrow \prod\Sb %
\tau \text{ dual cyclic}  \\ ^{\tau \approx _{S_n}\sigma }  \endSb \widehat{%
{\rm R}\left( T\right) _{\Lambda ,\tau }}  \label{complG}
\end{equation}

\begin{equation}
\widehat{{\rm R}\left( {\rm C}_G\left( \sigma \right) \right) _{\Lambda
,\sigma }}\otimes _{R({\rm C}_G(\sigma ))_\Lambda }{\rm R}\left( T\right)
_\Lambda \longrightarrow \widehat{{\rm R}\left( T\right) _{\Lambda ,\sigma }}
\label{complcisigma}
\end{equation}

\begin{lemma}
\label{decompose}The maps (\ref{complG}) and (\ref{complcisigma}) are
isomorphisms.
\end{lemma}

\TeXButton{Proof}{\proof}Since ${\rm R}\left( G\right) ={\rm R}\left(
T\right) ^{S_n}\rightarrow {\rm R}\left( T\right) $ is finite, the canonical
map $\widehat{{\rm R}\left( G\right) _{\Lambda ,\sigma }}\otimes _{{\rm R}%
\left( G\right) _\Lambda }{\rm R}\left( T\right) _\Lambda \rightarrow 
\widehat{{\rm R}\left( T\right) _\Lambda }^{{\frak m}_\sigma ^G}$ (where $%
\widehat{{\rm R}\left( T\right) _\Lambda }^{{\frak m}_\sigma ^G}$ denotes
the ${\frak m}_\sigma ^G$-adic completion of the ${\rm R}\left( G\right)
_\Lambda $-module ${\rm R}\left( T\right) _\Lambda $) is an isomorphism.
Moreover, ${\rm R}\left( G\right) _\Lambda =\left( {\rm R}\left( T\right)
_\Lambda \right) ^{S_n}$ implies that 
\[
\sqrt{{\frak m}_\sigma ^G{\rm R}\left( T\right) _\Lambda }=\bigcap_{\Sb \tau 
\text{ dual cyclic}  \\ ^{\tau \approx _{S_n}\sigma }  \endSb }\sqrt{{\frak m%
}_\tau ^T}=\bigcap_{\Sb \tau \text{ dual cyclic}  \\ ^{\tau \approx
_{S_n}\sigma }  \endSb }{\frak m}_\tau ^T 
\]
and by Corollary \ref{bunch} (i), ${\frak m}_\tau ^T+{\frak m}_{\tau
^{\prime }}^T={\rm R}\left( T\right) _\Lambda $ if $\tau \neq \tau ^{\prime
} $. By the Chinese remainder lemma, we conclude that (\ref{complG}) is an
isomorphism.

Arguing in the same way, we get that the canonical map 
\[
\widehat{{\rm R}\left( {\rm C}_G\left( \sigma \right) \right) _{\Lambda
,\sigma }}\otimes _{{\rm R}\left( {\rm C}_G\left( \sigma \right) \right)
_\Lambda }{\rm R}\left( T\right) _\Lambda \longrightarrow \prod\Sb \tau 
\text{ dual cyclic}  \\ ^{\tau \approx _{\Delta _\sigma }\sigma }  \endSb 
\widehat{{\rm R}\left( T\right) _{\Lambda ,\tau }} 
\]
is an isomorphism, where $\Delta _\sigma =S_n\cap {\rm C}_G\left( \sigma
\right) $ is the Weyl group of ${\rm C}_G\left( \sigma \right) $ with
respect to $T$ and we write $\tau \approx _{\Delta _\sigma }\sigma $ to
denote that $\tau $ and $\sigma $ are conjugate through an element of $%
\Delta _\sigma $. But $\Delta _\sigma \subset {\rm C}_G\left( \sigma \right) 
$ so that $\tau \approx _{\Delta _\sigma }\sigma $ iff $\tau =\sigma $ and
we conclude that (\ref{complcisigma}) is an isomorphism. 
\TeXButton{End Proof}{\endproof}

\begin{lemma}
\label{galois}For any essential dual cyclic subgroup $\sigma \subseteq G$,
the canonical morphism 
\[
\widehat{{\rm R}\left( G\right) _{\Lambda ,\sigma }}\longrightarrow \widehat{%
{\rm R}\left( {\rm C}_G\left( \sigma \right) \right) _{\Lambda ,\sigma }}
\]
is a finite \'etale Galois cover ({\bf \cite{SGA1}}, Exp. V) with Galois
group $w_G\left( \sigma \right) $.
\end{lemma}

\TeXButton{Proof}{\proof} Since ${\rm R}\left( T\right) $ is flat over ${\rm %
R}\left( G\right) ={\rm R}\left( T\right) ^{S_n}$, we have 
\[
\widehat{{\rm R}\left( G\right) _{\Lambda ,\sigma }}\simeq \widehat{{\rm R}%
\left( G\right) _{\Lambda ,\sigma }}\otimes _{{\rm R}\left( G\right)
_\Lambda }({\rm R}\left( T\right) _\Lambda ^{})^{S_n}\simeq \widehat{({\rm R}%
\left( G\right) _{\Lambda ,\sigma }}\otimes _{{\rm R}\left( G\right)
_\Lambda }{\rm R}\left( T\right) _\Lambda )^{S_n}\simeq (\prod_{\Sb \tau 
\text{ dual cyclic}  \\ ^{\tau \approx _{S_n}\sigma }  \endSb }\widehat{{\rm %
R}\left( T\right) _{\Lambda .\tau }})^{S_n}, 
\]
the last isomorphism being given in Lemma \ref{decompose}. By Lemma \ref
{smalllemma}, we get 
\[
\widehat{{\rm R}\left( G\right) _\sigma }\simeq \left( \widehat{{\rm R}%
\left( T\right) _\sigma }\right) ^{S_{n,\sigma }} 
\]
where $S_n$ acts on the set of dual cyclic subgroups of $T$ which are $S_n$%
-conjugated to $\sigma $ and $S_{n,\sigma }$ denotes the stabilizer of $%
\sigma $. Analogously, denoting by $\Delta _\sigma $ the Weyl group of ${\rm %
C}_G\left( \sigma \right) $, by Lemma \ref{lemma?} (ii) we have 
\begin{eqnarray*}
\widehat{{\rm R}\left( {\rm C}_G\left( \sigma \right) \right) _{\Lambda
,\sigma }} &\simeq &\widehat{{\rm R}\left( {\rm C}_G\left( \sigma \right)
\right) _{\Lambda ,\sigma }}\otimes _{{\rm R}\left( {\rm C}_G\left( \sigma
\right) \right) _\Lambda }\left( {\rm R}\left( T\right) _\Lambda \right)
^{\Delta _\sigma }\simeq \left( \widehat{{\rm R}\left( {\rm C}_G\left(
\sigma \right) \right) _{\Lambda ,\sigma }}\otimes _{{\rm R}\left( {\rm C}%
_G\left( \sigma \right) \right) _\Lambda }{\rm R}\left( T\right) _\Lambda
\right) ^{\Delta _\sigma } \\
\ &\simeq &(\widehat{{\rm R}\left( T\right) _{\Lambda ,\sigma }})^{\Delta
_\sigma }
\end{eqnarray*}
where the last isomorphism is given by Lemma \ref{decompose}. From the exact
sequence 
\[
1\rightarrow \Delta _\sigma \longrightarrow S_{n,\sigma }\longrightarrow
w_G\left( \sigma \right) \rightarrow 1 
\]
we conclude that 
\begin{equation}
\widehat{{\rm R}\left( G\right) _{\Lambda ,\sigma }}\simeq \left( \widehat{%
{\rm R}\left( {\rm C}_G\left( \sigma \right) \right) _{\Lambda ,\sigma }}%
\right) ^{w_G\left( \sigma \right) }.  \label{ok1}
\end{equation}
By {\bf \cite{SGA1}} Prop. 2.6, Exp. V, it is now enough to prove that the
stabilizers of geometric points (i.e. the inertia groups of points) of ${\rm %
Spec}\left( \widehat{{\rm R}\left( {\rm C}_G\left( \sigma \right) \right)
_{\Lambda ,\sigma }}\right) $ under the $w_G\left( \sigma \right) $-action
are trivial.

First of all, let us observe that ${\rm Spec}\left( \widetilde{{\rm R}}%
\left( \sigma \right) _\Lambda \right) $ is a closed subscheme of ${\rm Spec}%
\left( \widehat{{\rm R}\left( {\rm C}_G\left( \sigma \right) \right)
_{\Lambda ,\sigma }}\right) $. This can be seen as follows. It is obviously
enough to show that if $s\ $denotes the order of $\sigma $, the map 
\[
\pi _\sigma :{\rm R}\left( {\rm C}_G\left( \sigma \right) \right) _\Lambda
\longrightarrow {\rm R}\left( \sigma \right) _\Lambda =\frac{\Lambda \left[
t\right] }{\left( t^s-1\right) } 
\]
is surjective. First consider the case where $\sigma $ is contained in the
center of $G$. Since ${\rm R}\left( \sigma \right) _\Lambda $ is of finite
type over $\Lambda $, we may as well show that for any prime\footnote{%
Recall that $\sigma $ is essential hence $s$ is invertible in $\Lambda $.} $%
p\nmid s$ the induced map 
\[
\pi _{\sigma ,p}:{\rm R}\left( {\rm C}_G\left( \sigma \right) \right)
_\Lambda \otimes {\Bbb F}_p\longrightarrow {\rm R}\left( \sigma \right)
_\Lambda \otimes {\Bbb F}_p 
\]
is surjective. Note that if $E$ denote the standard $n$-dimensional
representation of $G$, $\pi _\sigma $ sends $\bigwedge^rE$ to $\binom nrt^r$%
. If $p\nmid n$, $\pi _{\sigma ,p}$ is surjective (in fact $\pi _\sigma
\left( E\right) =nt$ and $n$ is invertible in ${\Bbb F}_p$). If $p\mid n$,
let us write $n=qm$, with $q=p^i$ and $p\nmid m$. Since $\left( s,q\right)
=1 $, $t^q$ is a ring generator of ${\rm R}\left( \sigma \right) _\Lambda $
and to prove $\pi _{\sigma ,p}$ is injective it is enough to show that $%
p\nmid \binom nq$. This is elementary since the binomial expansion of 
\[
\left( 1+X\right) ^n=\left( 1+X^q\right) ^m 
\]
in ${\Bbb F}_p\left[ X\right] $, yields $\binom nq=m$ in ${\Bbb F}_p$. For a
general $\sigma \subseteq T$, let ${\rm C}_G\left( \sigma \right)
=\prod_{i=1}^l{\rm GL}_{d_i,k}$, where $\sum d_i=n$ and let $\sigma _i$
denote the image of $\sigma $ in ${\rm GL}_{d_i,k}$, $i=1,\ldots ,l$. Since $%
\sigma \subseteq \prod_{i=1}^l\sigma _i$ is an inclusion of diagonalizable
groups, the induced map 
\[
{\rm R}\left( \prod_{i=1}^l\sigma _i\right) =\bigotimes_{i=1}^l{\rm R}\left(
\sigma _i\right) \longrightarrow {\rm R}\left( \sigma \right) 
\]
is surjective (e.g. {\bf \cite{SGA3}}, tome II). But ${\rm R}\left( {\rm C}%
_G\left( \sigma \right) \right) _\Lambda \rightarrow {\rm R}\left( \sigma
\right) _\Lambda $ factors as 
\[
{\rm R}\left( {\rm C}_G\left( \sigma \right) \right) _\Lambda
=\bigotimes_{i=1}^l{\rm R}\left( {\rm GL}_{d_i,k}\right) _\Lambda
\longrightarrow \bigotimes_{i=1}^l{\rm R}\left( \sigma _i\right) _\Lambda
\longrightarrow {\rm R}\left( \sigma \right) _\Lambda 
\]
and also the first map is surjective (by the previous case, since $\sigma _i$
is contained in the center of ${\rm GL}_{d_i,k}$ and $\left| \sigma
_i\right| $ divides $\left| \sigma \right| $). This proves that ${\rm Spec}%
\left( \widetilde{{\rm R}}\left( \sigma \right) _\Lambda \right) $ is a
closed subscheme of ${\rm Spec}\left( \widehat{{\rm R}\left( {\rm C}_G\left(
\sigma \right) \right) _{\Lambda ,\sigma }}\right) $.

Since $\widehat{{\rm R}\left( {\rm C}_G\left( \sigma \right) \right)
_{\Lambda ,\sigma }}$ is the completion of ${\rm R}\left( {\rm C}_G\left(
\sigma \right) \right) _\Lambda $ along the ideal 
\[
\ker \left( {\rm R}\left( {\rm C}_G\left( \sigma \right) \right) _\Lambda
\rightarrow \widetilde{{\rm R}}\left( \sigma \right) _\Lambda \right) , 
\]
any nonempty closed subscheme of ${\rm Spec}\left( \widehat{{\rm R}\left( 
{\rm C}_G\left( \sigma \right) \right) _{\Lambda ,\sigma }}\right) $ meets
the closed subscheme ${\rm Spec}\left( \widetilde{{\rm R}}\left( \sigma
\right) _\Lambda \right) $. To prove that $w_G\left( \sigma \right) $ acts
freely on the geometric points of ${\rm Spec}\left( \widehat{{\rm R}\left( 
{\rm C}_G\left( \sigma \right) \right) _{\Lambda ,\sigma }}\right) $ it is
then enough to show that it acts freely on the geometric points of ${\rm Spec%
}\left( \widetilde{{\rm R}}\left( \sigma \right) _\Lambda \right) $.

Actually, more is true: the map $q:{\rm Spec}\left( \widetilde{{\rm R}}%
\left( \sigma \right) _\Lambda \right) \rightarrow {\rm Spec}\left( \Lambda
\right) $ is a $\left( {\Bbb Z}/s{\Bbb Z}\right) ^{*}$-torsor\footnote{%
Recall that the constant group scheme associated to $\left( {\Bbb Z}/s{\Bbb Z%
}\right) ^{*}$ is isomorphic to ${\rm Aut}_k\left( \sigma \right) $.}. In
fact, if ${\rm Spec}\left( \Omega \right) \rightarrow {\rm Spec}\left(
\Lambda \right) $ is a geometric point, the corresponding geometric fiber of 
$q$ is isomorphic to the spectrum of 
\[
\prod_{\alpha _i\in \widetilde{{\Bbb \mu }}_s\left( \Omega \right) }\frac{%
\Omega \left[ t\right] }{\left( t-\alpha _i\right) }\simeq \prod_{\alpha
_i\in \widetilde{{\Bbb \mu }}_s\left( \Omega \right) }\Omega 
\]
and $\left( {\Bbb Z}/s{\Bbb Z}\right) ^{*}$ acts by permutation on the
primitive roots $\widetilde{{\Bbb \mu }}_s\left( \Omega \right) $, by $%
\alpha \mapsto \alpha ^k$, $\left( k,s\right) =1$. In particular, the action
of the subgroup $w_G\left( \sigma \right) \subset \left( {\Bbb Z}/s{\Bbb Z}%
\right) ^{*}$ on ${\rm Spec}\left( \widetilde{{\rm R}}\left( \sigma \right)
_\Lambda \right) $ is free. \TeXButton{End Proof}{\endproof}

\smallskip\ 

\begin{proposition}
\label{firststep}The canonical morphism 
\[
K_{*}^{\prime }\left( X,G\right) _\Lambda \longrightarrow \prod_{\sigma \in 
{\cal C}\left( G\right) }\left( K_{*}^{\prime }\left( X^\sigma ,{\rm C}%
_G\left( \sigma \right) \right) _\sigma \right) ^{w_G\left( \sigma \right) }
\]
is an isomorphism.
\end{proposition}

\TeXButton{Proof}{\proof} By Lemma \ref{kurjev}, the canonical ring morphism 
\[
K_{*}^{\prime }\left( X,G\right) \otimes _{{\rm R}\left( G\right) }{\rm R}%
\left( T\right) \longrightarrow K_{*}^{\prime }\left( X,T\right) 
\]
is an isomorphism. Since ${\rm R}\left( G\right) \rightarrow {\rm R}\left(
T\right) $ is faithfully flat, it is enough to show that 
\[
K_{*}^{\prime }\left( X,T\right) _\Lambda \simeq K_{*}^{\prime }\left(
X,G\right) _\Lambda \otimes _{{\rm R}\left( G\right) _\Lambda }{\rm R}\left(
T\right) _\Lambda \longrightarrow \prod_{\sigma \in {\cal C}\left( G\right)
}\left( K_{*}^{\prime }\left( X^\sigma ,{\rm C}_G\left( \sigma \right)
\right) _\sigma \right) ^{w_G\left( \sigma \right) }\otimes _{{\rm R}\left(
G\right) _\Lambda }{\rm R}\left( T\right) _\Lambda 
\]
is an isomorphism. By Proposition \ref{concentration} (ii), we are left to
prove that 
\begin{equation}
\prod_{\sigma \in {\cal C}\left( G\right) }\left( K_{*}^{\prime }\left(
X^\sigma ,{\rm C}_G\left( \sigma \right) \right) _\sigma \right) ^{w_G\left(
\sigma \right) }\otimes _{{\rm R}\left( G\right) _\Lambda }{\rm R}\left(
T\right) _\Lambda \simeq \prod_{\Sb \sigma \text{ dual cyclic}  \\ \sigma
\subset T  \endSb }^{}K_{*}^{\prime }\left( X,T\right) _\sigma
\label{reduction}
\end{equation}
For any $\tau \in {\cal C}\left( G\right) $ ($\tau \subseteq T$, as usual),
we have 
\[
K_{*}^{\prime }\left( X^\tau ,{\rm C}_G\left( \tau \right) \right) _\tau
\otimes _{\widehat{{\rm R}\left( G\right) _{\Lambda ,\tau }}}\widehat{{\rm R}%
\left( T\right) _{\Lambda ,\tau }}\simeq \left( K_{*}^{\prime }\left( X^\tau
,{\rm C}_G\left( \tau \right) \right) _\tau \otimes _{\widehat{{\rm R}\left( 
{\rm C}_G\left( \tau \right) \right) _{\Lambda ,\tau }}}\widehat{{\rm R}%
\left( {\rm C}_G\left( \tau \right) \right) _{\Lambda ,\tau }}\right)
\otimes _{\widehat{{\rm R}\left( G\right) _{\Lambda ,\tau }}}\widehat{{\rm R}%
\left( T\right) _{\Lambda ,\tau }} 
\]
\[
\simeq \left( K_{*}^{\prime }\left( X^\tau ,{\rm C}_G\left( \tau \right)
\right) _\tau \otimes _{\widehat{R\left( G\right) _{\Lambda ,\tau }}}%
\widehat{{\rm R}\left( {\rm C}_G\left( \tau \right) \right) _{\Lambda ,\tau }%
}\right) \otimes _{\widehat{{\rm R}\left( {\rm C}_G\left( \tau \right)
\right) _{\Lambda ,\tau }}}\widehat{{\rm R}\left( T\right) _{\Lambda ,\tau }}
\]
By Lemma \ref{galois}, for any $\widehat{{\rm R}\left( G\right) _{\Lambda
,\tau }}$-module $M$, we have 
\[
M\otimes _{\widehat{{\rm R}\left( G\right) _{\Lambda ,\tau }}}\widehat{{\rm R%
}\left( {\rm C}_G\left( \tau \right) \right) _{\Lambda ,\tau }}\simeq
w_G\left( \tau \right) \times M 
\]
since a torsor is trivial when base changed along itself. Therefore 
\begin{equation}
K_{*}^{\prime }\left( X^\tau ,{\rm C}_G\left( \tau \right) \right) _\tau
\otimes _{\widehat{{\rm R}\left( G\right) _{\Lambda ,\tau }}}\widehat{{\rm R}%
\left( T\right) _{\Lambda ,\tau }}\simeq w_G\left( \tau \right) \times
\left( K_{*}^{\prime }\left( X^\tau ,{\rm C}_G\left( \tau \right) \right)
_\tau \otimes _{\widehat{{\rm R}\left( {\rm C}_G\left( \tau \right) \right)
_{\Lambda ,\tau }}}\widehat{{\rm R}\left( T\right) _{\Lambda ,\tau }}\right)
\label{lip}
\end{equation}
with $w_G\left( \tau \right) $ acting on l.h.s. by left multiplication on $%
w_G\left( \tau \right) $. Applying Lemma \ref{lemma?} (ii), to the l.h.s.,
we get 
\[
K_{*}^{\prime }\left( X^\tau ,{\rm C}_G\left( \tau \right) \right) _\tau
\otimes _{\widehat{{\rm R}\left( G\right) _{\Lambda ,\tau }}}\widehat{{\rm R}%
\left( T\right) _{\Lambda ,\tau }}\simeq w_G\left( \tau \right) \times
K_{*}^{\prime }\left( X^\tau ,T\right) _\tau 
\]
and taking invariants with respect to $w_G\left( \tau \right) $, 
\begin{equation}
\left( K_{*}^{\prime }\left( X^\tau ,{\rm C}_G\left( \tau \right) \right)
_\tau \otimes _{\widehat{{\rm R}\left( G\right) _{\Lambda ,\tau }}}\widehat{%
{\rm R}\left( T\right) _{\Lambda ,\tau }}\right) ^{w_G\left( \tau \right)
}\simeq K_{*}^{\prime }\left( X^\tau ,T\right) _\tau  \label{buonpunto}
\end{equation}
Comparing (\ref{reduction}) to (\ref{buonpunto}), we are reduced to proving
that for any $\sigma \in {\cal C}\left( G\right) $ there is an isomorphism 
\[
\left( K_{*}^{\prime }\left( X^\sigma ,{\rm C}_G\left( \sigma \right)
\right) _\sigma \right) ^{w_G\left( \sigma \right) }\otimes _{{\rm R}\left(
G\right) _\Lambda }{\rm R}\left( T\right) _\Lambda \simeq \prod_{\Sb \tau 
\text{ dual cyclic}  \\ ^{\tau \approx _{S_n}\sigma }  \endSb }\left(
K_{*}^{\prime }\left( X^\tau ,{\rm C}_G\left( \tau \right) \right) _\tau
\otimes _{\widehat{{\rm R}\left( G\right) _{\Lambda ,\tau }}}\widehat{{\rm R}%
\left( T\right) _{\Lambda ,\tau }}\right) ^{w_G\left( \tau \right) } 
\]
Since $\widehat{{\rm R}\left( T\right) _{\Lambda ,\tau }}$ is flat over $%
\widehat{{\rm R}\left( G\right) _{\Lambda ,\tau }}$ and $w_G\left( \tau
\right) $ acts trivially on it, we have ({\bf \cite{SGA1}}) 
\[
\left( K_{*}^{\prime }\left( X^\tau ,{\rm C}_G\left( \tau \right) \right)
_\tau \otimes _{\widehat{{\rm R}\left( G\right) _{\Lambda ,\tau }}}\widehat{%
{\rm R}\left( T\right) _{\Lambda ,\tau }}\right) ^{w_G\left( \tau \right)
}\simeq \left( K_{*}^{\prime }\left( X^\tau ,{\rm C}_G\left( \tau \right)
\right) _{\Lambda ,\tau }\right) ^{w_G\left( \tau \right) }\otimes _{%
\widehat{{\rm R}\left( G\right) _{\Lambda ,\tau }}}\widehat{{\rm R}\left(
T\right) _{\Lambda ,\tau }}. 
\]
By Lemma \ref{decompose}, we have isomorphisms 
\[
\left( K_{*}^{\prime }\left( X^\sigma ,{\rm C}_G\left( \sigma \right)
\right) _\sigma \right) ^{w_G\left( \sigma \right) }\otimes _{{\rm R}\left(
G\right) _\Lambda }{\rm R}\left( T\right) _\Lambda 
\]
\[
\simeq \left( K_{*}^{\prime }\left( X^\sigma ,{\rm C}_G\left( \sigma \right)
\right) _\sigma \right) ^{w_G\left( \sigma \right) }\otimes _{\widehat{{\rm R%
}\left( G\right) _{\Lambda ,\sigma }}}\widehat{{\rm R}\left( G\right)
_{\Lambda ,\sigma }}\otimes _{{\rm R}\left( G\right) _\Lambda }{\rm R}\left(
T\right) _\Lambda 
\]
\[
\simeq \prod_{\Sb \tau \text{ dual cyclic}  \\ ^{\tau \approx _{S_n}\sigma } 
\endSb }\left( K_{*}^{\prime }\left( X^\sigma ,{\rm C}_G\left( \sigma
\right) \right) _\sigma \right) ^{w_G\left( \sigma \right) }\otimes _{%
\widehat{{\rm R}\left( G\right) _{\Lambda ,\sigma }}}\widehat{{\rm R}\left(
T\right) _{\Lambda ,\tau }} 
\]
(recall that $\widehat{{\rm R}\left( G\right) _{\Lambda ,\sigma }}=\widehat{%
{\rm R}\left( G\right) _{\Lambda ,\tau }}$ for any $\tau \approx
_{S_n}\sigma $, since ${\frak m}_\sigma ^G={\frak m}_\tau ^G$). For each $%
\tau $, choosing an element $g\in S_n$ such that $g\sigma g^{-1}=\tau $,
determines an isomorphism $K_{*}^{\prime }\left( X^\sigma ,{\rm C}_G\left(
\sigma \right) \right) _\sigma \simeq K_{*}^{\prime }\left( X^\tau ,{\rm C}%
_G\left( \tau \right) \right) _\tau $ whose restriction to invariants 
\[
\left( K_{*}^{\prime }\left( X^\sigma ,{\rm C}_G\left( \sigma \right)
\right) _\sigma \right) ^{w_G\left( \sigma \right) }\simeq \left(
K_{*}^{\prime }\left( X^\tau ,{\rm C}_G\left( \tau \right) \right) _\tau
\right) ^{w_G\left( \tau \right) } 
\]
is independent on the choice of $g$. Therefore we have a canonical
isomorphism 
\[
\left( K_{*}^{\prime }\left( X^\sigma ,{\rm C}_G\left( \sigma \right)
\right) _\sigma \right) ^{w_G\left( \sigma \right) }\otimes _{{\rm R}\left(
G\right) _\Lambda }{\rm R}\left( T\right) _\Lambda \simeq \prod_{\Sb \tau 
\text{ dual cyclic}  \\ ^{\tau \approx _{S_n}\sigma }  \endSb }\left(
K_{*}^{\prime }\left( X^\sigma ,{\rm C}_G\left( \sigma \right) \right)
_\sigma \right) ^{w_G\left( \sigma \right) }\otimes _{\widehat{{\rm R}\left(
G\right) _{\Lambda ,\sigma }}}\widehat{{\rm R}\left( T\right) _{\Lambda
,\tau }} 
\]
\[
\simeq \prod_{\Sb \tau \text{ dual cyclic}  \\ ^{\tau \approx _{S_n}\sigma } 
\endSb }\left( K_{*}^{\prime }\left( X^\tau ,{\rm C}_G\left( \tau \right)
\right) _\tau \right) ^{w_G\left( \tau \right) }\otimes _{\widehat{{\rm R}%
\left( G\right) _{\Lambda ,\tau }}}\widehat{{\rm R}\left( T\right) _{\Lambda
,\tau }} 
\]
as desired.\TeXButton{End Proof}{\endproof}

Since $K_{*}\left( X,G\right) \simeq K_{*}^{\prime }\left( X,G\right) $ and $%
K_{*}\left( X^\sigma ,{\rm C}_G\left( \sigma \right) \right) \simeq
K_{*}^{\prime }\left( X^\sigma ,{\rm C}_G\left( \sigma \right) \right) $,
comparing Proposition \ref{firststep} with (\ref{mainGLn}), we see that the
proof of Theorem \ref{GLn} can be completed by the following

\begin{proposition}
\label{pincopallino}For any $\sigma \in {\cal C}\left( G\right) $, the
morphism given by Lemma \ref{construction} and induced by the product ${\rm C%
}_G\left( \sigma \right) \times \sigma \rightarrow {\rm C}_G\left( \sigma
\right) $ 
\[
\theta _{{\rm C}_G\left( \sigma \right) ,\sigma }:K_{*}^{\prime }\left(
X^\sigma ,{\rm C}_G\left( \sigma \right) \right) _\sigma \longrightarrow
K_{*}^{\prime }\left( X^\sigma ,{\rm C}_G\left( \sigma \right) \right) _{%
{\rm geom}}\otimes \widetilde{{\rm R}}\left( \sigma \right) _\Lambda 
\]
is an isomorphism.
\end{proposition}

\TeXButton{Proof}{\proof} To simplify the notation we write $\theta _\sigma $
for $\theta _{{\rm C}_G\left( \sigma \right) ,\sigma }$. As usual, we may
suppose $\sigma $ contained in $T$. Since ${\rm C}_G\left( \sigma \right) $
is isomorphic to a product of general linear groups over $k$, we can take $T$
as its maximal torus and by Lemma \ref{kurjev}, the canonical ring morphism 
\[
K_{*}^{\prime }\left( X,{\rm C}_G\left( \sigma \right) \right) \otimes _{%
{\rm R}\left( {\rm C}_G\left( \sigma \right) \right) }{\rm R}\left( T\right)
\longrightarrow K_{*}^{\prime }\left( X,T\right) 
\]
is an isomorphism. Moreover, ${\rm R}\left( {\rm C}_G\left( \sigma \right)
\right) \rightarrow {\rm R}\left( T\right) $ being faithfully flat, it is
enough to prove that $\theta _\sigma \otimes $ ${\rm id}_{{\rm R}\left(
T\right) }$ is an isomorphism. To prove this, let us consider the
commutative diagram 
\[
\begin{tabular}[t]{ccc}
$K_{*}^{\prime }\left( X^\sigma ,{\rm C}_G\left( \sigma \right) \right)
_\sigma \otimes _{{\rm R}\left( {\rm C}_G\left( \sigma \right) \right)
_\Lambda }{\rm R}\left( T\right) _\Lambda $ & $\stackrel{\theta _\sigma
\otimes {\rm id}}{\longrightarrow }$ & $\left( K_{*}^{\prime }\left( X,{\rm C%
}_G\left( \sigma \right) \right) _{{\rm geom}}\otimes \widetilde{{\rm R}}%
\left( \sigma \right) _\Lambda \right) \otimes _{{\rm R}\left( {\rm C}%
_G\left( \sigma \right) \right) _\Lambda }{\rm R}\left( T\right) _\Lambda $
\\ 
$^{\omega _\sigma }\downarrow $ &  & $\downarrow ^{\widetilde{\gamma _\sigma 
}}$ \\ 
$K_{*}^{\prime }\left( X^\sigma ,T\right) _\sigma $ & $\stackunder{\theta
_{T,\sigma }}{\longrightarrow }$ & $K_{*}^{\prime }\left( X^\sigma ,T\right)
_{{\rm geom}}\otimes \widetilde{{\rm R}}\left( \sigma \right) _\Lambda $%
\end{tabular}
\]
where :

- $K_{*}\left( X^\sigma ,{\rm C}_G\left( \sigma \right) \right) _{{\rm geom}%
}\otimes \widetilde{{\rm R}}\left( \sigma \right) _\Lambda $ is an ${\rm 
\QTR{mathrm}{R}}\left( {\rm C}_G\left( \sigma \right) \right) _\Lambda $%
-module via the coproduct ring morphism $\Delta _{{\rm C}_G\left( \sigma
\right) }:{\rm R}\left( {\rm C}_G\left( \sigma \right) \right) _\Lambda
\rightarrow {\rm R}\left( {\rm C}_G\left( \sigma \right) \right) _\Lambda
\otimes \widetilde{{\rm R}}\left( \sigma \right) _\Lambda $ (induced by the
product ${\rm C}_G\left( \sigma \right) \times \sigma \rightarrow {\rm C}%
_G\left( \sigma \right) $);

- $\omega _\sigma $ is the canonical map induced by the inclusion $%
T\hookrightarrow {\rm C}_G\left( \sigma \right) $ and is an isomorphism by
Lemma \ref{lemma?};

- $\theta _{T,\sigma }$ is an isomorphism as shown in the proof of Theorem 
\ref{torus};

- $\widetilde{\gamma _\sigma }$ sends $\left( x\otimes u\right) \otimes t$
to $\left( \Delta _T\left( t\right) \cdot x_{\mid T}\right) \otimes u$, for $%
x\in K_{*}\left( X^\sigma ,{\rm C}_G\left( \sigma \right) \right) _{{\rm geom%
}}$, $u\in \widetilde{{\rm R}}\left( \sigma \right) _\Lambda $, $t\in {\rm R}%
\left( T\right) _\Lambda $, $\Delta _T:{\rm R}\left( T\right) _\Lambda
\rightarrow {\rm R}\left( T\right) _\Lambda \otimes \widetilde{{\rm R}}%
\left( \sigma \right) _\Lambda $ being the coproduct induced by the product $%
T\times \sigma \rightarrow T$.

So we are left to prove that $\widetilde{\gamma _\sigma }$ is an isomorphism.

First of all, let us observe that if $R$ is a ring, $A\rightarrow A^{\prime
} $ a ring morphism and $M$ an $A$-module, there is a natural isomorphism 
\[
\left( M\otimes _{{\Bbb Z}}R\right) \otimes _{A\otimes _{{\Bbb Z}}B}\left(
A^{\prime }\otimes _{{\Bbb Z}}R\right) \longrightarrow \left( M\otimes
_AA^{\prime }\right) \otimes _{{\Bbb Z}}R 
\]
\[
\left( m\otimes r_1\right) \otimes \left( a^{\prime }\otimes r_2\right)
\longmapsto \left( m\otimes a^{\prime }\right) \otimes r_1r_2. 
\]
Applying this to $M=K_{*}\left( X^\sigma ,{\rm C}_G\left( \sigma \right)
\right) _{{\rm geom}}$, $R=\widetilde{{\rm R}}\left( \sigma \right) _\Lambda 
$, $A={\rm R}\left( {\rm C}_G\left( \sigma \right) \right) _\Lambda $, $%
A^{\prime }={\rm R}\left( T\right) _\Lambda $ and using Lemma \ref{lemma?},
we get a canonical isomorphism 
\begin{equation}
{\rm f}:K_{*}^{\prime }\left( X^\sigma ,T\right) _{{\rm geom}}\otimes 
\widetilde{{\rm R}}\left( \sigma \right) _\Lambda \longrightarrow \left(
K_{*}^{\prime }\left( X^\sigma ,{\rm C}_G\left( \sigma \right) \right) _{%
{\rm geom}}\otimes \widetilde{{\rm R}}\left( \sigma \right) _\Lambda \right)
\otimes _{{\rm R}\left( {\rm C}_G\left( \sigma \right) \right) _\Lambda
\otimes \widetilde{{\rm R}}\left( \sigma \right) _\Lambda }\left( {\rm R}%
\left( T\right) _\Lambda \otimes \widetilde{{\rm R}}\left( \sigma \right)
_\Lambda \right) ^{\prime }.  \label{claim1}
\end{equation}
where we have denoted by $\left( {\rm R}\left( T\right) _\Lambda \otimes 
\widetilde{{\rm R}}\left( \sigma \right) _\Lambda \right) ^{\prime }$ the $%
{\rm R}\left( {\rm C}_G\left( \sigma \right) \right) _\Lambda \otimes 
\widetilde{{\rm R}}\left( \sigma \right) _\Lambda $-algebra 
\[
{\rm res}\otimes {\rm id}:{\rm R}\left( {\rm C}_G\left( \sigma \right)
\right) _\Lambda \otimes \widetilde{{\rm R}}\left( \sigma \right) _\Lambda
\longrightarrow {\rm R}\left( T\right) _\Lambda \otimes \widetilde{{\rm R}}%
\left( \sigma \right) _\Lambda \text{.} 
\]
It is an elementary fact that there are mutually inverse isomorphisms $%
\alpha _{{\rm C}_G\left( \sigma \right) }$, $\beta _{{\rm C}_G\left( \sigma
\right) }$ and $\alpha _T$, $\beta _T$ fitting into the commutative diagrams

\medskip\ 

\begin{equation}
\text{\TeXButton{Ctriangle}
{\settriparms[-1`1`1;500]
\Ctriangle[\QTR{rm}{R}\left( \QTR{rm}{C}_G\left( \sigma \right) \right) _\Lambda \otimes \widetilde{\QTR{rm}{R}}\left( \sigma \right) _\Lambda  `\QTR{rm}{R}\left( \QTR{rm}{C}_G\left( \sigma \right) \right) _\Lambda  `\QTR{rm}{R}\left( \QTR{rm}{C}_G\left( \sigma \right) \right) _\Lambda \otimes \widetilde{\QTR{rm}{R}}\left( \sigma \right) _\Lambda  ;\Delta _{C_G\left( \sigma \right) }`\alpha _{C_G\left( \sigma \right) } `\QTR{rm}{id}\otimes 1 ]}%
}  \label{cisigma}
\end{equation}

\bigskip\ 

\medskip\ \ 
\begin{equation}
\text{\TeXButton{Ctriangle}
{\settriparms[-1`1`1;400]
\Ctriangle[\QTR{rm}{R}\left( T\right) _\Lambda \otimes \widetilde{\QTR{rm}{R}}\left( \sigma \right) _\Lambda `\QTR{rm}{R}\left( T\right) _\Lambda  `\QTR{rm}{R}\left( T\right) _\Lambda \otimes \widetilde{\QTR{rm}{R}}\left( \sigma \right) _\Lambda  ;\Delta _T`\alpha _T `\QTR{rm}{id}\otimes 1 ]}%
}  \label{ti}
\end{equation}
\medskip\ 

\noindent and compatible with restriction maps induced by $T\hookrightarrow 
{\rm C}_G\left( \sigma \right) $. This is exactly the dual assertion to the
general fact that ``an action $H\times Y\rightarrow Y$ is isomorphic over $X$
to the projection on the second factor ${\rm pr}_2:H\times Y\rightarrow Y$
'', for any group scheme $H$ and any algebraic space $Y$. From (\ref{cisigma}%
) we get an isomorphism 
\[
\widetilde{\alpha }:\left( {\rm R}\left( {\rm C}_G\left( \sigma \right)
\right) _\Lambda \otimes \widetilde{{\rm R}}\left( \sigma \right) _\Lambda
\right) ^{\prime }\otimes _{{\rm R}\left( {\rm C}_G\left( \sigma \right)
\right) _\Lambda }{\rm R}\left( T\right) _\Lambda \longrightarrow {\rm R}%
\left( T\right) _\Lambda \otimes \widetilde{{\rm R}}\left( \sigma \right)
_\Lambda 
\]
where $\left( {\rm R}\left( {\rm C}_G\left( \sigma \right) \right) _\Lambda
\otimes \widetilde{{\rm R}}\left( \sigma \right) _\Lambda \right) ^{\prime }$
denotes ${\rm R}\left( {\rm C}_G\left( \sigma \right) \right) _\Lambda $%
-algebra 
\[
\Delta _{{\rm C}_G\left( \sigma \right) }:{\rm R}\left( {\rm C}_G\left(
\sigma \right) \right) _\Lambda \rightarrow {\rm R}\left( {\rm C}_G\left(
\sigma \right) \right) _\Lambda \otimes \widetilde{{\rm R}}\left( \sigma
\right) _\Lambda . 
\]
Therefore, if we denote by $\left( {\rm R}\left( T\right) _\Lambda \otimes 
\widetilde{{\rm R}}\left( \sigma \right) _\Lambda \right) ^{\prime \prime }$
the ${\rm R}\left( {\rm C}_G\left( \sigma \right) \right) _\Lambda \otimes 
\widetilde{{\rm R}}\left( \sigma \right) _\Lambda $-algebra 
\[
\left( {\rm res}\otimes {\rm id}\right) \circ \alpha _T:{\rm R}\left( {\rm C}%
_G\left( \sigma \right) \right) _\Lambda \otimes \widetilde{{\rm R}}\left(
\sigma \right) _\Lambda \longrightarrow {\rm R}\left( T\right) _\Lambda
\otimes \widetilde{{\rm R}}\left( \sigma \right) _\Lambda \text{ ,} 
\]
the composition 
\[
\left( K_{*}^{\prime }\left( X^\sigma ,{\rm C}_G\left( \sigma \right)
\right) _{{\rm geom}}\otimes \widetilde{{\rm R}}\left( \sigma \right)
_\Lambda \right) ^{\prime }\otimes _{{\rm R}\left( {\rm C}_G\left( \sigma
\right) \right) _\Lambda }{\rm R}\left( T\right) _\Lambda 
\]
\[
=\left( K_{*}^{\prime }\left( X^\sigma ,{\rm C}_G\left( \sigma \right)
\right) _{{\rm geom}}\otimes \widetilde{{\rm R}}\left( \sigma \right)
_\Lambda \right) \tbigotimes\limits_{{\rm R}\left( {\rm C}_G\left( \sigma
\right) \right) _\Lambda \otimes \widetilde{{\rm R}}\left( \sigma \right)
_\Lambda }\left( \left( {\rm R}\left( {\rm C}_G\left( \sigma \right) \right)
_\Lambda \otimes \widetilde{{\rm R}}\left( \sigma \right) _\Lambda \right)
^{\prime }\otimes _{{\rm R}\left( {\rm C}_G\left( \sigma \right) \right)
_\Lambda }{\rm R}\left( T\right) _\Lambda \right) 
\]
\[
\stackrel{{\rm id}\otimes \widetilde{\alpha }}{\longrightarrow }\left(
K_{*}^{\prime }\left( X^\sigma ,{\rm C}_G\left( \sigma \right) \right) _{%
{\rm geom}}\otimes \widetilde{{\rm R}}\left( \sigma \right) _\Lambda \right)
\tbigotimes\limits_{{\rm R}\left( {\rm C}_G\left( \sigma \right) \right)
_\Lambda \otimes \widetilde{{\rm R}}\left( \sigma \right) _\Lambda }\text{ }%
\left( {\rm R}\left( T\right) _\Lambda \otimes \widetilde{{\rm R}}\left(
\sigma \right) _\Lambda \right) ^{\prime \prime }\rightarrow 
\]
\[
\stackrel{{\rm id}\otimes \beta _T}{\longrightarrow }\left( K_{*}^{\prime
}\left( X^\sigma ,{\rm C}_G\left( \sigma \right) \right) _{{\rm geom}%
}\otimes \widetilde{{\rm R}}\left( \sigma \right) _\Lambda \right)
\tbigotimes\limits_{{\rm R}\left( {\rm C}_G\left( \sigma \right) \right)
_\Lambda \otimes \widetilde{{\rm R}}\left( \sigma \right) _\Lambda }\text{ }%
\left( {\rm R}\left( T\right) _\Lambda \otimes \widetilde{{\rm R}}\left(
\sigma \right) _\Lambda \right) ^{\prime } 
\]
\[
\simeq K_{*}^{\prime }\left( X^\sigma ,T\right) _{{\rm geom}}\otimes 
\widetilde{{\rm R}}\left( \sigma \right) _\Lambda 
\]
is an isomorphism and it can be easily checked to coincide with $\widetilde{%
\gamma _\sigma }$. \TeXButton{End Proof}{\endproof}

\section{The main theorem: the general case}

In this Section, we use Theorem \ref{GLn} to deduce the same result for the
action of a linear algebraic $k$-group $G$, having finite stabilizers, on a
regular separated noetherian $k$-algebraic space $X$. We will write $\Lambda 
$ for $\Lambda _{\left( G,X\right) }$.

We start with a general fact

\begin{proposition}
\label{nilpotent}Let $X$ be a regular noetherian separated $k$-algebraic
space on which a linear algebraic $k$-group $G$ acts with finite
stabilizers. Then there exists an integer $N>0$ such that if $a_1,\ldots
,a_N\in K_0\left( X,G\right) _{{\rm geom}}$ have rank zero on each connected
component of $X$, then the multiplication by $\prod_{i=1}^Na_i$ on $%
K_{*}^{^{\prime }}\left( X,G\right) _{{\rm geom}}$ is zero.
\end{proposition}

In particular

\begin{corollary}
\label{nilpotent0}Let $X$ be a regular noetherian separated $k$-algebraic
space with a connected action of a linear algebraic $k$-group $G$ having
finite stabilizers. Then the geometric localization 
\[
{\rm rk}_{{\rm 0,geom}}:K_0\left( X,G\right) _{{\rm geom}}\longrightarrow
\Lambda 
\]
of the rank morphism has a nilpotent kernel.
\end{corollary}

{\bf Proof of Prop. \ref{nilpotent}.} Let us choose a closed immersion $%
G\hookrightarrow {\rm GL}_{n,k}$ (for some $n>0$). By Morita equivalence, 
\[
K_{*}^{\prime }\left( X,G\right) \simeq K_{*}^{\prime }\left( X\times ^G{\rm %
GL}_{n,k},{\rm GL}_{n,k}\right) 
\]
and 
\[
K_0\left( X,G\right) \simeq K_0\left( X\times ^G{\rm GL}_{n,k},{\rm GL}%
_{n,k}\right) \text{.} 
\]
Moreover, $\Lambda _{(X\times ^G{\rm GL}_{n,k},{\rm GL}_{n,k})}=\Lambda $.
Let $\xi =x/s\in K_{*}^{\prime }\left( X,G\right) _{{\rm geom}}$, with $x\in
K_{*}^{\prime }\left( X,G\right) _\Lambda $ and $s\in {\rm rk}^{-1}\left(
1\right) $ where ${\rm rk}:{\rm R}\left( G\right) \rightarrow \Lambda $ is
the rank morphism, $a_i=\alpha _i/s_i$, with $\alpha _i\in K_0\left(
X,G\right) _\Lambda $ and $s_i\in {\rm rk}^{-1}\left( 1\right) $ for $%
i=1,...,N$. Let us consider the elements $x/1$ in $K_{*}^{\prime }\left(
X\times ^G{\rm GL}_{n,k},{\rm GL}_{n,k}\right) _{{\rm geom}}$, $\alpha _i/1$
in $K_0\left( X\times ^G{\rm GL}_{n,k},{\rm GL}_{n,k}\right) _{{\rm geom}}$
for $i=1,\ldots ,N$. Since the canonical homomorphism 
\[
K_{*}^{\prime }\left( X\times ^G{\rm GL}_{n,k},{\rm GL}_{n,k}\right) _{{\rm %
geom}}\longrightarrow K_{*}^{\prime }\left( X,G\right) _{{\rm geom}} 
\]
is a morphism of modules over the ring morphism 
\[
K_0\left( X\times ^G{\rm GL}_{n,k},{\rm GL}_{n,k}\right) _{{\rm geom}%
}\longrightarrow K_0\left( X,G\right) _{{\rm geom}}, 
\]
if the theorem holds for $G={\rm GL}_{n,k}$ and $N$ is the corresponding
integer, the product $\prod_i\alpha _i/1$ in $K_0\left( X,G\right) _{{\rm %
geom}}$ annihilates $x/1\in K_{*}^{\prime }\left( X,G\right) _{{\rm geom}}$%
and a fortiori $\prod_ia_i$ annihilates $\xi $ in $K_{*}^{\prime }\left(
X,G\right) _{{\rm geom}}$. So, we may assume $G={\rm GL}_{n,k}$. Let $T$ be
the maximal torus of diagonal matrices in $G$. By Lemma \ref{lemma?} (i)
with $\sigma =1$, there are isomorphisms 
\[
\omega _{1,{\rm geom}}:K_0\left( X,{\rm GL}_{n,k}\right) _{{\rm geom}%
}\otimes _{{\rm R}\left( {\rm GL}_{n,k}\right) _\Lambda }{\rm R}\left(
T\right) _\Lambda \simeq K_0\left( X,T\right) _{{\rm geom}} 
\]
\[
K_{*}^{\prime }\left( X,{\rm GL}_{n,k}\right) _{{\rm geom}}\otimes _{{\rm R}%
\left( {\rm GL}_{n,k}\right) _\Lambda }{\rm R}\left( T\right) _\Lambda
\simeq K_{*}^{\prime }\left( X,T\right) _{{\rm geom}}. 
\]
Since ${\rm R}\left( {\rm GL}_{n,k}\right) \rightarrow {\rm R}\left(
T\right) $ is faithfully flat and the diagram 
\[
\begin{tabular}{ccc}
$K_0\left( X,{\rm GL}_{n,k}\right) _{{\rm geom}}\otimes _{{\rm R}\left( {\rm %
GL}_{n,k}\right) _\Lambda }{\rm R}\left( T\right) _\Lambda $ & $\stackrel{%
{\rm rk}_{{\rm 0,geom}}\otimes {\rm id}}{\longrightarrow }$ & $\Lambda
\otimes _{{\rm R}\left( {\rm GL}_{n,k}\right) }{\rm R}\left( T\right)
_\Lambda $ \\ 
$^{\omega _{1,{\rm geom}}}\downarrow $ &  & $\quad \qquad \downarrow ^{{\rm %
id}_\Lambda \otimes {\rm rk}_T}$ \\ 
$\qquad \qquad K_0\left( X,T\right) _{{\rm geom}}$ & $\stackunder{{\rm rk}_{%
{\rm 0,geom}}}{\longrightarrow }$ & $\Lambda $%
\end{tabular}
\]
commutes, we reduce ourselves to proving the proposition for $G=T$, a split
torus.

To handle this case, we proceed by noetherian induction on $X$. By {\bf \cite
{Th2}}, Prop. 4.10, there exists a $T$-invariant nonempty open subscheme $%
j:U\hookrightarrow X$, a closed diagonalizable subgroup $T^{\prime }$ of $T$
and a $T$-equivariant isomorphism 
\[
U\simeq T/T^{\prime }\times \left( U/T\right) . 
\]
Since $U$ is nonempty and $T$ acts on $X$ with finite stabilizers, $%
T^{\prime }$ is finite over $k$ and $K_{*}^{\prime }\left( U,T\right) \simeq
K_{*}^{\prime }\left( U/T\right) \otimes _{R\left( T\right) }R\left(
T^{\prime }\right) $, by Morita equivalence theorem ({\bf \cite{Th3}},
Proposition 6.2). Let $i:Z\hookrightarrow X$ be the closed complement of $U$
in $X$ and $N^{\prime }$ an integer satisfying the proposition for both $Z$
and $U$. Consider the geometric localization sequence 
\[
K_{*}^{\prime }\left( Z,T\right) _{{\rm geom}}\stackrel{i_{*}}{%
\longrightarrow }K_{*}^{\prime }\left( X,T\right) _{{\rm geom}}\stackrel{%
j^{*}}{\longrightarrow }K_{*}^{\prime }\left( U,T\right) _{{\rm geom}} 
\]
and let $\xi \in K_{*}^{\prime }\left( X,T\right) _{{\rm geom}}$. Let $%
a_1,\ldots ,a_{2N^{\prime }}\in K_0\left( X,T\right) _{{\rm geom}}$. By our
choice of $N^{\prime }$, 
\[
j^{*}\left( a_{N^{\prime }+1}\cdot \cdots \cdot a_{2N^{\prime }}\cup \xi
\right) =0, 
\]
thus $a_{N^{\prime }+1}\cdot \cdots \cdot a_{2N^{\prime }}\cap \xi
=i_{*}\left( \eta \right) $ for some $\eta $ in $K_{*}^{\prime }\left(
Z,T\right) _{{\rm geom}}$. By projection formula we get 
\[
a_1\cdot \cdots \cdot a_{2N^{\prime }}\cup \xi =i_{*}\left( i^{*}\left(
a_1\right) \cdot \cdots \cdot i^{*}\left( a_{N^{\prime }}\right) \cup \eta
\right) 
\]
which is zero by our choice of $N^{\prime }$ and by the fact that rank
morphisms commutes with pullbacks. Thus, $N\doteq 2N^{\prime }$ satisfies
our proposition. \TeXButton{End Proof}{\endproof}

\begin{remark}
\label{stack} By Corollary \ref{nilpotent0}, $K_{*}\left( X,G\right) _{{\rm %
geom}}$ is isomorphic to the localization of $K_{*}\left( X,G\right)
_\Lambda $ at the multiplicative subset $\left( {\rm rk}_0\right)
^{-1}\left( 1\right) $, where ${\rm rk}_0:K_0\left( X,G\right) _\Lambda
\rightarrow \Lambda $ is the rank morphism. Therefore, if $X$ is regular, $%
K_{*}\left( X,G\right) _{{\rm geom}}$ depends only on the quotient stack $%
\left[ X/G\right] $ (\cite{L-MB}) and not on its presentation as a quotient.
\end{remark}

\smallskip\ 

The main theorem of this paper is:

\begin{theorem}
\label{main}Let $X$ be a noetherian regular separated algebraic space over a
field $k$ and $G$ a linear algebraic $k$-group with a sufficiently rational
action on $X$ having finite stabilizers. Suppose moreover that for any
essential dual cyclic $k$-subgroup scheme $\sigma \subseteq G$, the quotient
algebraic space $G/{\rm \QTR{mathrm}{C}}_G\left( \sigma \right) $ is smooth
over $k$ (which is the case if, e.g., $G$ is smooth or abelian). Then ${\cal %
C}\left( G\right) $ is finite and the map defined in (\ref{due}) 
\[
\Psi _{X,G}:K_{*}\left( X,G\right) _\Lambda \longrightarrow \prod_{\sigma
\in {\cal C}\left( G\right) }\left( K_{*}\left( X^\sigma ,{\rm C}\left(
\sigma \right) \right) _{{\rm geom}}\otimes \widetilde{{\rm R}}\left( \sigma
\right) _\Lambda \right) ^{w_G\left( \sigma \right) }
\]
is an isomorphism of ${\rm R}\left( G\right) $-algebras.
\end{theorem}

\smallskip\ 

\begin{remark}
In the next subsection we will also give less restrictive hypotheses on $G$
under which Theorem \ref{main} still holds.

Note also that if $X$ has the ''$G$-equivariant resolution property'' (i.e.
any $G$-equivariant coherent sheaf is the $G$-equivariant epimorphic image
of a $G$-equivariant locally free coherent sheaf) then in Theorem \ref{main}
one can replace our $K_{*}$ with Quillen $K$-theory of $G$-equivariant
locally free coherent sheaves. This happens, for example, if $X$ is a scheme
and $G$ is smooth or finite ({\bf \cite{Th3}}).
\end{remark}

\subsection{Proof of Theorem 5.4}

Let us choose, for some $n$, a closed immersion $G\hookrightarrow {\rm GL}%
_{n,k}$ and consider the algebraic space quotient 
\[
Y\doteq {\rm GL}_{n,k}\times ^GX\text{.} 
\]
We claim that if the theorem holds for $Y$ with the induced ${\rm GL}_{n,k}$%
-action then it holds for $X$ with the given $G$-action. First of all, let
us note that $Y$ is separated. The action map $\psi :G\times \left( {\rm GL}%
_{n,k}\times X\right) \rightarrow \left( {\rm GL}_{n,k}\times X\right)
\times \left( {\rm GL}_{n,k}\times X\right) $ is proper (hence a closed
immersion) since its composition with the separated map 
\[
p_{123}:{\rm GL}_{n,k}\times X\times {\rm GL}_{n,k}\times X\rightarrow {\rm %
GL}_{n,k}\times X\times {\rm GL}_{n,k} 
\]
(here we use that $X$ is separated) is just ${\rm id}_X\times a$ where $a$
is the action map of $G$ on ${\rm GL}_{n,k}$, hence is proper ({\bf \cite
{EGAI}}, Remarque 5.1.7, which obviously carries over to algebraic spaces).
Let $P\doteq X\times {\rm GL}_{n,k}$. In the cartesian diagram, 
\[
\begin{tabular}{ccc}
$P\times _YP$ & $\stackrel{j}{\longrightarrow }$ & $P\times P$ \\ 
$\downarrow $ &  & $\quad \downarrow ^{\pi \times \pi }$ \\ 
$Y$ & $\stackunder{\Delta _Y}{\longrightarrow }$ & $Y\times Y$%
\end{tabular}
\]
$P\times _YP\simeq G\times P$, since $\pi :P\rightarrow Y$ is a $G$-torsor,
and $\pi $ is faithfully flat; therefore $\Delta _Y$ is a closed immersion
i.e. $Y$ is separated.

Note that $\Lambda _{\left( Y,{\rm GL}_{n,k}\right) }=\Lambda $.

Consider the morphism defined in (\ref{due}) 
\begin{equation}
\Psi _{X,G}:K_{*}\left( X,G\right) _\Lambda \longrightarrow \prod_{\sigma
\in {\cal C}\left( G\right) }\left( K_{*}\left( X^\sigma ,{\rm C}_G\left(
\sigma \right) \right) _{{\rm geom}}\otimes \widetilde{{\rm R}}\left( \sigma
\right) _\Lambda \right) ^{w_G\left( \sigma \right) };  \label{duebis}
\end{equation}
By Theorem \ref{GLn}, the map 
\[
\Psi _{Y,{\rm GL}_{n,k}}:K_{*}\left( Y,{\rm GL}_{n,k}\right) _\Lambda
\longrightarrow \prod_{\rho \in {\cal C}({\rm GL}_{n,k})}\left( K_{*}\left(
Y^\rho ,{\rm C}_{{\rm GL}_{n,k}}\left( \rho \right) \right) _{{\rm geom}%
}\otimes \widetilde{{\rm R}}\left( \rho \right) _\Lambda \right) ^{w_{{\rm GL%
}_{n,k}}\left( \rho \right) } 
\]
is an isomorphism and by Morita equivalence theorem ({\bf \cite{Th3}},
Proposition 6.2) $K_{*}\left( Y,{\rm GL}_{n,k}\right) _\Lambda \simeq
K_{*}\left( X,G\right) _\Lambda $. We will prove the theorem by constructing
an isomorphism 
\begin{equation}
\prod_{\rho \in {\cal C}({\rm GL}_{n,k})}\left( K_{*}\left( Y^\rho ,{\rm C}%
\left( \rho \right) \right) _{{\rm geom}}\otimes \widetilde{{\rm R}}\left(
\rho \right) \right) ^{w_{{\rm GL}_{n,k}}(\rho )}\longrightarrow
\prod_{\sigma \in {\cal C}\left( G\right) }\left( K_{*}\left( X^\sigma ,{\rm %
C}_G\left( \sigma \right) \right) _{{\rm geom}}\otimes \widetilde{{\rm R}}%
\left( \sigma \right) _\Lambda \right) ^{w_G\left( \sigma \right) }
\label{phantom}
\end{equation}
commuting with the $\Psi $'s and Morita isomorphisms.

Let $\alpha :{\cal C}\left( G\right) \rightarrow {\cal C}\left( {\rm GL}%
_{n,k}\right) $ be the natural map. If $Y^\rho \neq \emptyset $, there
exists a dual cyclic subgroup $\sigma \subseteq G$, ${\rm GL}_{n.k}$%
-conjugate to $\rho $ (and $X^\sigma \neq \emptyset $); therefore $Y^\rho
=\emptyset $ unless $\rho \in {\rm im}\left( \alpha \right) $ and we may
restrict the first product in (\ref{phantom}) to those $\rho $ in the image
of $\alpha $ and suppose ${\rm im}\left( \alpha \right) \subseteq {\cal C}%
\left( G\right) $ as well. The following proposition describes the $Y^\rho $%
's which appear:

\begin{proposition}
\label{disjointunion}Let $X$ be a noetherian regular separated algebraic
space over a field $k$ and $G$ a linear algebraic $k$-group with a
sufficiently rational action on $X$ having finite stabilizers. Suppose
moreover that for any essential dual cyclic $k$-subgroup scheme $\sigma
\subseteq G$, the quotient algebraic space $G/{\rm \QTR{mathrm}{C}}_G\left(
\sigma \right) $ is smooth over $k$. Let $G\hookrightarrow {\rm GL}_{n,k}$ a
closed embedding and $\rho \in {\rm im}\left( \alpha \right) $ an essential
dual cyclic subgroup and $Y\doteq {\rm GL}_{n,k}\times ^GX$ the algebraic
space quotient for the left diagonal action of $G$. If ${\cal C}_{{\rm GL}%
_{n,k},G}\left( \rho \right) \subseteq {\cal C}\left( G\right) $ denotes the
fiber $\alpha ^{-1}\left( \rho \right) $ , then:

(i) choosing for each $\sigma \in {\cal C}_{{\rm GL}_{n,k},G}\left( \rho
\right) $ an element $u_{\rho ,\sigma }\in {\rm GL}_{n,k}\left( k\right) $
such that $u_{\rho ,\sigma }\sigma u_{\rho ,\sigma }^{-1}=\rho $ (in the
obvious functor-theoretic sense), determines a unique isomorphism of
algebraic spaces over $k$%
\[
j_\rho :\coprod_{\sigma \in {\cal C}_{{\rm GL}_{n,k},G}(\rho )}{\rm N}_{{\rm %
GL}_{n,k}}\left( \sigma \right) \times ^{{\rm N}_G\left( \sigma \right)
}X^\sigma \longrightarrow Y^\rho ;
\]

(ii) ${\cal C}_{{\rm GL}_{n,k},G}\left( \rho \right) $ is finite.
\end{proposition}

\TeXButton{Proof}{\proof} (ii) follows from (i) since $Y^\rho $ is
quasi-compact. The proof of (i) requires several steps.

(a) {\em Definition of }$j_\rho ${\em .} If $\sigma \in {\cal C}_{{\rm GL}%
_{n,k},G}\left( \rho \right) $, let ${\cal N}_\sigma $ be the presheaf on
the category {\tt Sch}$_{/k}$ of $k$-schemes which associates to $%
T\rightarrow {\rm Spec}k$ the set 
\[
{\cal N}_\sigma \left( T\right) \doteq \frac{{\rm N}_{{\rm GL}_{n,k}}\left(
\sigma \right) \left( T\right) \times X^\sigma \left( T\right) }{{\rm N}%
_G\left( \sigma \right) \left( T\right) }; 
\]
since ${\rm N}_G\left( \sigma \right) $ acts freely on ${\rm N}_{{\rm GL}%
_{n,k}}\left( \sigma \right) \times X^\sigma $ (on the left), the flat sheaf
associated to ${\cal N}_\sigma $ is ${\rm N}_{{\rm GL}_{n,k}}\left( \sigma
\right) \times ^{{\rm N}_G\left( \sigma \right) }X^\sigma $. Let $\widehat{Y}%
^\rho $ be the presheaf on {\tt Sch}$_{/k}$ which associates to $%
T\rightarrow {\rm Spec}k$ the set 
\[
\widehat{Y}^\rho \left( T\right) \doteq \left\{ \left[ A,x\right] \in \frac{%
{\rm GL}_{n,k}\left( T\right) \times X\left( T\right) }{G\left( T\right) }%
\mid \forall T^{\prime }\rightarrow T\text{, }\forall r\in \rho \left(
T^{\prime }\right) \text{, }\left[ rA_{T^{\prime }},x_{T^{\prime }}\right]
=\left[ A_{T^{\prime }},x_{T^{\prime }}\right] \right\} ; 
\]
the flat sheaf associated to $\widehat{Y}^\rho $ is $Y^\rho $ (e.g. {\bf 
\cite{DG}}, II, \S 1, n. 3). If $u_{\rho ,\sigma }\in {\rm GL}_{n,k}\left(
k\right) $ is such that $u_{\rho ,\sigma }\sigma u_{\rho ,\sigma }^{-1}=\rho 
$ (in the obvious functor-theoretic sense), the presheaf map 
\begin{eqnarray*}
\widehat{j}_{\rho ,\sigma } &:&{\cal N}_\sigma \longrightarrow \widehat{Y}%
^\rho \\
\widehat{j}_{\rho ,\sigma }\left( T\right) &:&{\cal N}_\sigma \left(
T\right) \ni \left[ B,x\right] \longrightarrow \left[ u_{\rho ,\sigma
}B,x\right] \in \widehat{Y}^\rho \left( T\right)
\end{eqnarray*}
is easily checked to be well-defined. Let $j_{\rho ,\sigma }:{\rm N}_{{\rm GL%
}_{n,k}}\left( \sigma \right) \times ^{{\rm N}_G\left( \sigma \right)
}X^\sigma \rightarrow Y^\rho $ denote the associated sheaf map and define $%
j_\rho \doteq \coprod_{\sigma \in {\cal C}_{{\rm GL}_{n,k},G}\left( \rho
\right) }j_{\rho ,\sigma }$.

(b) $j_\rho ${\em \ induces a bijection on geometric points.}{\bf \ }This is
an elementary check. Let $\xi \in Y^\rho \left( \Omega \right) $ be a
geometric point. Then there exists a {\em fppf }cover $T_0\rightarrow {\rm %
Spec}\Omega $ and an element $\left[ A,x\right] \in \widehat{Y}$ $^\rho
\left( T_0\right) $ representing $\xi $. Therefore, for each $T\rightarrow
T_0$ and each $r\in \rho \left( T\right) $ there exists $g\in G\left(
T\right) $ such that: 
\begin{eqnarray*}
rA_Tg^{-1} &=&A_T \\
gx_T &=&x_T.
\end{eqnarray*}
Then, $A^{-1}\rho A$ defines (functorially over $T_0$) a dual cyclic
subgroup scheme $\sigma _0^{\prime }$ of $G_{\left( T_0\right) }$ over $T_0$%
. Since $\sigma _0^{\prime }$ is isomorphic to some $\mu _{n,T_0}$, it
descends to a dual cyclic subgroup $\sigma ^{\prime }$ of $G$ over $k$ which
is {\rm GL}$_{n,k}$-conjugate to $\rho $ since $T_0\rightarrow {\rm Spec}%
\Omega $ has a section and ${\rm GL}_{n,k}$ satisfies our rationality
condition $({\rm RC})$ (see Remark \ref{rmkRC} (i)). By definition of ${\cal %
C}_{{\rm GL}_{n,k},G}\left( \rho \right) $ there exists a unique $\sigma \in 
{\cal C}_{{\rm GL}_{n,k},G}\left( \rho \right) $ which is $G$-conjugated to $%
\sigma ^{\prime }$ over $k$ i.e. 
\[
g\sigma ^{\prime }g^{-1}=\sigma 
\]
(functorially) for some $g\in G\left( k\right) $. Since $\sigma \in {\cal C}%
_{{\rm GL}_{n,k},G}\left( \rho \right) $, there is an element $u\in
GL_{n,k}\left( k\right) $ such that $u\sigma u^{-1}=\rho $. Therefore $%
u^{-1}Ag^{-1}$ restricted to $T_0$ is in ${\rm N}_{{\rm GL}_{n,k}}\left(
\sigma \right) \left( T_0\right) ,$ $gx\in X^\sigma \left( T_0\right) $ and
if $\left[ u^{-1}Ag^{-1},gx\right] ^{\sim }$ denotes the element in $\left( 
{\rm N}_{{\rm GL}_{n,k}}\left( \sigma \right) \times ^{{\rm N}_G\left(
\sigma \right) }X^\sigma \right) \left( \Omega \right) $ represented by the
element $\left[ u^{-1}Ag^{-1},gx\right] $ in ${\cal N}_\sigma \left(
T_0\right) $, we have $j_{\rho ,\sigma }\left( \Omega \right) \left( \left[
u^{-1}Ag^{-1},gx\right] ^{\sim }\right) =\xi $ by definition of $j_{\rho
,\sigma }$. Thus $j_\rho \left( \Omega \right) $ is surjective.

Now, let $\eta \in $ $\left( {\rm N}_{{\rm GL}_{n,k}}\left( \sigma \right)
\times ^{{\rm N}_G\left( \sigma \right) }X^\sigma \right) \left( \Omega
\right) $ (respectively, $\eta ^{\prime }$ $\in \left( {\rm N}_{{\rm GL}%
_{n,k}}\left( \sigma ^{\prime }\right) \times ^{{\rm N}_G\left( \sigma
^{\prime }\right) }X^{\sigma ^{\prime }}\right) \left( \Omega \right) $) for 
$\sigma $ and $\sigma ^{\prime }$ in ${\cal C}_{{\rm GL}_{n,k},G}\left( \rho
\right) $. Choosing a common refinement, we can assume there exists a {\em %
fppf} cover $T_0\rightarrow {\rm Spec}\Omega $ such that $\eta $ (resp. $%
\eta ^{\prime }$) is represented by an element $\left[ B,y\right] \in $ $%
{\cal N}_\sigma \left( T_0\right) $ (resp. $\left[ B^{\prime },y^{\prime
}\right] \in {\cal N}_{\sigma ^{\prime }}\left( T_0\right) $). If $j_\rho
\left( \Omega \right) \left( \eta \right) =j_\rho \left( \Omega \right)
\left( \eta ^{\prime }\right) $, there exists a {\em fppf} cover $%
T_1\rightarrow T_0$ such that $\left[ u_{\rho ,\sigma }B,y\right] =\left[
u_{\rho ,\sigma ^{\prime }}B^{\prime },y^{\prime }\right] $ in ${\rm GL}%
_{n,k}\left( T_1\right) \times X\left( T_1\right) /G\left( T_1\right) $ i.e.
there is an element $g\in G\left( T_1\right) $ such that 
\begin{eqnarray*}
u_{\rho ,\sigma }Bg^{-1} &=&u_{\rho ,\sigma ^{\prime }}B^{\prime }\text{
\quad in }{\rm GL}_{n,k}\left( T_1\right) \\
gy &=&y^{\prime }\quad \text{ in }X\left( T_1\right) .
\end{eqnarray*}
Then it is easy to check that $\sigma =g^{-1}\sigma ^{\prime }g$ over $T_1$
and, as in the proof of surjectivity of $j_\rho \left( \Omega \right) $,
since $T_1\rightarrow {\rm Spec}\Omega $ has a section and $G$ satisfies our
rationality condition $({\rm RC})$, $\sigma $ and $\sigma ^{\prime }$ are $G$%
-conjugated over $k$ as well and therefore $\sigma =\sigma ^{\prime }$ as
elements in ${\cal C}_{{\rm GL}_{n,k},G}\left( \rho \right) $. In
particular, $g\in {\rm N}_G\left( \sigma \right) \left( T_1\right) $ and $%
\left[ B,y\right] =\left[ B^{\prime },y^{\prime }\right] $ in ${\cal N}%
_\sigma \left( T_1\right) $. Since $T_1\rightarrow {\rm Spec}\Omega $ is
still a {\em fppf} cover, we have $\eta =\eta ^{\prime }$ and $j_\rho \left(
\Omega \right) $ is injective.

(c) {\em Each }$j_{\rho ,\sigma }${\em \ is a closed and open immersion.} It
is enough to show that each $j_{\rho ,\sigma }$ is an open immersion because
in this case it is also a closed immersion, $Y^\rho $ being quasi-compact.
Since ${\rm N}_{{\rm GL}_{n,k}}\left( \rho \right) $ acts on both $%
\coprod_{\sigma \in {\cal C}_{{\rm GL}_{n,k},G}\left( \rho \right) }{\rm N}_{%
{\rm GL}_{n,k}}\left( \sigma \right) \times ^{{\rm N}_G\left( \sigma \right)
}X^\sigma $ and $Y^\rho $ and $j_\rho $ is equivariant, it will be enough to
prove that $j_{\rho ,\rho }$ is an open immersion. We will prove first that $%
j_{\rho ,\rho }$ is injective and unramified and then conclude the proof by
showing that it is also flat (in fact, an \'etale injective map is an open
immersion).

(c$_1$) $j_{\rho ,\rho }${\em \ is injective and unramified. }It is enough
to show that the inverse image under $j_{\rho ,\rho }$ of a geometric point
is a (geometric) point. Consider the commutative diagram 
\[
\begin{tabular}{ccc}
${\rm N}_{{\rm GL}_{n,k}}\left( \rho \right) \times X^\rho $ & $\stackrel{l}{%
\longrightarrow }$ & ${\rm GL}_{n,k}\times X$ \\ 
$^p\downarrow $ &  & $\downarrow ^\pi $ \\ 
${\rm N}_{{\rm GL}_{n,k}}\left( \rho \right) \times ^{{\rm N}_G\left( \rho
\right) }X^\rho $ & $\stackunder{i_\rho \circ j_{\rho ,\rho }}{%
\longrightarrow }$ & $Y$%
\end{tabular}
\]
where $l$ and $i_\rho :Y^\rho \hookrightarrow Y$ are the natural inclusions
and $p,\pi $ the natural projections. Let $y_0$ be a geometric point of $Y$
in the image of $i_\rho \circ j_{\rho ,\rho }$; using the action of ${\rm N}%
_{{\rm GL}_{n,k}}\left( \rho \right) $ on ${\rm N}_{{\rm GL}_{n,k}}\left(
\rho \right) \times ^{{\rm N}_G\left( \rho \right) }X^\rho $ and $Y^\rho $,
we may suppose that $y_0$ is of the form $\left[ 1,x_0\right] \in Y^\rho
\left( \Omega \right) $, with $\Omega $ an algebraically closed field and $%
x_0\in X^\rho \left( \Omega \right) $. Obviously, $\overline{\left(
1,x_0\right) }\in {\rm N}_{{\rm GL}_{n,k}}\left( \rho \right) \times ^{{\rm N%
}_G\left( \rho \right) }X^\rho \left( \Omega \right) $ is contained in $%
j_{\rho ,\rho }^{-1}\left( y_0\right) $ and, by faithful flatness of $p$, $%
j_{\rho ,\rho }^{-1}\left( y_0\right) =\overline{\left( 1,x_0\right) }$ if 
\begin{equation}
p^{-1}\left( \overline{\left( 1,x_0\right) }\right) =\pi ^{-1}\left(
y_0\right) \cap \left( {\rm N}_{{\rm GL}_{n,k}}\left( \rho \right) \times ^{%
{\rm N}_G\left( \rho \right) }X^\rho \right) .  \label{intersection}
\end{equation}
But $G\left( \Omega \right) \simeq \pi ^{-1}\left( y_0\right) $ via $%
g\mapsto \left( g^{-1},gx_0\right) $ and ${\rm N}_G\left( \rho \right)
\left( \Omega \right) \simeq p^{-1}\left( \overline{\left( 1,x_0\right) }%
\right) $ via $h\mapsto \left( h^{-1},hx_0\right) $, therefore (\ref
{intersection}) follows from ${\rm N}_{{\rm G}}\left( \rho \right) ={\rm N}_{%
{\rm GL}_{n,k}}\left( \rho \right) \cap G$.

(c$_2$) $j_{\rho ,\rho }${\em \ is flat. }This fact is proved in the next
subsection where we also single out a more general technical hypothesis for
the action of $G$ on $X$ under which Proposition \ref{disjointunion} still
holds.

\TeXButton{End Proof}{\endproof}

The remaining part of this subsection will be devoted to conclude the proof
of Theorem \ref{main} using Proposition \ref{disjointunion}. First we show
that Proposition \ref{disjointunion} (ii) allows one to define a canonical
ismorphism 
\[
\prod_{\rho \in {\cal C}({\rm GL}_{n,k})}\left( K_{*}\left( Y^\rho ,{\rm C}_{%
{\rm GL}_{n,k}}\left( \rho \right) \right) _{{\rm geom}}\otimes \widetilde{%
{\rm \QTR{mathrm}{R}}}\left( \rho \right) _\Lambda \right) ^{w_{{\rm GL}%
_{n,k}}(\rho )}\simeq 
\]
\[
\prod_{\sigma \in {\cal C}({\rm G})}\left( K_{*}\left( {\rm N}_{{\rm GL}%
_{n,k}}\left( \sigma \right) \times ^{{\rm N}_{{\rm G}}\left( \sigma \right)
}X^\sigma ,{\rm C}_{{\rm GL}_{n,k}}\left( \sigma \right) \right) _{{\rm geom}%
}\otimes \widetilde{{\rm \QTR{mathrm}{R}}}\left( \sigma \right) _\Lambda
\right) ^{w_{{\rm GL}_{n,k}}\left( \sigma \right) }\text{ ;}
\]
next we show, using Lemma \ref{smalllemma}, that each factor in the r.h.s.
is isomorphic to 
\[
(K_{*}\left( {\rm C}_{{\rm GL}_{n,k}}\left( \sigma \right) \times ^{{\rm C}_{%
{\rm G}}\left( \sigma \right) }X^\sigma ,{\rm C}_{{\rm GL}_{n,k}}\left(
\sigma \right) \right) _{{\rm geom}}\otimes \widetilde{{\rm \QTR{mathrm}{R}}}%
\left( \sigma \right) _\Lambda ^{})^{w_{{\rm G}}\left( \sigma \right) }.
\]
The conclusion (i.e. the isomorphism (\ref{phantom})) is then accomplished
by establishing, for any regular noetherian separated algebraic space $Z$ on
which $G$ acts with finite stabilizers, a ''geometric'' Morita equivalence 
\[
K_{*}\left( {\rm GL}_{n,k}\times ^GZ,{\rm GL}_{n,k}\right) _{{\rm geom}%
}\simeq K_{*}\left( Z,G\right) _{{\rm geom}}\text{ .}\label{kinky}
\]

\smallskip\ 

First of all, note that the choice of a family $\left\{ u_{\rho ,\sigma
}\mid \sigma \in {\cal C}_{{\rm GL}_{n,k},G}\left( \rho \right) \right\} $
of elements $u_{\rho ,\sigma }\in {\rm GL}_{n,k}\left( k\right) $ such that $%
u_{\rho ,\sigma }\sigma u_{\rho ,\sigma }^{-1}=\rho $, which uniquely
defines $j_\rho $ in Proposition \ref{disjointunion}, also determines a
unique family of isomorphisms 
\[
\left\{ {\rm int}\left( u_{\rho ,\sigma }\right) :{\rm \QTR{mathrm}{C}}_{%
{\rm GL}_{n,k}}\left( \rho \right) \rightarrow {\rm \QTR{mathrm}{C}}_{{\rm GL%
}_{n,k}}\left( \sigma \right) \mid \sigma \in {\cal C}_{{\rm GL}%
_{n,k},G}\left( \rho \right) \right\} 
\]
(where ${\rm int}(u_{\rho ,\sigma })$ denotes conjugation by $u_{\rho
,\sigma }$) and this family gives us an action of ${\rm C}_{{\rm GL}%
_{n,k}}\left( \rho \right) $ on 
\[
\coprod_{\sigma \in {\cal C}_{{\rm GL}_{n,k},G}(\rho )}{\rm N}_{{\rm GL}%
_{n,k}}\left( \sigma \right) \times ^{{\rm N}_G\left( \sigma \right)
}X^\sigma 
\]
(since ${\rm N}_{{\rm GL}_{n,k}}\left( \sigma \right) $, and then ${\rm C}_{%
{\rm GL}_{n,k}}\left( \sigma \right) $, acts naturally on ${\rm N}_{{\rm GL}%
_{n,k}}\left( \sigma \right) \times ^{{\rm N}_G\left( \sigma \right)
}X^\sigma $ by left multiplication on ${\rm N}_{{\rm GL}_{n,k}}\left( \sigma
\right) $). With this action, $j_\rho $ becomes a ${\rm C}_{{\rm GL}%
_{n,k}}\left( \rho \right) $-equivariant isomorphism and since ${\rm int}%
\left( u_{\rho ,\sigma }\right) $ induces an isomorphism ${\rm \QTR{mathrm}{R%
}}({\rm C}_{{\rm GL}_{n,k}}\left( \rho \right) )\simeq {\rm \QTR{mathrm}{R}}(%
{\rm C}_{{\rm GL}_{n,k}}\left( \sigma \right) )$ commuting with rank
morphisms, $j_\sigma $ induces an isomorphism 
\[
K_{*}\left( Y^\rho ,{\rm C}_{{\rm GL}_{n,k}}\left( \rho \right) \right) _{%
{\rm geom}}\otimes \widetilde{{\rm \QTR{mathrm}{R}}}\left( \rho \right)
_\Lambda \simeq \prod_{\sigma \in {\cal C}_{{\rm GL}_{n,k},G}\left( \rho
\right) }K_{*}\left( {\rm N}_{{\rm GL}_{n,k}}\left( \sigma \right) \times ^{%
{\rm N}_{{\rm G}}\left( \sigma \right) }X^\sigma ,{\rm C}_{{\rm GL}%
_{n,k}}\left( \sigma \right) \right) _{{\rm geom}}\otimes \widetilde{{\rm 
\QTR{mathrm}{R}}}\left( \sigma \right) _\Lambda 
\]
which, by definition of the action of ${\rm N}_{{\rm GL}_{n,k}}\left( \rho
\right) $ on each ${\rm N}_{{\rm GL}_{n,k}}\left( \sigma \right) \times ^{%
{\rm N}_G\left( \sigma \right) }X^\sigma $, induces an isomorphism 
\begin{equation}
\left( K_{*}\left( Y^\rho ,{\rm C}_{{\rm GL}_{n,k}}\left( \rho \right)
\right) _{{\rm geom}}\otimes \widetilde{{\rm \QTR{mathrm}{R}}}\left( \rho
\right) _\Lambda \right) ^{w_{{\rm GL}_{n,k}}\left( \rho \right) }\simeq
\label{almostend}
\end{equation}
\[
\prod_{\sigma \in {\cal C}_{{\rm GL}_{n,k},G}\left( \rho \right) }\left(
K_{*}\left( {\rm N}_{{\rm GL}_{n,k}}\left( \sigma \right) \times ^{{\rm N}_{%
{\rm G}}\left( \sigma \right) }X^\sigma ,{\rm C}_{{\rm GL}_{n,k}}\left(
\sigma \right) \right) _{{\rm geom}}\otimes \widetilde{{\rm \QTR{mathrm}{R}}}%
\left( \sigma \right) _\Lambda \right) ^{w_{{\rm GL}_{n,k}}\left( \sigma
\right) } 
\]
Now, if $j_\rho ^{\prime }$ is induced, as in Proposition \ref{disjointunion}%
, by another choice of a family $\left\{ v_{\rho ,\sigma }\mid \sigma \in 
{\cal C}_{{\rm GL}_{n,k},G}\left( \rho \right) \right\} $ of elements $%
v_{\rho ,\sigma }\in {\rm GL}_{n,k}\left( k\right) $ such that $v_{\rho
,\sigma }\sigma v_{\rho ,\sigma }^{-1}=\rho $, $v_{\rho ,\sigma
}^{-1}u_{\rho ,\sigma }\in {\rm N}_{{\rm GL}_{n,k}}\left( \sigma \right)
\left( k\right) $ and there is a commutative diagram 
\[
\begin{tabular}{ccc}
${\rm N}_{{\rm GL}_{n,k}}\left( \sigma \right) \times ^{{\rm N}_{{\rm G}%
}\left( \sigma \right) }X^\sigma $ & $\stackrel{\left( v_{\rho ,\sigma
}^{-1}u_{\rho ,\sigma }\right) \cdot }{\longrightarrow }$ & ${\rm N}_{{\rm GL%
}_{n,k}}\left( \sigma \right) \times ^{{\rm N}_{{\rm G}}\left( \sigma
\right) }X^\sigma $ \\ 
$\qquad \qquad \qquad \quad _{j_\rho }\searrow $ &  & $\swarrow _{j_\rho
^{\prime }}\qquad \qquad \qquad $ \\ 
& $Y^\rho $ & 
\end{tabular}
\]
Therefore, the isomorphism (\ref{almostend}) on the invariants is actually 
{\em independent }on the choice of the family $\left\{ u_{\rho ,\sigma }\mid
\sigma \in {\cal C}_{{\rm GL}_{n,k},G}\left( \rho \right) \right\} $. Since $%
{\cal C}_{{\rm GL}_{n,k},G}\left( \rho \right) =\alpha ^{-1}\left( \rho
\right) $ and, as already observed, $Y^\rho =\emptyset $ unless $\rho \in 
{\rm im}\left( \alpha \right) $, this gives us a canonical isomorphism 
\[
\prod_{\rho \in {\cal C}({\rm GL}_{n,k})}\left( K_{*}\left( Y^\rho ,{\rm C}_{%
{\rm GL}_{n,k}}\left( \rho \right) \right) _{{\rm geom}}\otimes \widetilde{%
{\rm \QTR{mathrm}{R}}}\left( \rho \right) _\Lambda \right) ^{w_{{\rm GL}%
_{n,k}}(\rho )}\simeq 
\]
\[
\prod_{\sigma \in {\cal C}({\rm G})}\left( K_{*}\left( {\rm N}_{{\rm GL}%
_{n,k}}\left( \sigma \right) \times ^{{\rm N}_{{\rm G}}\left( \sigma \right)
}X^\sigma ,{\rm C}_{{\rm GL}_{n,k}}\left( \sigma \right) \right) _{{\rm geom}%
}\otimes \widetilde{{\rm \QTR{mathrm}{R}}}\left( \sigma \right) _\Lambda
\right) ^{w_{{\rm GL}_{n,k}}\left( \sigma \right) } 
\]

Now, let us fix $\sigma \in {\cal C}({\rm G})$ and let us choose a set $%
{\cal A}\subset {\rm N}_{{\rm GL}_{n,k}}\left( \sigma \right) \left(
k\right) $ such that the classes in $w_{{\rm GL}_{n,k}}\left( \sigma \right) 
$ of the elements in ${\cal A}$ constitute a set of representatives for the $%
w_{{\rm G}}\left( \sigma \right) $-orbits in $w_{{\rm GL}_{n,k}}\left(
\sigma \right) $; ${\cal A}$ is a finite set. Since 
\[
{\rm C}_{{\rm GL}_{n,k}}\left( \sigma \right) \times ^{{\rm C}_{{\rm G}%
}\left( \sigma \right) }X^\sigma \hookrightarrow {\rm N}_{{\rm GL}%
_{n,k}}\left( \sigma \right) \times ^{{\rm N}_{{\rm G}}\left( \sigma \right)
}X^\sigma 
\]
is an open and closed immersion, the morphism 
\begin{eqnarray*}
\coprod_{{\cal A}}{\rm C}_{{\rm GL}_{n,k}}\left( \sigma \right) \times ^{%
{\rm C}_{{\rm G}}\left( \sigma \right) }X^\sigma &\longrightarrow &{\rm N}_{%
{\rm GL}_{n,k}}\left( \sigma \right) \times ^{{\rm N}_{{\rm G}}\left( \sigma
\right) }X^\sigma \\
\left[ C,x\right] _{A_i\in {\cal A}} &\longmapsto &\left[ A_i\cdot C,x\right]
\end{eqnarray*}
(in the obvious functor-theoretic sense), which is easily checked to induce
an isomorphism on geometric points, is an isomorphism. Therefore, there is
an induced isomorphism 
\begin{eqnarray*}
&&\ \ \ \ \prod_{{\cal A}}K_{*}\left( {\rm C}_{{\rm GL}_{n,k}}\left( \sigma
\right) \times ^{{\rm C}_{{\rm G}}\left( \sigma \right) }X^\sigma ,{\rm C}_{%
{\rm GL}_{n,k}}\left( \sigma \right) \right) _{{\rm geom}}\otimes \widetilde{%
{\rm \QTR{mathrm}{R}}}\left( \sigma \right) _\Lambda \\
\ &\simeq &K_{*}\left( {\rm N}_{{\rm GL}_{n,k}}\left( \sigma \right) \times
^{{\rm N}_{{\rm G}}\left( \sigma \right) }X^\sigma ,{\rm C}_{{\rm GL}%
_{n,k}}\left( \sigma \right) \right) _{{\rm geom}}\otimes \widetilde{{\rm 
\QTR{mathrm}{R}}}\left( \sigma \right) _\Lambda
\end{eqnarray*}
Since $w_{{\rm GL}_{n,k}}\left( \sigma \right) $ acts transitively on ${\cal %
A}$ with stabilizer $w_{{\rm G}}\left( \sigma \right) $, by Lemma \ref
{smalllemma}, we get a canonical isomorphism; 
\[
\left( K_{*}\left( {\rm N}_{{\rm GL}_{n,k}}\left( \sigma \right) \times ^{%
{\rm N}_{{\rm G}}\left( \sigma \right) }X^\sigma ,{\rm C}_{{\rm GL}%
_{n,k}}\left( \sigma \right) \right) _{{\rm geom}}\otimes \widetilde{{\rm 
\QTR{mathrm}{R}}}\left( \sigma \right) _\Lambda \right) ^{w_{{\rm GL}%
_{n,k}}\left( \sigma \right) } 
\]
\[
\simeq \left( K_{*}\left( {\rm C}_{{\rm GL}_{n,k}}\left( \sigma \right)
\times ^{{\rm C}_{{\rm G}}\left( \sigma \right) }X^\sigma ,{\rm C}_{{\rm GL}%
_{n,k}}\left( \sigma \right) \right) _{{\rm geom}}\otimes \widetilde{{\rm 
\QTR{mathrm}{R}}}\left( \sigma \right) _\Lambda \right) ^{w_{{\rm G}}\left(
\sigma \right) }. 
\]
Since by Morita equivalence ({\bf \cite{Th3}}, Proposition 6.2), 
\begin{equation}
K_{*}\left( {\rm C}_{{\rm GL}_{n,k}}\left( \sigma \right) \times ^{{\rm C}_{%
{\rm G}}\left( \sigma \right) }X^\sigma ,{\rm C}_{{\rm GL}_{n,k}}\left(
\sigma \right) \right) \simeq K_{*}\left( X^\sigma ,{\rm C}_{{\rm G}}\left(
\sigma \right) \right) ,  \label{mmorita}
\end{equation}
to conclude the proof of Theorem \ref{main}, we need only to show that the
natural morphism 
\begin{equation}
K_{*}\left( {\rm C}_{{\rm GL}_{n,k}}\left( \sigma \right) \times ^{{\rm C}_{%
{\rm G}}\left( \sigma \right) }X^\sigma ,{\rm C}_{{\rm GL}_{n,k}}\left(
\sigma \right) \right) _{{\rm geom}}\simeq K_{*}\left( X^\sigma ,{\rm C}_{%
{\rm G}}\left( \sigma \right) \right) _{{\rm geom}}  \label{moritageom}
\end{equation}
induced by (\ref{mmorita}) is still an isomorphism. Since the diagram 
\[
\begin{tabular}{lll}
$K_{*}\left( {\rm GL}_{n,k}\times ^{{\rm C}_{{\rm G}}\left( \sigma \right)
}X^\sigma ,{\rm GL}_{n,k}\right) _{{\rm geom}}$ & $\stackrel{\alpha }{%
\longrightarrow }$ & $K_{*}\left( {\rm C}_{{\rm GL}_{n,k}}\left( \sigma
\right) \times ^{{\rm C}_{{\rm G}}\left( \sigma \right) }X^\sigma ,{\rm C}_{%
{\rm GL}_{n,k}}\left( \sigma \right) \right) _{{\rm geom}}$ \\ 
& $_\gamma \searrow $ & $\qquad \qquad \downarrow ^\beta $ \\ 
&  & $K_{*}\left( X^\sigma ,{\rm C}_{{\rm G}}\left( \sigma \right) \right) _{%
{\rm geom}}$%
\end{tabular}
\]
is commutative and, by Morita equivalence 
\begin{eqnarray*}
K_{*}\left( {\rm C}_{{\rm GL}_{n,k}}\left( \sigma \right) \times ^{{\rm C}_{%
{\rm G}}\left( \sigma \right) }X^\sigma ,{\rm C}_{{\rm GL}_{n,k}}\left(
\sigma \right) \right) &\simeq &K_{*}\left( {\rm GL}_{n,k}\times ^{{\rm C}_{%
{\rm GL}_{n,k}}\left( \sigma \right) }{\rm C}_{{\rm GL}_{n,k}}\left( \sigma
\right) \times ^{{\rm C}_{{\rm G}}\left( \sigma \right) }X^\sigma ,{\rm GL}%
_{n,k}\right) \\
&\simeq &\ \ K_{*}\left( {\rm GL}_{n,k}\times ^{{\rm C}_{{\rm GL}%
_{n,k}}\left( \sigma \right) }X^\sigma ,{\rm GL}_{n,k}\right) ,
\end{eqnarray*}
to show $\beta $ is an isomorphism it is enough to prove that for any
regular separated algebraic space $Z$ on which $G$ acts with finite
stabilizers, Morita equivalence induces an isomorphism 
\begin{equation}
K_{*}\left( {\rm GL}_{n,k}\times ^GZ,{\rm GL}_{n,k}\right) _{{\rm geom}%
}\simeq K_{*}\left( Z,G\right) _{{\rm geom}},  \label{kinky}
\end{equation}
since in this case both $\alpha $ and $\gamma $ are isomorphisms.

Let $\pi :{\rm \QTR{mathrm}{R}}\left( {\rm GL}_{n,k}\right) \rightarrow {\rm 
\QTR{mathrm}{R}}\left( G\right) $ is the restriction morphism, $\rho :{\rm 
\QTR{mathrm}{R}}\left( G\right) \rightarrow K_0\left( Z,G\right) $ the
pullback along $Z\rightarrow {\rm Spec}k$, ${\rm rk}^{\prime }:{\rm R}\left( 
{\rm GL}_{n,k}\right) \rightarrow \Lambda $ and ${\rm rk}:{\rm R}\left(
G\right) \rightarrow \Lambda $ the rank morphisms, $S^{\prime }\doteq \left( 
{\rm rk}^{\prime }\right) ^{-1}\left( 1\right) $, $S\doteq \left( {\rm rk}%
\right) ^{-1}\left( 1\right) $ and $T\doteq \pi \left( S^{\prime }\right)
\subseteq S$; the following diagram commutes 
\[
\begin{tabular}{ccc}
$T^{-1}K_0\left( Z,G\right) _\Lambda $ & $\stackrel{{\rm rk}_{0,T}}{%
\rightarrow }$ & $\Lambda $ \\ 
$\downarrow $ & $\quad \nearrow _{{\rm rk}_{{\rm geom}}}$ &  \\ 
$K_0\left( Z,G\right) _{{\rm geom}}$ &  & 
\end{tabular}
\]
where ${\rm rk}_{{\rm geom}}$ and ${\rm rk}_{0,T}$ denote the localizations
of the rank morphism ${\rm rk}_0:K_0\left( Z,G\right) _\Lambda \rightarrow
\Lambda $. By Morita equivalence the natural map (which commutes with the
induced rank morphisms) 
\[
K_0\left( {\rm GL}_{n,k}\times ^{{\rm G}}Z,{\rm GL}_{n,k}\right) _{{\rm geom}%
}\longrightarrow T^{-1}K_0\left( Z,G\right) _\Lambda 
\]
is an isomorphism and then, by Proposition \ref{nilpotent}, $\ker \left( 
{\rm rk}_{0,T}:T^{-1}K_0\left( Z,G\right) _\Lambda \rightarrow \Lambda
\right) $ is nilpotent. Now, if $s\in S$, ${\rm rk}_{0,T}\left( \rho \left(
s\right) /1\right) ={\rm rk}\left( s\right) =1$ and therefore 
\[
T^{-1}K_0\left( Z,G\right) _\Lambda \rightarrow K_0\left( Z,G\right) _{{\rm %
geom}} 
\]
and 
\[
K_0\left( {\rm GL}_{n,k}\times ^{{\rm G}}Z,{\rm GL}_{n,k}\right) _{{\rm geom}%
}\rightarrow K_0\left( Z,G\right) _{{\rm geom}} 
\]
are both isomorphisms. Since $K_{*}\left( {\rm GL}_{n,k}\times ^{{\rm G}}Z,%
{\rm GL}_{n,k}\right) _\Lambda $ is naturally a $K_0\left( {\rm GL}%
_{n,k}\times ^{{\rm G}}Z,{\rm GL}_{n,k}\right) _\Lambda $-module and an $%
{\rm R}\left( {\rm GL}_{n,k}\right) _\Lambda $-module via the pullback ring
morphism 
\[
\rho ^{\prime }:{\rm R}\left( {\rm GL}_{n,k}\right) _\Lambda \rightarrow
K_0\left( {\rm GL}_{n,k}\times ^{{\rm G}}Z,{\rm GL}_{n,k}\right) _\Lambda 
\text{ ,} 
\]
we have: 
\[
K_{*}^{\prime }\left( {\rm GL}_{n,k}\times ^{{\rm G}}Z,{\rm GL}_{n,k}\right)
_{{\rm geom}} 
\]
\[
\simeq K_{*}^{\prime }\left( {\rm GL}_{n,k}\times ^{{\rm G}}Z,{\rm GL}%
_{n,k}\right) _\Lambda \otimes _{K_0\left( {\rm GL}_{n,k}\times ^{{\rm G}}Z,%
{\rm GL}_{n,k}\right) _\Lambda }K_0\left( {\rm GL}_{n,k}\times ^{{\rm G}}Z,%
{\rm GL}_{n,k}\right) _{{\rm geom}} 
\]
\[
\simeq K_{*}^{\prime }\left( Z,G\right) _\Lambda \otimes _{K_0\left(
Z,G\right) _\Lambda }K_0\left( {\rm GL}_{n,k}\times ^{{\rm G}}Z,{\rm GL}%
_{n,k}\right) _{{\rm geom}}\simeq K_{*}^{\prime }\left( Z,G\right) _\Lambda
\otimes _{K_0\left( Z,G\right) _\Lambda }K_0\left( Z,G\right) _{{\rm geom}} 
\]
\[
\simeq K_{*}^{\prime }\left( Z,G\right) _{{\rm geom}} 
\]
which proves (\ref{moritageom}) and conclude the proof of theorem \ref{main}.

\subsection{Hypotheses on $G\label{hypothesesonG}$}

In this subsection we conclude the proof of Proposition \ref{disjointunion},
showing that (this is part (c$_2$) of the proof) 
\[
j_{\rho ,\rho }:{\rm N}_{{\rm GL}_{n,k}}\left( \rho \right) \times ^{{\rm N}%
_G\left( \rho \right) }X^\rho \longrightarrow Y^\rho 
\]
{\em \ }\noindent is flat. This is the only step in the proof of Proposition 
\ref{disjointunion} where we make use of the hypothesis that the quotient
algebraic space $G/{\rm C}_G\left( \rho \right) $ is smooth over $k$.
Actually, our proof will work under the following weaker hypothesis. Let $S$
denote the spectrum of the dual numbers over $k$%
\[
S={\rm Spec}\left( k\left[ \varepsilon \right] \right) 
\]
and for any $k$-group scheme $H$, we let $\overline{{\rm H}}^{{\rm 1}}\left(
S,H\right) $ denote the $k$-vector space of isomorphism classes of pairs $%
\left( P\rightarrow S,y\right) $ where $P\rightarrow S$ is an $H$-torsor and 
$y$ is a $k$-rational point on the closed fiber of $P$. Then Proposition \ref
{disjointunion}, and hence Theorem \ref{main}, still holds with hypothesis

(S) ''for any essential dual cyclic subgroupscheme $\sigma \subseteq G$, the
quotient $G/${\rm C}$_G\left( \sigma \right) $ is smooth''

\smallskip\ 

\noindent over $k$''replaced by the following:

(S') ``for any essential dual cyclic $k$-subgroup scheme $\sigma \subseteq G$%
, we have 
\[
\dim \overline{{\rm H}}^{{\rm 1}}\left( S,{\rm C}_G\left( \sigma \right)
\right) =\dim (\overline{{\rm H}}^{{\rm 1}}\left( S,G\right) )^\sigma \quad 
\text{''.} 
\]

\smallskip\ 

First we will prove that $j_{\rho ,\rho }$ is flat assuming ({\bf S}')
holds. Then we will show that ({\bf S}) implies ({\bf S}'). Since $p:{\rm N}%
_{{\rm GL}_{n,k}}\left( \rho \right) \times X^\rho \rightarrow {\rm N}_{{\rm %
GL}_{n,k}}\left( \rho \right) \times ^{{\rm N}_G\left( \rho \right) }X^\rho $
is faithfully flat, it is enough to prove that $\underline{j}_\rho \doteq
j_{\rho ,\rho }\circ p$ is flat. Let $\pi :{\rm GL}_{n,k}\times X\rightarrow
Y$ be the projection and 
\[
f:{\rm GL}_{n,k}\times X\times G\longrightarrow {\rm GL}_{n,k}\times X 
\]
\[
\left( A,x,g\right) \longmapsto \left( Ag^{-1},gx\right) . 
\]
Consider the following cartesian squares 
\begin{equation}
\begin{tabular}{ccccc}
$U$ & $\stackrel{u_\rho }{\longrightarrow }$ & $\pi ^{-1}\left( Y^\rho
\right) $ & $\hookrightarrow $ & ${\rm GL}_{n,k}\times X$ \\ 
$\downarrow $ &  & $\downarrow $ &  & $\downarrow ^\pi $ \\ 
${\rm N}_{{\rm GL}_{n,k}}\left( \rho \right) \times X^\rho $ & $\stackunder{%
\underline{j}_\rho }{\longrightarrow }$ & $Y^\rho $ & $\hookrightarrow $ & $%
Y $%
\end{tabular}
\label{cart1}
\end{equation}
Since $\pi $ is faithfully flat, it is enough to prove that $u_\rho $ is
flat. But the squares 
\begin{equation}
\begin{tabular}{ccccc}
$U$ & $\hookrightarrow $ & ${\rm GL}_{n,k}\times X\times G$ & $\stackrel{f}{%
\longrightarrow }$ & ${\rm GL}_{n,k}\times X$ \\ 
$\downarrow $ &  & $^{{\rm pr}_{12}}\downarrow $ &  & $\downarrow ^\pi $ \\ 
${\rm N}_{{\rm GL}_{n,k}}\left( \rho \right) \times X^\rho $ & $%
\hookrightarrow $ & ${\rm GL}_{n,k}\times X$ & $\stackunder{\pi }{%
\longrightarrow }$ & $Y$%
\end{tabular}
\label{cart2}
\end{equation}
are cartesian and (in the obvious functor-theoretic sense) 
\[
U=\left\{ \left( A,x,g\right) \in {\rm GL}_{n,k}\times X\times G\mid
A^{-1}\rho A=\rho ,\text{ }x\in X^\rho \right\} \simeq {\rm N}_{{\rm GL}%
_{n,k}}\left( \rho \right) \times X^\rho \times G. 
\]
Moreover, if $P\doteq \left\{ A\in {\rm GL}_{n,k}\mid A^{-1}\rho A\subseteq
G\right\} $, the map 
\[
\pi ^{-1}\left( Y^\rho \right) =\left\{ \left( A,x\right) \in {\rm N}_{{\rm %
GL}_{n,k}}\left( \rho \right) \times X\mid A^{-1}\rho A\subseteq G,\text{ }%
x\in X^{A^{-1}\rho A}\right\} \longrightarrow P\times X^\rho 
\]
\[
\left( A,x\right) \longmapsto \left( A,Ax\right) 
\]
is an isomorphism. Therefore, we are reduced to prove that the map 
\[
v_\rho :{\rm N}_{{\rm GL}_{n,k}}\left( \rho \right) \times X^\rho \times
G\longrightarrow P\times X^\rho 
\]
\[
\left( A,x,g\right) \longmapsto \left( Ag^{-1},Ax\right) 
\]
is flat. But, since the diagram 
\[
\begin{tabular}{ccc}
${\rm N}_{{\rm GL}_{n,k}}\left( \rho \right) \times X^\rho \times G$ & $%
\stackrel{v_\rho }{\longrightarrow }$ & $P\times X^\rho $ \\ 
$^{{\rm pr}_{13}}\downarrow $ &  & $\downarrow ^{{\rm pr}_1}$ \\ 
${\rm N}_{{\rm GL}_{n,k}}\left( \rho \right) \times G$ & $\stackunder{\Theta
_\rho }{\longrightarrow }$ & $P$%
\end{tabular}
\]
where $\Theta _\rho \left( A,g\right) \doteq \left( Ag^{-1}\right) $, is
easily checked to be cartesian, it is enough to show that $\digamma _\rho $
is flat. To do this, let us observe that $\rho $ acts by conjugation on $%
{\rm GL}_{n,k}/G$ (quotient by the $G$-action on ${\rm GL}_{n,k}$ by right
multiplication) and we have a cartesian diagram 
\[
\begin{tabular}{ccc}
$P$ & $\hookrightarrow $ & ${\rm GL}_{n,k}$ \\ 
$^\tau \downarrow $ &  & $\downarrow $ \\ 
$\left( {\rm GL}_{n,k}/G\right) ^\rho $ & $\hookrightarrow $ & ${\rm GL}%
_{n,k}/G$%
\end{tabular}
; 
\]
then $\tau $ is a $G$-torsor and $\Theta _\rho $ is $G$-equivariant. Thus,
the commutative diagram in which the vertical arrows are $G$-torsors 
\[
\begin{tabular}{ccc}
${\rm N}_{{\rm GL}_{n,k}}\left( \rho \right) \times G$ & $\stackrel{\Theta
_\rho }{\longrightarrow }$ & $P$ \\ 
$^{{\rm pr}_1}\downarrow $ &  & $\downarrow ^\tau $ \\ 
${\rm N}_{{\rm GL}_{n,k}}\left( \rho \right) $ & $\stackunder{\chi _\rho }{%
\longrightarrow }$ & $\left( {\rm GL}_{n,k}/G\right) ^\rho $%
\end{tabular}
\]
where $\chi _\rho \left( A\right) \doteq \left[ A\right] \in {\rm GL}%
_{n,k}/G $, is cartesian and we may reduce ourselves to prove that $\chi
_\rho $ is flat. Now observe that ${\rm N}_{{\rm GL}_{n,k}}\left( \rho
\right) $ acts on the left of both ${\rm N}_{{\rm GL}_{n,k}}\left( \rho
\right) $ and $\left( {\rm GL}_{n,k}/G\right) ^\rho $ in such a way that $%
\chi _\rho $ is ${\rm N}_{{\rm GL}_{n,k}}\left( \rho \right) $-equivariant.
Therefore it is enough to prove that $\chi _\rho $ is flat when restricted
to the connected component of the identity in ${\rm N}_{{\rm GL}%
_{n,k}}\left( \rho \right) $ i.e. that the map 
\[
\chi _\rho ^{\prime }:{\rm \QTR{mathrm}{C}}_{{\rm GL}_{n,k}}\left( \rho
\right) \longrightarrow \left( {\rm GL}_{n,k}/G\right) ^\rho 
\]
is flat. Now, ${\rm \QTR{mathrm}{C}}_{{\rm GL}_{n,k}}\left( \rho \right)
=\left( {\rm GL}_{n,k}\right) ^\rho $, where $\rho $ acts by conjugation,
and both $\left( {\rm GL}_{n,k}\right) ^\rho $ and $\left( {\rm GL}%
_{n,k}/G\right) ^\rho $ are smooth by {\bf \cite{Th1}} Prop. 3.1 (since $%
{\rm GL}_{n,k}$ and ${\rm GL}_{n,k}/G$ are smooth) and each fiber of $\chi
_\rho ^{\prime }$ has dimension equal to $\dim \left( {\rm C}_G\left( \rho
\right) \right) $ because $\chi _\rho ^{\prime }$ is ${\rm \QTR{mathrm}{C}}_{%
{\rm GL}_{n,k}}\left( \rho \right) $-equivariant for the natural actions and
all the fibers are obtained from $\left( \chi _\rho ^{\prime }\right)
^{-1}\left( \left[ 1\right] \right) ={\rm \QTR{mathrm}{C}}_G\left( \rho
\right) $ by the ${\rm \QTR{mathrm}{C}}_{{\rm GL}_{n,k}}\left( \rho \right) $%
-action. Therefore, $\chi _\rho ^{\prime }$ is flat if 
\begin{equation}
\dim \left( {\rm \QTR{mathrm}{C}}_{{\rm GL}_{n,k}}\left( \rho \right)
\right) =\dim \left( {\rm \QTR{mathrm}{C}}_G\left( \rho \right) \right)
+\dim \left( \left( {\rm GL}_{n,k}/G\right) ^\rho \right) .  \label{dim}
\end{equation}
Note that, in any case, 
\begin{equation}
\dim \left( {\rm \QTR{mathrm}{C}}_{{\rm GL}_{n,k}}\left( \rho \right)
\right) \leq \dim \left( {\rm \QTR{mathrm}{C}}_G\left( \rho \right) \right)
+\dim \left( \left( {\rm GL}_{n,k}/G\right) ^\rho \right) .  \label{dis}
\end{equation}
Since ${\rm GL}_{n,k}$ is smooth, $\dim (\left( {\rm GL}_{n,k}/G\right)
^\rho )=\dim _k({\rm T}_1\left( {\rm GL}_{n,k}/G\right) ^\rho )$, where $%
{\rm T}_1$ denotes the tangent space at the class of $1\in {\rm GL}_{n,k}$.
Moreover, since $\overline{{\rm H}}^{{\rm 1}}\left( S,{\rm GL}_{n,k}\right)
=0$, there is an exact sequence of $k$-vector spaces 
\begin{equation}
0\rightarrow {\rm Lie}\left( G\right) \longrightarrow {\rm Lie}\left( {\rm GL%
}_{n,k}\right) \longrightarrow {\rm T}_1\left( {\rm GL}_{n,k}/G\right)
\longrightarrow \overline{{\rm H}}^{{\rm 1}}\left( S,G\right) \rightarrow 0
\label{prestar}
\end{equation}
which, $\rho $ being linearly reductive over $k$, yields an exact sequence
of $\rho $-invariants 
\begin{equation}
0\rightarrow {\rm Lie}\left( G\right) ^\rho \longrightarrow {\rm Lie}\left( 
{\rm GL}_{n,k}\right) ^\rho \longrightarrow {\rm T}_1\left( {\rm GL}%
_{n,k}/G\right) ^\rho \longrightarrow \overline{{\rm H}}^{{\rm 1}}\left(
S,G\right) ^\rho \rightarrow 0  \label{star}
\end{equation}
But ${\rm GL}_{n,k}$ is smooth, so 
\[
\dim _k({\rm Lie}\left( {\rm GL}_{n,k}\right) ^\rho )=\dim \left( {\rm GL}%
_{n,k}\right) ^\rho =\dim ({\rm \QTR{mathrm}{C}}_{{\rm GL}_{n,k}}\left( \rho
\right) ) 
\]
and, since ${\rm Lie}\left( {\rm C}_G\left( \rho \right) \right) ={\rm Lie}%
\left( G\right) \cap {\rm Lie}({\rm C}_{{\rm GL}_{n,k}}\left( \rho \right) )$%
, we get 
\[
\dim _k({\rm Lie}\left( G\right) ^\rho )=\dim _k({\rm Lie}\left( {\rm C}%
_G\left( \rho \right) \right) \text{.} 
\]
By (\ref{star}), we get 
\[
\dim _k(\overline{{\rm H}}^{{\rm 1}}\left( S,G\right) ^\rho )=\dim _k({\rm T}%
_1\left( {\rm GL}_{n,k}/G\right) ^\rho )-\dim ({\rm \QTR{mathrm}{C}}_{{\rm GL%
}_{n,k}}\left( \rho \right) )+\dim _k({\rm Lie}\left( {\rm C}_G\left( \rho
\right) \right) = 
\]
\begin{equation}
\dim (\left( {\rm GL}_{n,k}/G\right) ^\rho )-\dim \left( {\rm \QTR{mathrm}{C}%
}_{{\rm GL}_{n,k}}\left( \rho \right) \right) +\dim \left( {\rm \QTR{mathrm}{%
C}}_G\left( \rho \right) \right) +\dim _k({\rm Lie}\left( {\rm C}_G\left(
\rho \right) \right) -\dim ({\rm C}_G\left( \rho \right) )\text{;}
\label{halfbistar}
\end{equation}
hence (\ref{dim}) is satisfied if 
\begin{equation}
\dim _k(\overline{{\rm H}}^{{\rm 1}}\left( S,G\right) ^\rho )=\dim _k({\rm %
Lie}\left( {\rm C}_G\left( \rho \right) \right) -\dim ({\rm C}_G\left( \rho
\right) )\text{.}  \label{bistar}
\end{equation}
But 
\[
\dim _k(\overline{{\rm H}}^{{\rm 1}}\left( S,{\rm C}_G\left( \rho \right)
\right) )=\dim _k({\rm Lie}\left( {\rm C}_G\left( \rho \right) \right) -\dim
({\rm C}_G\left( \rho \right) )\text{,} 
\]
by the exact sequence (analogous to (\ref{prestar}) with $G$ replaced by $%
{\rm C}_G\left( \rho \right) $) 
\[
0\rightarrow {\rm Lie}\left( {\rm C}_G\left( \rho \right) \right)
\longrightarrow {\rm Lie}\left( {\rm GL}_{n,k}\right) \longrightarrow {\rm T}%
_1\left( {\rm GL}_{n,k}/{\rm C}_G\left( \rho \right) \right) \longrightarrow 
\overline{{\rm H}}^{{\rm 1}}\left( S,{\rm C}_G\left( \rho \right) \right)
\rightarrow 0 
\]
hence (\ref{bistar}) holds by hypothesis ({\bf S'}).

We complete the proof of Proposition \ref{disjointunion} showing that ({\bf S%
}) implies ({\bf S'}). Since ${\rm C}_G\left( \rho \right) \subseteq G$, we
have a natural map 
\[
\epsilon :\overline{{\rm H}}^{{\rm 1}}\left( S,{\rm C}_G\left( \rho \right)
\right) \longrightarrow \overline{{\rm H}}^{{\rm 1}}\left( S,G\right) ^\rho 
\text{\quad ;} 
\]
and by (\ref{halfbistar}) and (\ref{dis}), we get 
\begin{equation}
\dim _k(\overline{{\rm H}}^{{\rm 1}}\left( S,G\right) ^\rho )\geq \dim _k(%
\overline{{\rm H}}^{{\rm 1}}\left( S,{\rm C}_G\left( \rho \right) \right) )%
\text{.}  \label{dis2}
\end{equation}
Now, if ({\bf S}) holds, i.e. if $G/{\rm C}_G\left( \rho \right) $ is
smooth, and $\left[ P\rightarrow S,y\right] $ is a class in $\overline{{\rm H%
}}^{{\rm 1}}\left( S,G\right) ^\rho $, $P/{\rm C}_G\left( \rho \right)
\rightarrow S$ is smooth and $y$ induces a point in the closed fiber of $P/%
{\rm C}_G\left( \rho \right) \rightarrow S$; we may reduce the structure
group to ${\rm C}_G\left( \rho \right) $, thus showing that $\epsilon $ is
surjective. By (\ref{dis2}) we conclude that $\epsilon $ is an isomorphism
and this implies ({\bf S'}).

\subsection{Final remarks}

\begin{proposition}
\label{finitegroup}Let $X$ be a noetherian regular separated algebraic space
over $k$ and $G$ a finite group acting on $X$. There is a canonical
isomorphism of $R(G)$-algebras 
\[
K_{*}\left( X,G\right) _{{\rm geom}}\otimes {\Bbb Z}\left[ 1/\left| G\right|
\right] \simeq K_{*}\left( X\right) ^G\otimes {\Bbb Z}\left[ 1/\left|
G\right| \right] .
\]
\end{proposition}

\TeXButton{Proof}{\proof} Since $\ker ({\rm rk}:K_0\left( X\right)
\rightarrow {\rm H}^0(X,{\Bbb Z}\left[ 1/\left| G\right| \right] ))$ is
nilpotent by Corollary \ref{nilpotent0}, the canonical homomorphism 
\[
\pi ^{*}:K_{*}\left( X,G\right) \longrightarrow K_{*}\left( X\right) ^G 
\]
induces a ring homomorphism (still denoted by $\pi ^{*}$): 
\[
\pi ^{*}:K_{*}\left( X,G\right) _{{\rm geom}}\otimes {\Bbb Z}\left[ 1/\left|
G\right| \right] \longrightarrow K_{*}\left( X\right) ^G\otimes {\Bbb Z}%
\left[ 1/\left| G\right| \right] . 
\]
Moreover, the functor 
\[
\pi _{*}:{\cal F}\longmapsto \bigoplus_{g\in G}g^{*}{\cal F} 
\]
defined on coherent ${\cal O}_X$-modules, induces a homomorphism 
\[
\pi _{*}:K_{*}^{\prime }\left( X\right) ^G\otimes {\Bbb Z}\left[ 1/\left|
G\right| \right] \longrightarrow K_{*}^{\prime }\left( X,G\right) _{{\rm geom%
}}\otimes {\Bbb Z}\left[ 1/\left| G\right| \right] 
\]
and (recalling that $K_{*}\left( X,G\right) \simeq K_{*}^{\prime }\left(
X,G\right) $) we obviously get: 
\[
\pi ^{*}\pi _{*}({\cal F)}=\left| G\right| \cdot {\cal F}\text{.} 
\]
On the other hand, we have: 
\[
\pi _{*}\pi ^{*}({\cal F})\simeq {\cal F}\otimes \pi _{*}{\cal O}_X\text{ .} 
\]
But ${\rm rk}(\pi _{*}{\cal O}_X)=\left| G\right| $ and therefore $\pi
_{*}\pi ^{*}$ is an isomorphism too, because of Corollary \ref{nilpotent0}. 
\TeXButton{End Proof}{\endproof}

As a corollary of this result and of Theorem \ref{main}, we recover Theorem
1 of {\bf \cite{Vi2}} which was proved there in a completely different way.

We conclude the paper with a conjecture expressing the fact that $%
K_{*}(X,G)_{{\rm geom}}$ should be the $K$-theory of the quotient $X/G$, if $%
X/G$ is regular, after inverting the orders of all the essential dual cyclic
subgroups of $G$:

\begin{conjecture}
\label{conj}Let $X$ be a noetherian regular separated algebraic space over a
field $k$ and $G$ a linear algebraic $k$-group acting on $X$ with finite
stabilizers in such a way that the quotient $X/G$ exists as a regular
algebraic space. Let $N$ denote the least common multiple of the orders of
all the essential dual cyclic subgroups of $G$ and $\Lambda ={\Bbb Z}\left[
1/N\right] $. If $p:X\rightarrow X/G$ is the quotient map, the composition 
\[
K_{*}\left( X/G\right) _\Lambda \stackrel{p^{*}}{\longrightarrow }%
K_{*}\left( X,G\right) _\Lambda \longrightarrow K_{*}\left( X,G\right) _{%
{\rm geom}}
\]
is an isomorphism.
\end{conjecture}

\begin{remark}
\label{bert}Bertrand Toen pointed out to us that if $X/G$ is smooth it
follows from the results of \cite{EG-GRR} that the composition 
\[
K_0(X/G)\otimes {\Bbb Q}\longrightarrow K_0(X,G)\otimes {\Bbb Q}%
\longrightarrow K_0(X,G)_{{\rm geom}}\otimes {\Bbb Q}
\]
is an isomorphism.
\end{remark}

\section{Appendix: Higher equivariant $K$-theory of noetherian regular
separated algebraic spaces}

In this Appendix we describe the $K$-theories we use in the paper and their
relationships. We essentially follow the example of {\bf \cite{Th-Tr}},
Section 3. Of {\bf \cite{Th-Tr}} we also adopt the language.

Let us remark that it is strongly probable that there exist equivariant
versions of most of the results in {\bf \cite{Th-Tr}}, Section 3. In
particular, there should exist a higher $K$-theory of $G$-equivariant
cohomologically bounded pseudocoherent complexes on $Z$ (respectively, of $G$%
-equivariant perfect complexes on $Z$) for any quasi-compact algebraic space 
$Z$ having most of the alternative models described in {\bf \cite{Th-Tr}},
3.5-3.12. The arguments below can also be considered as a first step toward
an extension of {\bf \cite{Th-Tr}}, 3.11-3.12 to the equivariant case on
algebraic spaces. However, to keep the paper to a reasonable size, we have
decided to give only the results we need and moreover we have made almost no
attempt to optimize the hypotheses.

We would also like to mention the paper {\bf \cite{J}} (in particular
Section 1) in which, among many other results, the general techniques of 
{\bf \cite{Th-Tr}} are used as guidelines for the $K$-theory of arbitrary
Artin stacks.

\smallskip 

We work in a slightly more general situation than required in the rest of
the paper.

Let $S$ be a separated noetherian scheme, $G$ be a group scheme affine over $%
S$ which is finitely presented, separated and flat over $S$. We denote by $G$%
-${\sl AlgSp}_{reg}$ be the category of regular noetherian algebraic spaces
separated over $S$ with an action of $G$ over $S$ and equivariant maps.

\begin{definition}
If $X\in G$-${\sl AlgSp}_{reg}$, we denote by $K_{*}\left( X,G\right) $
(respectively, $K_{*}^{\prime }(X,G)$, resp., $K_{*}^{naive}(X,G)$) the
Waldhausen $K$-theory of the complicial biWaldhausen category ({\bf \cite
{Th-Tr}}) ${\cal W}_{1,X}$ of complexes of quasi-coherent $G$-equivariant $%
{\cal O}_X$-Modules with bounded coherent cohomology (respectively, the
Quillen $K$-theory of the abelian category of $G$-equivariant coherent $%
{\cal O}_X$-Modules, resp., the Quillen $K$-theory of the exact category of $%
G$-equivariant locally free coherent ${\cal O}_X$-Modules).
\end{definition}

\smallskip\ 

\begin{proposition}
\label{a1}Let $Z\to S$ be a morphism of noetherian algebraic spaces such
that the diagonal $Z\to Z\times _SZ$ is affine, $H\to S$ an affine group
space acting on $Z$. Let ${\cal F}$ be an equivariant quasicoherent sheaf on 
$Z$ of finite flat dimension; then there exists a flat equivariant
quasicoherent sheaf ${\cal F}^{\prime }$ on $Z$ together with a surjective $H
$-equivariant homomorphism ${\cal F}^{\prime }\to {\cal F}$.

In particular, if $Z$ is regular this holds for all equivariant
quasicoherent sheaves ${\cal F}$ on $Z$.
\end{proposition}

The hypotheses of the previous Proposition insure that the usual morphism $%
Z\times _SH\to Z\times _SZ$ is affine. In fact, the projection $Z\times
_SH\to Z$ is obviously affine, the projection $Z\times _SZ\to Z$ has affine
diagonal, so this follows from the elementary fact that if $Z\to U\to V$ are
morphisms of algebraic spaces, $Z\to V$ is affine and $U\to V$ has affine
diagonal, then $Z\to U$ is affine. Consider the quotient stack ${\cal Z}%
=[Z/H]$ ({\bf \cite{L-MB}}); the argument above implies that the diagonal $%
{\cal Z}\to {\cal Z}\times _S{\cal Z}$ is affine. Since an $H$-equivariant
quasicoherent ${\cal O}_Z$-Module is the same as a quasicoherent Module over 
${\cal Z}$, now Proposition \ref{a1} follows from the following more general
result:

\begin{proposition}
\label{a2}Let $S$ be a noetherian algebraic space, ${\cal X}$ a noetherian
algebraic stack over $S$ with affine diagonal. Let ${\cal F}$ be a
quasicoherent sheaf of finite flat dimension on ${\cal X}$; then there
exists a flat quasicoherent sheaf ${\cal F}^{\prime }$ on ${\cal X}$
together with a surjective homomorphism ${\cal F}^{\prime }\to {\cal F}$.
\end{proposition}

\TeXButton{Proof}{\proof}Take an affine scheme $U$ with a flat morphism $%
f:U\to {\cal X}$; then $f$ is affine, and in particular the pushforward $%
f_{*}$ on quasicoherent sheaves is exact. Consider a quasicoherent sheaf $%
{\cal F}$ on ${\cal X}$ of finite flat dimension, with the adjunction map $%
{\cal F}\to f_{*}f^{*}{\cal F}$. This map is injective; call ${\cal Q}$ its
cokernel. Clearly the flat dimension of $f_{*}f^{*}{\cal F}$ is the same as
the flat dimension of ${\cal F}$; we claim that the flat dimension of ${\cal %
Q}$ is at most equal to the flat dimension of ${\cal F}$. Now, if there were
a section ${\cal X}\to U$ of $f$, then the sequence 
\[
0\to {\cal F}\to f_{*}f^{*}{\cal F}\to {\cal Q}\to 0 
\]
would split, and this would be clear. However, to compute the flat dimension
of ${\cal Q}$ we can pull back to any flat surjective map to ${\cal X}$; in
particular after pulling back to $U$ we see that $f$ acquires a section, and
the statement is checked. Now $U$ is an affine scheme, so we can take a flat
quasicoherent sheaf ${\cal P}$ on $U$ with a surjective map $u:{\cal P}\to
f^{*}{\cal F}$. Call ${\cal F}^{\prime }$ the kernel of the composition $%
f_{*}{\cal P}\to f_{*}f^{*}{\cal F}\to {\cal Q}$; then ${\cal F}^{\prime }$
surjects onto ${\cal F}$, and fits into an exact sequence $0\to {\cal F}%
^{\prime }\to f_{*}{\cal P}\to {\cal Q}\to 0$. But $f_{*}{\cal P}$ is flat
over ${\cal X}$, so the flat dimension of ${\cal F}^{\prime }$ is less than
the flat dimension of ${\cal Q}$, unless ${\cal Q}$ is flat. But since the
flat dimension of ${\cal Q}$ is at most equal to the flat dimension of $%
{\cal F}$, we see that the flat dimension of ${\cal F}^{\prime }$ is less
than the flat dimension of ${\cal F}$, unless ${\cal F}$ is flat. The proof
is completed with a straightforward induction on the flat dimension of $%
{\cal F}$. \TeXButton{End Proof}{\endproof}

\smallskip\ 

\begin{theorem}
\label{a3} Let $X$ be an object in $G$-${\sl AlgSp}_{reg}$. The obvious
inclusions of the following complicial biWaldhausen categories induce
homotopy equivalences on the Waldhausen $K$-theory spectra $K^{(i)}(X)\doteq
K({\cal W}_{i,X})$, $i=1,2,3$. In particular, the corresponding Waldhausen $K
$-theories $K_{*}^{(i)}(X,G)$ coincide.

{\rm (i)} ${\cal W}_{1,X}=($complexes of quasi-coherent $G$-equivariant $%
{\cal O}_X$-Modules with bounded coherent cohomology$);$

{\rm (ii)} ${\cal W}_{2,X}=\left( \text{bounded complexes in }G-{\tt Coh}%
_X\right) ;$

{\rm (iii)} ${\cal W}_{3,X}=($complexes of flat quasi-coherent $G$%
-equivariant ${\cal O}_X$-Modules with bounded coherent cohomology$)$

Moreover the Waldhausen $K$-theory of any of the categories above coincide
with Quillen $K$-theory $K_{*}^{\prime }(X,G)$ of $G$-equivariant coherent $%
{\cal O}_X$-Modules.
\end{theorem}

\TeXButton{Proof}{\proof}By {\bf \cite{Th6}}, 1.13, the inclusion of ${\cal W%
}_{2,X}$ in ${\cal W}_{1,X}$ induces an equivalence of $K$-theory spectra.
Proposition \ref{a1}, together with {\bf \cite{Th-Tr}}, Lemma 1.9.5 (applied
to ${\cal D}=($ flat $G$-equivariant ${\cal O}_X$-Modules$)$ and ${\cal A}%
=(G $-equivariant ${\cal O}_X$-Modules$)$), implies that for any object $%
E^{*}$ in ${\cal W}_{1,X}$ there exists an object $F^{*}$ in ${\cal W}_{3,X}$
and a quasi-isomorphism $F^{*}\widetilde{\longrightarrow }E^{*}$. Therefore,
by {\bf \cite{Th-Tr}}, 1.9.7 and 1.9.8, the inclusion of ${\cal W}_{3,X}$ in 
${\cal W}_{1,X}$ induces an equivalence of $K$-theory spectra.

The last statement of the Theorem follows immediately from {\bf \cite{Th6}},
1.13, p. 518. \TeXButton{End Proof}{\endproof}

\smallskip\ 

Since any complex in ${\cal W}_{3,X}$ is degreewise flat and $X$ is regular
(hence boudedness of cohomology is preserved under tensor product\footnote{%
In fact this is a non equivariant statement and a local property in the flat
topology, so it reduces to the same statement for regular affine schemes
which is elementary (see also {\bf \cite{SGA6}}).}), the tensor product of
complexes makes the Waldhausen $K$-theory spectrum of ${\cal W}_{3,X}$ into
a functor $K^{(3)}$ from $G$-${\sl AlgSp}_{reg}$ to ring spectra, with
product 
\[
K^{(3)}\wedge K^{(3)}\longrightarrow K^{(3)},
\]
exactly as described in {\bf \cite{Th-Tr}} 3.15. In particular, by Theorem 
\ref{a3}, $K_{*}$ is a functor from $G$-${\sl AlgSp}_{reg}$ to graded rings.
In the same way, tensor product with complexes in ${\cal W}_{3,X}$ gives a
pairing 
\[
K^{(3)}\wedge K^{(1)}\longrightarrow K^{(1)}
\]
between the corresponding functors from $G$-${\sl AlgSp}_{reg}$ to spectra
so that $K_{*}^{(1)}(X,G)$ becomes a module over the ring $K_{*}^{(3)}(X,G)$
functorially in $(X,G)\in G$-${\sl AlgSp}_{reg}$. We denote the
corresponding cap-product by 
\[
\cap :K_{*}^{(3)}(X,G)\otimes K_{*}^{(1)}(X,G)\longrightarrow
K_{*}^{(1)}(X,G)
\]
which becomes the ring product in $K_{*}(X,G)$ with the identifications
allowed by Theorem \ref{a3}. Note that there is an obvious ring morphism $%
\eta :K_{*}^{naive}(X,G)\rightarrow K_{*}^{(3)}(X,G)$ and if 
\[
\cap ^{naive}:K_{*}^{naive}(X,G)\otimes K_{*}^{^{\prime
}}(X,G)\longrightarrow K_{*}^{^{\prime }}(X,G)
\]
denotes the usual ''naive'' cap-product on Quillen $K$-theories, there is a
commutative diagram 
\[
\begin{tabular}{lll}
$K_{*}^{naive}(X,G)\otimes K_{*}^{^{\prime }}(X,G)$ & $\stackrel{\cap
^{naive}}{\longrightarrow }$ & $K_{*}^{^{\prime }}(X,G)$ \\ 
\qquad \qquad $^{\eta \otimes u}\downarrow $ &  & $\qquad \downarrow ^u$ \\ 
$K_{*}^{(3)}(X,G)\otimes K_{*}^{(1)}(X,G)$ & $\stackunder{\cap }{%
\longrightarrow }$ & $K_{*}^{(1)}(X,G)$%
\end{tabular}
\]
where $u$ is the isomorphism of Theorem \ref{a3}. Because of that we will
simply write $\cap $ for both the naive and non-naive cap-products.

Note that, as shown in {\bf \cite{Th6}} 1.13, p.519, $K_{*}^{\prime }(-,G)$
(and therefore $K_{*}(X,G)$ under our hypotheses) is a covariant functor for
proper maps in $G$-${\sl AlgSp}_{reg}$; on the other hand, since any map in $%
G$-${\sl AlgSp}_{reg}$ has finite $Tor$-dimension, $K_{*}(-,G)$ is a
controvariant functor from $G$-${\sl AlgSp}_{reg}$ to (graded) rings. In
fact, if $f:X\longrightarrow Y$ is a morphism in $G$-${\sl AlgSp}_{reg}$ the
same argument in {\bf \cite{Th-Tr}}, 3.14.1 shows that there is an induced
pullback exact functor $f^{*}:{\cal W}_{3,Y}\rightarrow {\cal W}_{3,X}$ and
then we use theorem \ref{a3} to identify $K_{*}^{(3)}(-,G)$ with $K_{*}(-,G)$%
.

\begin{proposition}
\label{a4}{\rm (Projection formula) }Let $j:Z\longrightarrow X$ be a closed
immersion in $G$-${\sl AlgSp}_{reg}$. Then, if $\alpha $ is in $K_{*}(X,G)$
and $\beta $ in $K_{*}^{\prime }\left( Z,G\right) $, we have 
\[
j_{*}(j^{*}(\alpha )\cap \beta )=\alpha \cap j_{*}\left( \beta \right) 
\]
in $K_{*}^{\prime }\left( X,G\right) $.
\end{proposition}

\TeXButton{Proof}{\proof} Since $j$ is affine, $j_{*}$ is exact on
quasi-coherent Modules and therefore induces an exact functor of complicial
biWaldhausen categories $j_{*}:{\cal W}_{1,Z}\rightarrow {\cal W}_{1,X}$
(the condition of bounded coherent cohomology being preserved by regularity
of $Z$ and $X$). Therefore, the maps 
\begin{eqnarray*}
(\alpha ,\beta ) &\longmapsto &j_{*}(j^{*}(\alpha )\cap \beta ),\text{ } \\
(\alpha ,\beta ) &\longmapsto &\text{ }\alpha \cap j_{*}(\beta )
\end{eqnarray*}
from $K_{*}(X,G)\times K_{*}(Z,G)$ to $K_{*}^{\prime }(X,G)\simeq K_{*}(X,G)$
are induced by the exact functors ${\cal W}_{3,X}\times {\cal W}%
_{1,Z}\longrightarrow {\cal W}_{1,X}$%
\begin{eqnarray}
(F^{*},E^{*}) &\longmapsto &j_{*}(j^{*}(F^{*})\otimes E^{*}),\text{ }
\label{bousfield} \\
(F^{*},E^{*}) &\longmapsto &\text{ }F^{*}\otimes j_{*}(E^{*}).  \nonumber
\end{eqnarray}
But for any equivariant quasi-coherent sheaf ${\cal F}$ on $X$ and ${\cal G}$
on $Z$, there is a natural (hence equivariant) isomorphism 
\[
j_{*}(j^{*}{\cal F}\otimes {\cal G})\simeq {\cal F}\otimes j_{*}{\cal G}
\]
which, again by naturality, induces an isomorphism between the two functors
in (\ref{bousfield}); therefore we conclude by {\bf \cite{Th-Tr}}, 1.5.4. 
\TeXButton{End Proof}{\endproof}{\rm \ }

\begin{remark}
Since we only need the projection formula for (regular) closed immersion in
this paper we have decided to state the result only in this case. However,
since by {\bf \cite{Th6}} 1.13 p. 519, $K_{*}(X,G)$ coincides also with
Waldhausen $K$-theory of the category ${\cal W}_{4,X}$ of complexes of $G$%
-equivariant quasi-coherent injective Modules on $X$ with bounded coherent
cohomology, and therefore, by Theorem \ref{a3}, it also coincides with
Waldhausen $K$-theory of the category ${\cal W}_{5,X}$ complexes of $G$%
-equivariant quasi-coherent flasque Modules on $X$ with bounded coherent
cohomology, for any proper map $f:X\rightarrow Y$ in $G$-${\sl AlgSp}_{reg},$%
we have an exact functor $f_{*}:{\cal W}_{5,X}\rightarrow {\cal W}_{5,Y}$
which therefore gives a ''model'' for the push-forward $f_{*}:$ $%
K_{*}(X,G)\rightarrow K_{*}(Y,G)$ (cf. {\bf \cite{Th-Tr}}, 3.16). Now, the
proof of {\bf \cite{Th-Tr}} 3.17 should also give a proof of Proposition \ref
{a4} with $j$ replaced by any proper map in $G$-${\sl AlgSp}_{reg}$, because
it only uses {\bf \cite{Th-Tr}} 2.5.5, which obviously holds for $X$ and $Y$
noetherian algebraic spaces, and {\bf \cite{SGA4}} XVII, 4.2, which should
give a canonical $G$-equivariant Godement flasque resolution of any complex
of $G$-equivariant Modules on any algebraic space in $G$-${\sl AlgSp}_{reg}$%
, since it is developed in a general topos.
\end{remark}

\smallskip\ 

It is very probable that theorem \ref{a3} and therefore the functoriality
with respect to morphisms of finite $Tor$-dimension still hold without the
regularity assumption on the algebraic spaces. On the other hand, it should
also be true that with $G$ and $X$ as above (therefore $X$ regular), the
Waldhausen $K$-theory of the category of $G$-equivariant perfect complexes
on $X$ coincides with $K_{*}^{\prime }(X,G)$. This last statement should
follow, (with a bit of work to identify $K_{*}^{\prime }(X,G)$ with the
Waldhausen $K$-theory of $G$-equivariant pseudocoherent complexes with
bounded cohomology on $X$) from {\bf \cite{J}}, 1.6.2.

\medskip\ 

E-mail addresses of the authors:

Gabriele Vezzosi: {\tt vezzosi@dm.unibo.it}

Angelo Vistoli: {\tt vistoli@dm.unibo.it}


\bigskip 

\bigskip 

\begin{thebibliography}{Th-Tr}
\bibitem[SGA4]{SGA4}  M.\ Artin, A. Grothendieck, J. L. Verdier, {\em %
Th\'eorie des topos et cohomologie \'etale des sch\'emas}, LNM {\bf 305},
Springer-Verlag, New York, 1973.

\bibitem[SGA6]{SGA6}  P. Berthelot, A. Grothendieck, L. Illusie, {\em %
Th\'eorie des intersections et th\'eor\`eme de Riemann-Roch}, LNM {\bf 225},
Springer-Verlag, New York, 1971.\ 

\bibitem[DG]{DG}  M. Demazure, P. Gabriel, {\em Groupes alg\'ebriques}, Tome
I, Masson, Paris, 1970.

\bibitem[SGA3]{SGA3}  M. Demazure, A. Grothendieck, {\em Sch\'emas en groupes%
}, LNM {\bf 151}, {\bf 152}, {\bf 153}, Springer- Verlag, New York, 1970.

\bibitem[EG]{EG-GRR}  D. Edidin, W. Graham, {\em Riemann-Roch for
equivariant Chow groups}, Preprint math.AG/9905081, 13 May 1999.

\bibitem[SGA1]{SGA1}  A. Grothendieck, {\em Rev\^etements \'etales et groupe
fondamental}, LNM {\bf 224}, Springer-Verlag, New York, 1971.

\bibitem[EGAI]{EGAI}  A. Grothendieck, J. A. Dieudonn\'e, {\em El\'ements de
G\'eom\'etrie Alg\'ebrique I}, Springer-Verlag, Berlin, 1971.

\bibitem[J]{J}  R. Joshua, {\em Higher intersection theory on algebraic
stacks II}, Preprint available at
http://www.math.ohio-state.edu/~joshua/pub.html.

\bibitem[L-MB]{L-MB}  G. Laumon, L. Moret-Bailly, {\em Champs Alg\'ebriques}%
, Springer Verlag, 2000.

\bibitem[Me]{Me}  A. S. Merkurjev, {\em Comparison of equivariant and
ordinary }$K${\em -theory of algebraic varieties}, To appear in St.
Petersburg Math. J.

\bibitem[Q]{Q}  D. Quillen, {\em Higher algebraic }$K${\em -theory I}, LNM 
{\bf 341}, Springer-Verlag, New York, 1973.

\bibitem[Se]{Se}  J. P. Serre, {\em Linear representations of finite groups}%
, Springer-Verlag, New York, 1977.

\bibitem[Sr]{Sr}  V. Srinivas, {\em Algebraic }$K${\em -theory}, Second
edition, Birkhauser, Boston, 1996.

\bibitem[Th1]{Th2}  R. W. Thomason, {\em Comparison of equivariant algebraic
and topological }$K${\em -theory}, Duke Math. J. {\bf 53} (1986), 795-825.

\bibitem[Th2]{Th6}  R.W. Thomason, {\em Lefschetz-Riemann-Roch theorem and
coherent trace formula}, Inv. Math. {\bf 85 }(1986), 515-543.

\bibitem[Th3]{Th3}  R. W. Thomason, {\em Algebraic }$K${\em -theory of group
scheme actions}, in ``Algebraic Topology and Algebraic $K$-theory'', Ann. of
Math. Stud. {\bf 113} (1987), 539-563.

\bibitem[Th4]{Th5}  R. W. Thomason, {\em Equivariant algebraic vs.
topological }$K${\em -homology Atiyah-Segal style}, Duke Math. J., 56
(1988), 589-636.

\bibitem[Th5]{Th1}  R. W. Thomason, {\em Une formule de Lefschetz en }$K$%
{\em -th\'eorie \'equivariante alg\'ebrique}, Duke Math. J. {\bf 68} (1992),
447-462.

\bibitem[Th6]{Th4}  R. W. Thomason, {\em Les }$K${\em -groupes d'un sch\'ema
\'eclat\'e et une formule d'intersection exc\'edentaire}, Inv. Math. {\bf 112%
} (1993), 195-215.

\bibitem[Th-Tr]{Th-Tr}  R. W. Thomason, T. Trobaugh, {\em Higher algebraic }$%
K${\em -theory of schemes and of derived categories}, Grothendieck
Festschrift, vol. III, 247-435, Birkhauser.

\bibitem[To1]{To1}  B. Toen, {\em Th\'eor\`emes de Riemann-Roch pour les
champs de Deligne-Mumford}, $K$-Theory {\bf 18} (1999), 33-76.

\bibitem[To2]{To2}  B. Toen, {\em Notes on }$G${\em -theory of
Deligne-Mumford stacks}, Preprint math.AG/9912172.

\bibitem[Vi1]{Vi1}  A. Vistoli, {\em Equivariant Grothendieck groups and
equivariant Chow groups}, in E. Ballico, C. Ciliberto, F. Catanese (eds.),
Classification of irregular varieties, LNM {\bf 1515}, Springer Verlag, New
York, 1992, pp. 112-133.

\bibitem[Vi2]{Vi2}  A. Vistoli, {\em Higher equivariant }$K${\em -theory for
finite group actions}, Duke Math. J., {\bf 63} (1991), 399-419.
\end{thebibliography}
\end{document}